\documentclass[a4paper,11pt,leqno,twoside]{article}
\usepackage{amssymb}
\usepackage[mathscr]{eucal}

\def\NM{{\mathbb{N}}}

\def\MM{{\mathbb{M}}}

\def\HM{{\mathbb{H}}}

\def\PM{{\mathbb{P}}}

\def\QM{{\mathbb{Q}}}
\def\FM{{\mathbb{F}}}
\def\ZM{{\mathbb{Z}}}
\def\WM{{\mathbb{W}}}
\def\CM{{\mathbb{C}}}
\def\AM{{\mathbb{A}}}
\def\XM{{\mathbb{X}}}
\def\MM{{\mathbb{M}}}
\def\MG{{\mathfrak M}}
\def\BG{{\mathfrak B}}

\def\ZG{{\mathfrak Z}}

\def\AC{{\mathcal A}}

\def\CC{{\mathcal C}}
\def\HC{{\mathcal H}}
\def\NC{{\mathcal N}}

\def\IC{{\mathcal I}}
\def\OC{{\mathcal O}}
\def\MC{{\mathcal M}}

\def\PC{{\mathcal P}}

\def\LC{{\mathcal L}}

\def\FC{{\mathcal F}}

\def\GC{{\mathcal G}}
\def\UC{{\mathcal U}}

\def\TC{{\mathcal T}}
\def\DC{{\mathcal D}}

\def\ssi{si et seulement si}

\def\levi{sous-groupe de Levi }

\def\simto{\buildrel\hbox{$\sim$}\over\longrightarrow}

\def\leq{\leqslant}
\def\geq{\geqslant}
\def\injo{\hookrightarrow}
\def\id{\mathop{\mathrm{Id}}\nolimits}

\def\ba{\backslash}
\def\wt{\widetilde}
\def\wh{\widehat}
\def\o#1{\overline{#1}}

\def\application#1#2#3#4#5{\begin{array}{rcl}
                            #1 \;\;\; #2 & \to &  #3 \\
                              #4 & \mapsto & #5 
                            \end{array}} 
\def\cas#1#2#3#4#5{\begin{array}{rcl} #1 \; & = &
    \left\{\begin{array}{rcl} #2 & \hbox{ si } & #3 \\
                             #4 & \hbox{ si } & #5 \end{array}
                         \right. \end{array}}

\def\To#1{\buildrel\hbox{\tiny{$#1$}}\over\longrightarrow}
\def\to{\rightarrow}

\def\ker{\mathop{\hbox{\sl ker}\,}}

\def\im{\mathop{\hbox{\sl im}\,}}

\def\hom#1#2#3{\hbox{\sl Hom}_{#3}\>\!\left(#1,#2\right)}
\def\endo#1#2{\hbox{\sl End}_{#1}\>\!\left(#2\right)}
\def\aut#1#2{\hbox{\sl Aut}_{#1}(#2)}

\def\Mo#1#2{\mathop{\hbox {\sl Mod}_{#1}(#2)}}%cat des RG-modules lisses
\def\Irr#1#2{\mathop{\hbox {\sl Irr}_{#1}\left(#2\right)}}%RG-modules irreductibles
\def\Cu#1#2{\mathop{\hbox {\sl Cusp}_{#1}\left(#2\right)}}%   cuspidales
\def\Disc#1#2{\mathop{\hbox {\sl Disc}_{#1}\left(#2\right)}}%RG-modules irreductibles
\def\Ell#1#2{\mathop{\hbox {\sl Ell}_{#1}\left(#2\right)}}
\def\ind#1#2#3{\hbox {\sl Ind}_{#1}^{#2}\>\!\left(#3\right)}  %induction
\def\cInd#1#2{\hbox {\sl ind}_{#1}^{#2}}
\def\cind#1#2#3{\hbox {\sl ind}_{#1}^{#2}\>\!\left(#3\right)} %ind a supports compacts
\def\ip#1#2#3{\hbox {\sl i}_{#1}^{#2}\>\!(#3)}  %induction parabolique
\def\Ip#1#2{\hbox {\sl i}_{#1}^{#2}}

\def\dim{\mathop{\mbox{\sl dim}}\nolimits}

\def\gr{\mathop{\mbox{\sl gr}}\nolimits}

\def\gal{\mathop{\mbox{\sl Gal}}\nolimits}

\def \ext#1#2#3#4{\mbox{\sl Ext}^{#1}_{#4}\>\!\left(#2,#3\right)}
\def \tor#1#2#3#4{\mbox{\sl Tor}^{#4}_{#1}\>\!(#2,#3)}

\def \limi#1{\lim\limits_{\displaystyle\longrightarrow\atop {#1}}}
\def \limp#1{\lim\limits_{\displaystyle\longleftarrow\atop {#1}}}
\def \limproj{{\lim\limits_{\longleftarrow}}}

\makeatletter
\renewcommand{\subsubsection}{\@startsection{subsubsection}{3}{0mm}{-\baselineskip}{-0.01\baselineskip}{\bf}}%\normalfont\normalsize\itshape}}
\makeatother

\def\ali{\subsubsection{}}
\def\alin#1{\ali{\sl #1}\ : }

\def\ini{\setcounter{equation}{\value{subsubsection}}\addtocounter{subsubsection}{1}}

\newtheorem{theo}[subsubsection]{Th{\'e}or{\`e}me}
\newtheorem{lemme}[subsubsection]{Lemme}
\newtheorem{prop}[subsubsection]{Proposition}
\newtheorem{coro}[subsubsection]{Corollaire}

\newtheorem{fact}[subsubsection]{Fait}

\newtheorem{conj}[subsubsection]{Conjecture}
\newtheorem{DEf}[subsubsection]{D{\'e}finition}

\newtheorem{rema}[subsubsection]{Remarque}  %\renewcommand{\therema}{}
\newtheorem{nota}[subsubsection]{Notation}  %\renewcommand{\thenota}{}
\newtheorem{conv}[subsubsection]{Convention}  %\renewcommand{\thenota}{}

\newtheorem{conjintro}{Conjecture}

\newtheorem{theointro}[conjintro]{Th{\'e}or{\`e}me}

\newtheorem{propsec}[subsection]{Proposition}

\newcommand{\findem}{\hfill$\Box$\par\medskip}
\newcommand{\dem}{\noindent {\sl Preuve :} \rm }

\newenvironment{proof}{\dem}{\findem}

\usepackage[french]{babel}
\usepackage[arrow,matrix]{xy} 
\usepackage[latin1]{inputenc}

\title{Th\'eorie de Lubin-Tate non-ab\'elienne et repr\'esentations elliptiques}
\author{J.-F. Dat}

\begin{document}
\maketitle
\bibliographystyle{plain}

\def\la{\langle}
\def\ra{\rangle}
\def\knr{{\wh{K^{nr}}}}
\def\ka{\wh{K^{ca}}}

\abstract{Non abelian Lubin-Tate theory studies the cohomology of some moduli spaces for $p$-divisible groups, the broadest  definition of which is due to Rapoport-Zink, aiming both at providing explicit realizations of local Langlands functoriality and at studying bad reduction of Shimura varieties. In this paper we consider the most famous examples ; the so-called Drinfeld and Lubin-Tate towers. In the Lubin-Tate case, Harris and Taylor proved that the supercuspidal part of the cohomology realizes both the local Langlands and Jacquet-Langlands correspondences, as conjectured by Carayol. Recently, Boyer computed the remaining part of the cohomology and exhibited two defects : first, the representations of $GL_d$ which appear are of a very particular and restrictive form ; second, the Langlands correspondence is not realized anymore. In this paper, we study the cohomology complex in a suitable equivariant derived category, and show how it encodes Langlands correspondence for {\em elliptic} representations. Then we transfer this result to the Drinfeld tower via an enhancement of a theorem of Faltings due to Fargues. We deduce that Deligne's weight-monodromy conjecture is true for varieties uniformized by Drinfeld's coverings of his symmetric spaces. This  completes the computation of local $L$-factors of some unitary Shimura varieties.    }

\def\mdro{\MC_{Dr,0}}
\def\mdrn{\MC_{Dr,n}}
\def\mdr{\MC_{Dr}}
\def\mlto{\MC_{LT,0}}
\def\mltn{\MC_{LT,n}}
\def\mlt{\MC_{LT}}
\def\mltK{\MC_{LT,K}}

\tableofcontents

\section{Introduction}

\def\dd{D_d^\times}

\subsection{Un peu d'histoire}

Soit $K$ un corps local de caract{\'e}ristique r{\'e}siduelle $p$, $K^{ca}$ une cl\^oture alg\'ebrique et $W_K$ le groupe de Weil associ\'e. Dans leur article \cite{LT1} de 1965, Lubin et Tate se sont inspir\'es de la th\'eorie de la multiplication complexe des courbes elliptiques pour expliciter de mani\`ere exclusivement locale la loi de r\'eciprocit\'e  du corps de classes d'Artin pour $K$. Ils ont pour cela \'etudi\'e certains groupes formels munis d'une action de l'anneau des entiers  $\OC_K$ de $K$. Sur une cl\^oture alg\'ebrique $k^{ca}$ du corps r\'esiduel $k$ de $\OC_K$, les $\OC_K$-modules formels de dimension $1$ sont classifi\'es par leur ``hauteur" $d$. Les auteurs montrent que celui de hauteur $d=1$ se rel\`eve uniquement \`a isomorphisme pr\`es sur la compl\'etion $\wh{K^{nr}}=\WM(k^{ca})$ de l'extension non-ramifi\'ee maximale de $K$ ; le $\OC_K$-module form\'e par les $K^{ca}$-points de $\OC_K$-torsion d'un tel rel\`evement  est isomorphe \`a $K/\OC_K$ et muni d'une action de l'inertie $I_K\subset W_K$, d'o\`u un morphisme $I_K\To{} \OC_K^\times$ qui s'av\`ere induire l'isomorphisme du corps de classes (restreint \`a l'inertie).

La ``th\'eorie de Lubin-Tate non-ab\'elienne", ainsi baptis\'ee par Carayol dans \cite{CarAA}, vise \`a expliciter de mani\`ere locale certaines {\em correspondances} appartenant au vaste programme de
g{\'e}n{\'e}ralisation non-ab{\'e}lienne de la th{\'e}orie du corps de classes propos{\'e}
par Langlands, d\`es 1967. On s'int\'eresse ici en  particulier \`a 
\begin{itemize}
        \item la {\em correspondance de Langlands} qui pour tout entier $d>0$ met en bijection les classes de repr{\'e}senta\-tions {\em lisses} irr{\'e}ductibles de $GL_d(K)$ et les classes de repr{\'e}sentations {\em continues} de dimension $d$ de $W_K$ ; nous la noterons $\pi \mapsto \sigma_d(\pi)$ et renvoyons \`a \cite{LRS} \cite{HaTay} \cite{HeLang} \cite{HeSMF} et \ref{corLan}.
        \item La {\em correspondance de Jacquet-Langlands} qui pour tout entier $d>0$ met en bijection les classes de repr\'esentations {\em lisses} irr\'eductibles du groupe des inversibles de l'alg\`ebre \`a division $D_d$ de centre $K$ et invariant $1/d$ avec les ``s\'eries discr\`etes" de $GL_d(K)$. Nous la noterons $\rho\mapsto JL_d(\rho)$ ainsi que $\pi\mapsto LJ_d(\pi)$ son ``inverse", et renvoyons \`a \cite{DKV}
\cite{Badu1} et  \ref{JacLan}.
\end{itemize}
Les coefficients des repr\'esentations sont ici $l$-adiques pour un premier $l\neq p$. Rappelons n\'eanmoins que ces correspondances sont d\'efinies initialement en termes de repr\'esentations {\em complexes} (et en rempla{\c c}ant $W_K$ par le groupe de Weil-Deligne), qu'elles sont caract\'eris\'ees par  des propri{\'e}t{\'e}s de pr{\'e}servation d'invariants de
nature arithm{\'e}tico-analytique, mais qu'elles  se transf\`erent (presque) sans ambigu\"it\'e \`a tout corps de coefficients abstraitement isomorphe \`a $\CM$, {\em cf} \ref{corLan}.

Comme dans la th\'eorie de Lubin-Tate ab\'elienne, la r\'ealisation explicite de ces correspondances est obtenue gr\^ace aux points de torsion d'un certain $\OC_K$-module formel, sauf que celui-ci va d\'esormais vivre sur un $\knr$-espace analytique de dimension $d-1$.  Il y a en fait deux constructions d'un tel $\OC_K$-module formel.

\begin{itemize}
\item La premi\`ere est une g\'en\'eralisation directe de la situation ab\'elienne ; lorsque $d>1$, Lubin et Tate dans \cite{LT2} montrent que ``le" $\OC_K$-module formel $\HM_d$ de dimension $1$ et hauteur $d$ sur $k^{ca}$ ne se rel\`eve plus de mani\`ere unique \`a $\wh{K^{nr}}$ mais que son foncteur des d\'eformations est repr\'esentable par un anneau de s\'eries formelles \`a $d-1$ variables sur $\wh{K^{nr}}$. Ainsi la boule unit\'e ouverte $\knr$-analytique de dimension $d-1$ est munie d'un $\OC_K$-module formel ``universel" dont les points de torsion forment un 
$\OC_K$-module (ind)\'etale dont les fibres sont isomorphes \`a  $(K/\OC_K)^d$. Le classifiant des trivialisations de cet $\OC_K$-module est donc un (pro)rev\^etement galoisien de groupe $GL_d(\OC_K)$ %, que nous noterons $\MC_{LT}^{d,ca}$. Par 
qui, par fonctorialit\'e, est aussi muni d'une action commutante du groupe des inversibles de l'anneau des endomorphismes de $\HM_d$, lequel se trouve \^etre l'anneau des entiers $\OC_{D_d}$ de  $D_d$. Notons $\MC^{d,ca}_{LT}$ le changement de base de ce pro-rev\^etement \`a la compl\'etion $\wh{K^{ca}}$ de $K^{ca}$ ; celui-ci est donc muni d'une action de $GL_d(\OC_K)\times \OC_{D_d}^\times \times I_K$, et Drinfeld a expliqu\'e dans \cite{Drinell}, {\em cf \ref{actionLT}} comment prolonger cette action\footnote{On triche un peu ici ; il faut d\'eformer par isog\'enies pour avoir une action du produit triple en entier.} \`a  $GL_d(K)\times \dd \times W_K$.

\item La deuxi\`eme construction est enti\`erement d\^ue \`a Drinfeld, qui a introduit dans \cite{Drincov} un autre probl\`eme de ``d\'eformations" de $\OC_K$-modules formels, o\`u les r\^oles de $GL_d(K)$ et $\dd$ sont invers\'es. Il consid\`ere des $\OC_K$-modules formels de dimension $d$ et hauteur $d^2$  munis d'une action de $\OC_{D_d}$ ``sp\'eciale". Sur $k^{ca}$, il existe un tel objet, disons $\XM_d$ qui est unique \`a quasi-isog\'enie pr\`es, et dont le groupe des quasi-isog\'enies est $GL_d(K)$. Le probl\`eme de ``d\'eformations par quasi-isog\'enies" de $\XM_d$ est repr\'esentable par un sch\'ema formel dont la fibre g\'en\'erique est le fameux espace sym\'etrique de Drinfeld $\Omega^{d-1}$ ($\PM^{d-1}_{\knr}$ priv\'e des hyperplans $K$-rationnels) qui, par le jeu des trivialisations des points de torsion du $\OC_{D_d}$-module universel, se voit donc muni d'un pro-rev\^etement de groupe $\OC_{D_d}^\times$. Notons encore $\MC_{Dr}^{d,ca}$ le changement de base \`a $\wh{K^{ca}}$ d'icelui ; il est muni d'une action de $I_K\times \OC_{D_d}^\times \times GL_d(K)$ que l'on peut prolonger \`a $GL_d(K)\times \dd\times W_K$.
\end{itemize}

On dispose maintenant de plusieurs th\'eories cohomologiques pour lin\'eariser $\MC^{d,ca}_{LT}$ et $\MC^{d,ca}_{Dr}$ ; nous adopterons ici la ``cohomologie $l$-adique \`a supports compacts" de Berkovich. Notons que dans le cas $LT$, elle est duale des cycles \'evanescents \'etudi\'es par Harris-Taylor et Boyer.
Les $\o\QM_l$-espaces obtenus sont de dimension infinie et l'action de $GL_d(K)\times \dd$ y est lisse, par des r\'esultats g\'en\'eraux de Berkovich. Pendant longtemps, ces espaces ont \'et\'e \'etudi\'es de chaque c\^ot\'e ({\em LT} et {\em Dr}) de mani\`ere ind\'ependante. Mais  dans un tr\`es court article \cite{FaltDrin}, Faltings a esquiss\'e une preuve de ce que les math\'ematiciens impliqu\'es commen{\c c}aient \`a soup{\c c}onner : {\em il existe des isomorphismes $GL_d(K)\times\dd\times W_K$-\'equivariants}\footnote{Il faut faire attention \`a la normalisation des actions pour obtenir l'\'equivariance.}
%  vraiment l'isomorphisme ci-dessus dans le cas de $\MC^{d,ca}_{LT}$, il 
%convient de conjuguer par l'automorphisme ext{\'e}rieur $g\mapsto
%{^t}g^{-1}$  de $G_d$ l'action d{\'e}finie par Carayol dans \cite{CarAA}.}

$$ H^i_c(\MC^{d,ca}_{LT},\o\QM_l) \simeq H^i_c(\MC^{d,ca}_{Dr},\o\QM_l) .$$
Les arguments de Faltings sont repris et compl\'et\'es par Fargues dans \cite{FarFal}, ainsi que par Genestier et V. Lafforgue en \'egales caract\'eristiques (travail en cours). Dans la suite de cette introduction, on fixe l'entier $d$ et on se contente  de d\'esigner par $H^i_c(\MC^{ca},\o\QM_l)$ ou m\^eme $H^i_c$ ces espaces de cohomologie.

 %Dans ce texte, on s'int{\'e}resse {\`a} deux tours d'espaces $\ka$-analytiques
%d{\'e}finies par Drinfeld.
%\begin{itemize}
%\item la tour $(\MC^{ca}_{Dr,n})_{n\in\NM}$ au-dessus du $K$-espace sym{\'e}trique
%de Drinfeld de dimension $d-1$, d{\'e}finie dans \cite{Drincov} et
%rappel{\'e}e en \ref{tourDr},
%\item  la tour $(\MC^{ca}_{LT,n})_{n\in\NM}$ 
%``de Lubin-Tate'' au-dessus de la boule ouverte de dimension $d-1$ sur $K$,
%d{\'e}finie dans \cite{Drinell} et rappel{\'e}e en \ref{tourLT}.
%\end{itemize}
%Chacun des pro-$\ka$-espaces rigides $\MC^{ca}_{Dr}$ et $\MC^{ca}_{LT}$ obtenus en
%passant {\`a} la limite est muni
%d'une action de $G_d\times \dd\times W_K$, donc leur ``cohomologie {\'e}tale
%$l$-adique {\`a} supports compacts'' fournit des repr{\'e}sentations de ce produit triple. La
%d{\'e}finition habituelle de celle-ci est simplement, pour $?=Dr$ ou $LT$, 
%$$H^i_c(\MC^{ca}_?,\o\QM_l):= \limi{n} H^i_c(\MC^{ca}_{?,n},\o\QM_l),$$
%compte tenu de la finitude des morphismes de transition.
%On sait par des r{\'e}sultats g{\'e}n{\'e}raux de Berkovich que l'action de
%$G_d\times \dd$ sur ces espaces est lisse. 

Rappelons que par la correspondance de Langlands, les repr\'esentations {\em irr\'eductibles} de $W_K$ correspondent aux repr\'esentations {\em supercuspidales} de $GL_d(K)$, lesquelles sont caract\'eris\'ees par la propri\'et\'e d'\^etre {\em projectives} et {\em injectives} dans la cat\'egorie des repr\'esentations lisses de $GL_d(K)$ \`a caract\`ere central. La partie {\em supercuspidale} des espaces $H^*_c(\MC^{ca},\o\QM_l)$ est donc \`a la fois la plus importante et la plus maniable \`a \'etudier. Elle a \'et\'e explor\'ee dans les quatre articles \cite{HaCusp} \cite{Boyer1} \cite{HaTay}
\cite{Hausb}, chacun de ces articles concernant le c\^ot\'e {\em LT} ou le c\^ot\'e {\em Dr} en caract{\'e}ristique nulle ou en {\'e}gales
caract{\'e}ristiques. Ils d\'emontrent en particulier la fameuse conjecture 
que Carayol a \'enonc\'ee dans \cite{CarAA} \`a la suite de travaux pionniers de Deligne et Drinfeld pour $d=2$ :
%Le r{\'e}sultat final le plus spectaculaire, anciennement conjecture
%de Deligne-Drinfeld-Carayol \cite{CarAA} est le suivant :
 soit $\pi$ une
repr{\'e}sentation supercuspidale $l$-adique de $GL_d(K)$, de caract{\`e}re central
d'ordre fini ; alors on a, en normalisant convenablement les actions, {\em cf} \cite{HaICM}, 
$$\begin{array}{rcl} \hom{H_c^i(\MC^{ca}_{},\o\QM_l)}{\pi}{GL_d(K)} \; &
  \mathop{\simeq}\limits_{\dd\times W_K} &
%\hom{H_c^i(\MC^{ca}_{LT},\o\QM_l)}{\pi}{GL_d(K)} \\ &
 % \mathop{\simeq}\limits_{\dd\times W_K} &
    \left\{\begin{array}{rcl} LJ_d(\pi)\otimes
  \sigma_d(\pi)(?) & \hbox{ si } & i=d-1 \\
                             0 & \hbox{ si } & i\neq d-1 \end{array}
                         \right. \end{array}$$ % \;\; \hbox{ et }$$
%$$\begin{array}{rcl} \hom{H_c^i(\MC^{ca}_{LT},\o\QM_l)}{JL_d(\pi)}{D} \; &
%  \mathop{\simeq}\limits_{G_d\times W_K} &
%    \left\{\begin{array}{rcl} \pi\otimes
%  \sigma_d(\pi)(?) & \hbox{ si } & i=d-1 \\
%                             0 & \hbox{ si } & i\neq d-1 \end{array}
%                         \right. \end{array}
%$$
o{\`u} $(?)$ d{\'e}signe une certaine torsion {\`a} la Tate. Pr\'ecisons aussi que Harris et Taylor dans \cite{HaTay} prouvent cette formule en m\^eme temps que l'existence m\^eme de la correspondance de Langlands. 
%\footnote{Ce ne sont pas exactement
%  les pr{\'e}dictions de \cite{CarAA}, mais l'espace not{\'e} $\MC^{ca}_{LT}$ ici
%  est , ni la m{\^e}me cohomologie, ni la
%vraiment l'isomorphisme ci-dessus dans le cas de $\MC^{ca}_{LT}$, il
%convient de conjuguer par l'automorphisme ext{\'e}rieur $g\mapsto
%^tg^{-1}$  de $G_d$ l'action d{\'e}finie par Carayol dans \cite{CarAA}.}
Les m{\'e}thodes employ{\'e}es dans ces articles reposent sur la propri{\'e}t{\'e}
d'uniformisation $p$-adique de certaines vari{\'e}t{\'e}s globales (de Shimura
ou de Drinfeld) par ces
espaces, les suites spectrales de type Hochschild-Serre ou de cycles
{\'e}vanescents
associ{\'e}es, et des arguments de formule des traces. Elles ne permettent donc en g{\'e}n{\'e}ral que d'obtenir des
informations sur la somme altern{\'e}e des $H^i_c$. 
%N{\'e}anmoins certains
%calculs explicites sont connus, les premiers 

Pendant longtemps,  seul le calcul de Schneider-Stuhler dans \cite{SS1} pour l'espace sym\'etrique de Drinfeld fournissait des renseignements sur la partie non-cuspidale de la cohomologie. Ce calcul a inspir\'e \`a Harris une conjecture (non publi\'ee) sur la forme explicite des groupes de  cohomologie individuels du c\^ot\'e $Dr$. %  sur un corps $p$-adique. Tout  
R\'ecemment, Boyer \cite{Boyer2} a annonc\'e une preuve de cette conjecture, en travaillant du c\^ot\'e $LT$. % et en utilisant le th\'eor\`eme de Faltings-Fargues. % et en \'egales caract\'eristiques. Un peu auparavant, Faltings avait annonc\'e dans \cite{FaltDrin} que les  cohomologies des deux tours devaient \^etre isomorphes.
%Comme ses  arguments sont tr\`es r\'esum\'es et valables seulement sur un corps $p$-adique, ils seront compl\'et\'es et \'etendus en  \'egales caract\'eristiques dans un travail en cours de Fargues, Genestier et Lafforgue.
%{\'e}tant d{\^u}s depuis longtemps {\`a} Schneider-Stuhler \cite{SS1}, les plus r{\'e}cents {\`a} Boyer \cite{Boyer2}. Ce
%dernier d{\'e}crit compl{\`e}tement les $H^i_c$ dans le cas $\MC^{ca}_{LT}$ sur un
%corps de fonctions et montre en particulier l'analogue d'une
%pr{\'e}diction qu'avait faite auparavant Harris (non publi{\'e})  dans le cas de
%$\MC^{ca}_{Dr}$  sur un corps $p$-adique.
Son r\'esultat peut s'exprimer ``qualitativement" comme ceci :
{\em Soit $\pi$ une repr\'esentation irr\'eductible de contribution non-nulle \`a $H^i_c$, alors $\pi$ est {\em elliptique} et 
$$ \hom{H^i_c(\MC^{ca},\o\QM_l)}{\pi}{GL_d(K)} \mathop{\simeq}\limits_{\dd\times W_K}
LJ_d(\pi)\otimes\sigma'(\pi)(?)  $$
o{\`u} $(?)$ d{\'e}signe une certaine torsion {\`a} la Tate.}

Quelques explications sur cette formule : $\pi$ est dite 
{\em elliptique} si la semi-simplifi\'ee de sa correspondante de Langlands est de la forme $\sigma_d(\pi)^{ss}=\sigma'(\pi)\oplus\sigma'(\pi)(1)\oplus \cdots\oplus \sigma'(\pi)(d')$ pour une certaine repr{\'e}sentation $l$-adique irr{\'e}ductible $\sigma'(\pi)$ (voir \ref{defell} pour d'autres caract{\'e}risations plus intrins\`eques). 
%la restriction de son caract{\`e}re-distribution aux {\'e}l{\'e}ments semi-simples r{\'e}guliers elliptiques
%est non-nulle (voir \ref{defell} pour d'autres caract{\'e}risations). 
Pour une telle repr\'esentation, on note $LJ_d(\pi)$ le transfert de Jacquet-Langlands de l'unique s\'erie discr\`ete dont la correspondante de Langlands a la m\^eme semi-simplifi\'ee que $\sigma_d(\pi)$.

\subsection{Cet article}

La description obtenue par Boyer montre deux d\'efauts de la cohomologie des espaces $\MC^{ca}$ :
\begin{itemize}
        \item parmi les repr\'esentations elliptiques, seules {\em certaines} apportent une contribution non-nulle (et d'ail\-leurs, ce ph\'enom\`ene appara\^it d\'eja dans le calcul de Schneider et Stuhler pour $\Omega$).
        \item si la correspondance de Jacquet-Langlands (convenablement \'etendue) est bien r\'ealis\'ee, il n'en va pas de m\^eme de la correspondance de Langlands.
On voit en particulier que l'op{\'e}rateur de monodromie (celui que l'on sait associer {\`a} toute
repr{\'e}sentation $l$-adique continue de dimension finie de $W_K$ par le th{\'e}or{\`e}me
de Grothendieck) est toujours nul sur les composantes isotypiques des
$H^i_c$.
\end{itemize}

L'id\'ee principale de ce texte, qui fait suite \`a \cite{Dat1}, est que pour corriger ces d\'efauts il faut enrichir la cohomologie, ou plut\^ot lui restituer sa richesse perdue, en consid\'erant ``le'' complexe de cohomologie
$R\Gamma_c(\MC^{ca},\o\QM_l)$
 vu comme objet de la cat{\'e}gorie
d{\'e}riv{\'e}e de la cat{\'e}gorie ab{\'e}lienne des $\o\QM_l(GL_d(K)\times
\dd)$-repr{\'e}sentations lisses.  Oublions un instant les 
difficult{\'e}s de d{\'e}finition et d'{\'e}tude que cela pose et {\'e}non{\c c}ons
notre th\'eor\`eme principal  :
\begin{theointro} \label{main} Pour toute repr\'esentation lisse irr\'eductible $\pi$ de $GL_d(K)$, on a
$$\HC^*\left(R\hom{R\Gamma_c(\MC^{ca},\o\QM_l)}{\pi}{D^b(GL_d(K))}\right)
\mathop{\simeq}\limits_{\dd\times
  W_K}LJ_d(\pi)\otimes\sigma_d(\pi)|-|^{\frac{d-1}{2}}.$$
\end{theointro}
Faisons quelques remarques sur les  deux termes de cet isomorphisme :
\begin{itemize}
        \item Du c\^ot\'e gauche : les complexes $R\Gamma_c \in D^b_{\o\QM_l}(GL_d(K))$ et (par cons\'equent) $R\hom{R\Gamma_c}{\pi}{D^b(GL_d)}\in D^b(\o\QM_l)$, sont munis d'une action de $\dd\times W_K$. La notation $\HC^*$ d\'esigne l'\'equivalence de cat\'egorie 
          $\HC^*: \CC^\bullet \in D^b(\o\QM_l) \mapsto \bigoplus_{i\in\ZM}
\HC^i(\CC^\bullet)$ entre $D^b(\o\QM_l)$ et la cat\'egorie des espaces vectoriels \`a graduation de support fini, convenablement triangul\'ee. Ainsi le terme de gauche est un espace vectoriel gradu\'e muni d'une action de $\dd\times W_K$, mais on oublie la graduation pour se retrouver avec une honn\^ete repr\'esentation lin\'eaire de $\dd\times W_K$.
  \item Du c\^ot\'e droit : seule la notation $LJ_d$ demande une explication que voici, {\em cf} \ref{JacLan} : la correspondance de Jacquet-Langlands fournit un plongement des groupes de Grothendieck $R(\dd) \injo R(GL_d(K))$ qui induit un isomorphisme $R(\dd)\simto \o{R}(GL_d(K))$ o\`u $\o{R}$ d\'esigne le quotient par les combinaisons lin\'eaires d'induites paraboliques. D'o\`u une r\'etraction canonique $R(GL_d(K))\To{} R(\dd)$ que nous notons\footnote{Cette notation est emprunt\'ee \`a Badulescu qui a d\'efini dans \cite{Badu2} des r\'etractions similaires dans des situations plus compliqu\'ees o\`u $D_d$ est remplac\'ee par une alg\`ebre centrale simple quelconque.} $LJ_d$. On v\'erifie alors que pour toute $\pi$ irr\'eductible, $LJ_d(\pi)$ est non-nulle \ssi\ $\pi$ est elliptique, et dans ce cas, la notation est coh\'erente avec celle introduite plus haut, au signe pr\`es. En particulier, le th\'eor\`eme \ref{main} ne fournit une construction de la correspondance de Langlands {\em que pour les elliptiques}.
\end{itemize}

Il est raisonnable de penser que les espaces $\MC^{ca}$ ne {\em peuvent pas} donner  de r\'ealisation, en quelque sens cohomologique que ce soit\footnote{comme l'indiquent les calculs de traces de \cite{Faltrace}}, de la correspondance de Langlands pour les repr\'esentations non-elliptiques, et donc de se demander quels espaces pourraient fournir une telle r\'ealisation. Les candidats sont naturellement \`a chercher parmi les espaces de Rapoport-Zink \cite{RZ} qui, rappelons-le, font toujours intervenir deux groupes. Dans cette perspective, on peut reformuler le th\'eor\`eme \ref{main} en introduisant le groupe
$GD:=(GL_d(K)\times \dd)/\Delta$ o\`u $\Delta=\{(z,z),z\in K^\times\}$ sous la forme : 
%On sait que, quitte {\`a} conjuguer l'action de $G_d$  initialement d{\'e}finie
%par Drinfeld par un automorphisme ext{\'e}rieur convenable, l'action
%g{\'e}om{\'e}trique de $G_d\times \dd$ sur les espaces 
%de modules $\MC^{ca}_{Dr}$ et $\MC^{ca}_{LT}$ se factorise par ce groupe. De plus,
%toute repr{\'e}sentation irr{\'e}ductible de $GD$ est de la forme
%$\pi\otimes \rho$ o{\`u} $\pi$ et $\rho$ sont des repr{\'e}sentations
%irr{\'e}ductibles de $G_d$ et $\dd$ dont les caract{\`e}res centraux sont
%``inverses l'un de l'autre''. 
{\em pour toutes repr\'esentations irr\'eductibles $\pi$ de $GL_d(K)$ et $\rho$ de $\dd$, on a
$$\HC^*\left(R\hom{R\Gamma_c(\MC^{ca},\o\QM_l)}{\pi\otimes
    \rho}{D^b(G D)}\right)
\mathop{\simeq}\limits_{W_K} \left\{\begin{array}{rl} 
  \sigma_d(\pi)|-|^{\frac{d-1}{2}} & \hbox{ si }  \rho=LJ_d(\pi^\vee) \\
                             0 & \hbox{ sinon } \end{array}\right.. $$ }
La notation $?^\vee$ d{\'e}signe la contragr{\'e}diente de la repr{\'e}sentation
$?$. Dans le cas qui nous int\'eresse ici, les deux formes sont \'equivalentes car les repr\'esentations de $\dd$ sont ``essentiellement" semi-simples, mais pour des raisons de sym\'etrie, c'est la seconde qu'il semble plus naturel de vouloir g\'en\'eraliser.

\bigskip

D\'ecrivons maintenant la strat\'egie que nous suivons : la premi\`ere chose \`a faire est bien-s\^ur de d\'efinir convenablement les deux $R\Gamma_c$, ce que nous faisons au paragraphe  \ref{defRG} et passerons sous silence ici, pour ne pas effrayer les \'eventuels lecteurs. Les d\'efinitions sont faites pour que, d'apr\`es Fargues \cite{FarFal}, il existe  un isomorphisme   $R\Gamma_c(\MC^{ca}_{LT},\o\QM_l) \simeq R\Gamma_c(\MC^{ca}_{Dr},\o\QM_l)$ dans  $D^b(GD)$, compatible  aux actions de $W_K$.

Une fois le complexe $R\Gamma_c$ bien d{\'e}fini, on est en pr{\'e}sence d'un
objet {\em a priori} compliqu{\'e}. Le miracle est que cet objet est en
fait ``aussi simple qu'il peut l'\^etre'' puiqu'il est {\em scindable} :
\begin{propsec}\label{prop1} %{Proposition \ref{}}
Il existe un isomorphisme (pas unique) dans $D^b(GD)$,
 $$R\Gamma_c(\MC^{ca},\o\QM_l)
    \simto \bigoplus_{i\in \NM} H^i_c(\MC^{ca},\o\QM_l)[-i].$$
 \end{propsec}
Ce fait repose sur la description de Boyer et peut se prouver de deux mani\`eres : soit par un argument de th{\'e}orie des repr{\'e}sentations de $GL_d(K)$ (le calcul complet des
groupes d'extensions $\ext{i}{\pi}{\pi'}{G_d}$ pour les couples de repr{\'e}sentations
elliptiques $\pi$ et $\pi'$ de $GL_d(K)$), soit en  utilisant l'action d'un rel\`evement de Frobenius sur les $H^i_c$.

%Notons $\o{GL_d}:=GL_d(K)/\varpi^\ZM = GD/\dd$ et f
Fixons une irr\'eductible $\rho$ de $\dd$, notons $\omega_\rho :K^\times \To{} \o\QM_l$ son caract\`ere central, et $D^b_{\omega_\rho}(GL_d)$ la cat\'egorie d\'eriv\'ee born\'ee de la cat\'egorie des repr\'esentations lisses de $GL_d$ de caract\`ere central $\omega_\rho$. 
Le complexe $R\Gamma_c[\rho] := R\Gamma_c \otimes^L_{\dd} \rho \in D^b_{\omega_\rho}({GL_d})$  est scindable en $R\Gamma_c[\rho]\simeq \bigoplus_i H^i_c[\rho][-i]$ o\`u $H^i_c[\rho]:=H^i_c\otimes_{\dd} \rho$. % d\'esigne la partie $\rho^\vee$-isotypique  de la cohomologie $H^i_c$.
Notons que du c\^ot\'e $Dr$, ce complexe  s'interpr\`ete g\'eom\'etriquement comme le complexe de cohomologie $R\Gamma_c(\Omega^{d-1},\LC_\rho)$ du syst\`eme local $l$-adique $\LC_\rho$  associ\'e \`a $\rho$ via les structures de niveau de Drinfeld.

La propri\'et\'e de scindage permet de d\'ecrire l'alg\`ebre des endomorphismes de $R\Gamma_c[\rho]$ sous la forme :
$$ \endo{D^b_{\omega_\rho}({GL_d})}{R\Gamma_c[\rho]} \simeq \bigoplus_{i\geq j}
\ext{i-j}{H^i_c[\rho]}{H^j_c[\rho]}{GL_d,\omega_\rho},$$
le produit sur le $\o\QM_l$-espace vectoriel de droite {\'e}tant donn{\'e} par
le $\cup$-produit.
Gr\^ace \`a la th\'eorie de Bushnell-Kutzko et au th\'eor\`eme 1.3 de \cite{Dat1}, on peut d\'ecrire compl\`etement ce $\cup$-produit. Notons 
$(\sigma_\rho, V_{\sigma_\rho})$ la repr{\'e}sentation $l$-adique de
$W_K$ associ{\'e}e {\`a} $\rho$ 
par correspondance de Langlands et fonctorialit{\'e} de Jacquet-Langlands
: il s'agit d'une repr{\'e}sentation ind{\'e}composable de dimension $d$ et
l'alg{\`e}bre engendr{\'e}e par
l'image de $W_K$ dans $\endo{\o\QM_l}{V_{\sigma_\rho}}$ est une alg{\`e}bre
triangulaire par blocs que nous noterons $\AC_\rho$.
Le calcul de $\cup$-produits mentionn\'e ci-dessus montre alors ({\em cf} \ref{propnoneq})

\begin{propsec} \label{prop2}
Tout scindage de $R\Gamma_c$ induit un isomorphisme de $\o\QM_l$-alg{\`e}bres
$$\endo{D^b_{\omega_\rho}(GL_d)}{R\Gamma_c[\rho]} \simto %\prod_{\rho\in
\AC_{\rho^\vee}.$$
\end{propsec}

On veut ensuite {\'e}lucider l'action de $W_K$ sur le membre de gauche.
Nous commen{\c c}ons par montrer que l'inertie agit potentiellement par son quotient $l$-adique et que l'on peut d\'efinir la ``partie unipotente de la monodromie" $N_\rho\in \endo{D^b_{\omega_\rho}(GL_d)}{R\Gamma_c[\rho]}$. 
 Pour d{\'e}crire pr{\'e}cis{\'e}ment ce $N_\rho$, nous aurons besoin
de minorer son ordre de nilpotence. {\em \`A ce point nous avons besoin des r\'esultats globaux de Boyer} : nous utiliserons sa description du gradu\'e pour la filtration de monodromie du complexe des cycles \'evanescents de certaines vari\'et\'es de Shimura (ou de Drinfeld). En notant  $\gamma_\rho : W_K \To{}
\endo{D^b_{\omega_\rho}(GL_d)}{R\Gamma_c[\rho]}^\times$ le morphisme  de groupes donnant l'action de $W_K$ sur
$R\Gamma_c[\rho]$, on obtient la description
suivante, {\em cf} \ref{propeq}
\begin{propsec} \label{prop3}
 Il existe un scindage de $R\Gamma_c$ pour lequel l'isomorphisme d'alg{\`e}bres pr\'ec\'edent rend le diagramme suivant commutatif :
$$\xymatrix{
\endo{D^b_{\omega_\rho}(GL_d)}{R\Gamma_c[\rho]}
\ar[r]^-\sim   & \AC_{\rho^\vee} \\ W_K \ar[u]^{\gamma_\rho} \ar[ur]^{\sigma_{\rho^\vee}} &
}$$
\end{propsec}
{\`A} partir de cette proposition, le th{\'e}or{\`e}me \ref{main} se prouve
en utilisant {\`a} nouveau le calcul explicite des extensions entre
repr{\'e}sentations elliptiques et de leurs $\cup$-produits.

Dans le cas $\rho=1$, le complexe \`a \'etudier est simplement $R\Gamma_c(\Omega_K^{d-1},\o\QM_l)$. Nous l'avons d\'eja \'etudi\'e dans \cite{Dat1} en suivant la m\^eme strat\'egie et les m\^emes \'etapes que ci-dessus. Le point essentiel o\`u les arguments divergent est l'\'etude de l'ordre de nilpotence du $N_\rho$. Dans \cite{Dat1}, nous avons utilis\'e la suite spectrale de Rapoport-Zink et l'uniformisation $p$-adique, alors qu'ici nous utilisons les r\'esultats de Boyer, et donc {\em in fine} la g\'eom\'etrie du c\^ot\'e Lubin-Tate. 
Notons d'ailleurs que toutes les \'etapes qui pr\'ec\`edent s'appliquent \`a la tour de Lubin-Tate ind\'ependamment du th\'eor\`eme de Faltings-Fargues. Ce dernier intervient en dernier lieu pour transf\'erer cette \'etude du c\^ot\'e Drinfeld, ce qui nous am\`ene \`a notre deuxi\`eme r\'esultat principal

\begin{theointro} \label{main2}
Les vari\'et\'es uniformis\'ees par les rev\^etements de l'espace sym\'etrique de Drinfeld $\Omega$ satisfont la conjecture monodromie-poids de Deligne.
\end{theointro}

Ceci termine le calcul du facteur $L$ local des vari\'et\'es de Shimura unitaires \'etudi\'ees dans \cite{HaCusp}, \cite{CarAA} et \cite{RapAA} en une place o\`u l'alg\`ebre \`a division et involution globale associ\'ee reste totalement ramifi\'ee.

\bigskip

D\'ecrivons bri\`evement le contenu des diff\'erentes parties.
Dans la  partie \ref{elliptiques} on d\'efinit et caract\'erise les
repr\'esentations elliptiques et on explicite leur comportement \`a
travers les correspondances de Langlands et Jacquet-Langlands.
 Dans la  partie \ref{drinfeld} on pr\'ecise la version des espaces de Drinfeld et Lubin-Tate utilis\'ee pour les \'enonc\'es cohomologiques principaux, puis on d\'efinit leur complexe de cohomologie, ainsi que plusieurs variantes utiles par la suite. On y \'enonce la version du th\'eor\`eme de Faltings-Fargues pertinente pour cet article.
La  partie \ref{real} contient la preuve du
th\'eor\`eme \ref{main} sous l'hypoth\`ese  que l'op\'erateur $N_\rho$ introduit ci-dessus est d'ordre assez grand ; la strat\'egie est la m\^eme que dans la partie $4.2$ de \cite{Dat1}.
Enfin la partie \ref{monodromie} s'occupe de monodromies : d'une part on prouve l'estimation n\'ecessaire \`a la partie \ref{real} en se raccrochant aux r\'esultats globaux de Boyer au prix d'acrobaties techniques m\^elant topos fibr\'es, faisceaux pervers et formalisme $l$-adique. D'autre part on prouve le th\'eor\`eme \ref{main2} ci-dessus. Enfin l'appendice contient un r\'esultat technique, mais d'int\'er\^et ind\'ependant utilis\'e dans la partie \ref{monodromie}.

\section{Repr{\'e}sentations elliptiques et correspondances}\label{elliptiques}

Dans cette partie, nous rappelons des faits bien connus sur les
correspondances de Langlands et Jacquet-Langlands, d'autres un peu moins
connus mais qui le sont certainement des sp{\'e}cialistes, puis nous explicitons  le cas des
repr{\'e}sentations elliptiques. Nous abr\`egerons $G_d:=GL_d(K)$.

\subsection{Correspondance de Jacquet-Langlands locale} \label{JacLan}

\alin{Notations}
Il sera commode de consid{\'e}rer  des repr{\'e}sentations {\`a}
coefficients dans un corps $C$ abstraitement isomorphe au corps des
complexes $\CM$. Ce corps pourra parfois {\^e}tre $\CM$ ou un $\o\QM_l$,
selon les besoins topologiques qu'on aura. Pour tout groupe $H$
localement profini, nous notons $\Irr{C}{H}$ l'ensemble des classes 
d'isomorphisme de $C$-repr{\'e}sentations {\em lisses } irr{\'e}ductibles de
$H$. Pour le groupe $G_d$, on isole certains sous-ensembles
remarquables :
\begin{enumerate}
\item On d{\'e}signe par  $\Cu{C}{G_d}\subseteq \Irr{C}{G_d}$ le sous-ensemble
  des repr{\'e}sentations {\em
    cuspidales}, {\em i.e.} dont les coefficients sont {\`a} support
  compact-modulo-le-centre. Parmi les irr\'eductibles, elles se caract\'erisent aussi comme \'etant  les objets injectifs et projectifs de la cat\'egorie des repr\'esentations lisses de $G_d$ \`a caract\`ere central.  
\item On note $\Disc{\CM}{G_d} \subseteq
  \Irr{\CM}{G_d}$ le sous-ensemble des repr{\'e}sentations ``de la s{\'e}rie
  discr{\`e}te'', {\em i.e.} dont les coefficients sont de carr{\'e} int{\'e}grable (au sens
  complexe) modulo-le-centre. 
Malgr{\'e} cette d{\'e}finition de nature analytique, il se trouve que la
notion de ``s{\'e}rie discr{\`e}te'' de $G_d$  est invariante par
automorphismes du corps $\CM$ ;  c'est une cons{\'e}quence des th{\'e}or{\`e}mes
de multiplicit{\'e}s limites de \cite{Rog}, ou plus prosa{\"\i}quement une
cons{\'e}quence de la classification de Bernstein-Zelevinski \cite[Thm
9.3]{Zel}. Cela nous permet d'isoler sans ambigu{\"\i}t{\'e} un sous-ensemble 
$\Disc{C}{G_d} \subseteq
  \Irr{C}{G_d}$ dont nous appellerons les membres ``s{\'e}ries discr{\`e}tes''
  par abus de langage.
\end{enumerate}

On peut consid{\'e}rer la correspondance de Jacquet-Langlands comme le
reflet spectral de la correspondance ``g{\'e}om{\'e}trique''
bien connue 
entre classes de conjugaison elliptiques semi-simples r{\'e}guli{\`e}res de
$G_d$ et de $\dd$, donn{\'e}e par l'{\'e}galit{\'e} des  polyn{\^o}mes caract\'eristiques. L'{\'e}nonc{\'e} spectral
classique concerne les repr{\'e}sentations complexes  : 
%notons $E^2(?)$ 
%l'ensemble des classes de repr{\'e}sentations complexes irr{\'e}ductibles du
%groupe $?=G,D$ dont
%les coefficients sont de carr{\'e} int{\'e}grable modulo le centre (aussi
%appel{\'e}es ``s{\'e}ries discr{\`e}tes''). Remarquons que
%$E^2(D)=\Irr{\CM}{D}$. 

\begin{theo} (Correspondance de Jacquet-Langlands, \cite{DKV},\cite{Badu1}) \label{JL1}
  Il existe une bijection $$JL_d:\; \Irr{\CM}{\dd}\simto \Disc{\CM}{G_d}$$ caract{\'e}ris{\'e}e par
  l'{\'e}galit{\'e} de caract{\`e}res $\chi_{JL_d(\rho)}(g) =
  (-1)^{d-1}\chi_{\rho}(x)$ pour toutes classes elliptiques $g\in G_d,
  x\in\dd$ se  correspondant ({\em i.e.} ayant m{\^e}me polyn{\^o}me minimal
  de degr{\'e} $d$).
\end{theo}

Remarquons que l'{\'e}galit{\'e} de caract{\`e}res de  cet {\'e}nonc{\'e}
 peut se tester sur les fonctions localement constantes {\`a} support compact dans
l'ensemble (ouvert) des elliptiques r{\'e}guliers, et ne fait donc
intervenir que le caract{\`e}re-distribution des repr{\'e}sentations de
$\Disc{\CM}{G_d}$. Ce caract{\`e}re-distribution est d{\'e}fini sur n'importe
quel corps de coefficients, en particulier sur $C$. Compte tenu de ce
que l'ensemble $\Disc{\CM}{G_d}$ et la condition d'{\'e}galit{\'e} des
caract{\`e}res sont stables par l'action
des automorphismes du corps $\CM$, l'{\'e}nonc{\'e} ci-dessus se transf{\`e}re
sans ambigu{\"\i}t{\'e} {\`a} un {\'e}nonc{\'e} formellement analogue sur le corps $C$. 
Nous noterons encore $JL_d:\; \Irr{C}{\dd}\simto \Disc{C}{G_d}$ la
bijection obtenue.

%(c'est une cons{\'e}quence des th{\'e}or{\`e}mes
%de multiplicit{\'e}s limites de \cite{Rog}). 
%Elle se transf{\`e}re donc sans
%probl{\`e}mes {\`a} tout corps $Q$ alg{\'e}briquement clos de caract{\'e}ristique
%nulle. 
%De plus, le caract{\`e}re d'une repr{\'e}sentation admissible de $G_d$ en
%un {\'e}l{\'e}ment elliptique se d{\'e}finit alg{\'e}briquement et fait sens sur un
%corps g{\'e}n{\'e}ral. La bijection $JL$ est visiblement compatible {\`a} l'action des
%automorphismes de $\CM$ de chaque c{\^o}t{\'e} ; l'{\'e}nonc{\'e} ci-dessus a donc un
%sens clair sur tout corps 
% de coefficients $Q$ alg{\'e}briquement isomorphe {\`a} $\CM$ (par exemple
% $Q=\o\QM_l$, pour $l$ premier).
%De plus toute ``s{\'e}rie discr{\`e}te'' s'obtient de mani{\`e}re unique comme ``Steinberg
%g{\'e}n{\'e}ralis{\'e}e'' $St(\tau,e)$ o{\`u} $e$ divise $d$ et $\tau$ est une
%repr{\'e}sentation cuspidale de $G_dL_{d/e}(K)$.

\alin{Groupes de Grothendieck}
Notons maintenant $R(G_d)$ et $R(\dd)$ les groupes de Grothendieck des
$C$-repr{\'e}sentations de longueur finies. La bijection de Jacquet-Langlands 
induit une injection $R(\dd)\injo R(G_d)$, v{\'e}rifiant
l'{\'e}galit{\'e} de caract{\`e}res du th{\'e}or{\`e}me \ref{JL1} sur les classes de conjugaison
elliptiques se correspondant. On veut d{\'e}finir une r{\'e}traction pour cette injection,
v{\'e}rifiant la m{\^e}me {\'e}galit{\'e} de caract{\`e}res. 
Soit $\chi$ un caract{\`e}re lisse de $K^\times$, notons $R(G_d,\chi)$,
resp. $R(\dd,\chi)$,  le
sous-groupe de $R(G_d)$, resp. $R(\dd)$, engendr{\'e} par les irr{\'e}ductibles de caract{\`e}re
central $\chi$. 
Les d{\'e}compositions selon le caract{\`e}re central $R(G_d)=\bigoplus_\chi
R(G_d,\chi)$, resp. $R(\dd)=\bigoplus_\chi R(\dd,\chi)$, sont respect{\'e}es par
l'application $JL$. De plus, les groupes $R(G_d,\chi)$ et $R(\dd,\chi)$
 sont munis de formes bilin{\'e}aires enti{\`e}res, voir \cite{SS2},
$$\application{\la\;.\;\ra:\;}{R(G_d,\chi)\times
  R(G_d,\chi)}{\ZM}{(\pi,\pi')}{\sum_i (-1)^i\dim(\ext{i}{\pi}{\pi'}{G_d,\chi})}$$
et respectivement pour $\dd$. 
Ces formes bilin{\'e}aires en induisent une sur la somme
directe $R(G_d)$, resp. $R(\dd)$, que l'on note de la m{\^e}me mani{\`e}re.
Dans le cas de $\dd$, cette forme est non
d{\'e}g{\'e}n{\'e}r{\'e}e et on a simplement $\la\rho,\rho'\ra=\dim\hom{\rho}{\rho'}{\dd}$.
Dans le cas de $G_d$, la situation est sensiblement plus compliqu{\'e}e.
Notons 
$R_I(G_d)$ le sous-groupe de $R(G_d)$ engendr{\'e} par les induite
paraboliques de repr{\'e}sentations de sous-groupes de Levi propres et 
$\o{R}(G_d)$ le quotient $R(G_d)/R_I(G_d)$.

\begin{lemme}
 $\o{R}(G_d)$ est libre sur $\ZM$ et  la forme bilin{\'e}aire $\la\;,\;\ra$  s'y descend 
 et y est non-d{\'e}g{\'e}n{\'e}r{\'e}e, une base orthonormale {\'e}tant donn{\'e}e par les
 images de s{\'e}ries discr{\`e}tes. 
\end{lemme}
\begin{proof} 
Pour cette preuve, nous pouvons identifier $C$ et $\CM$.
%Pour nous ramener {\`a} des r{\'e}sultats connus, nous supposerons
%  que le corps de coefficients $C$ est $\CM$, ce qui suffit puisqu'on
%  suppose que $C$ est abstraitement isomorphe {\`a} $\CM$.
%  Dans ce cas, 
L'{\'e}nonc{\'e} est alors une cons{\'e}quence des trois propri{\'e}t{\'e}s suivantes :
  \begin{enumerate}
  \item La classification de Zelevinski \cite{Zel}(ou celle dite du quotient de
    Langlands) montre que
$$R(G_d) =\bigoplus_{(M,\sigma,\psi)} \ZM [\ip{M}{G_d}{\sigma\psi}] $$
o{\`u} $[?]$ est l'{\'e}l{\'e}ment de $R(G_d)$ associ{\'e} {\`a} la repr{\'e}sentation $?$, les
triplets $(M,\sigma,\psi)$ sont form{\'e}s d'un \levi standard de $G_d$
(c'est-{\`a}-dire un produit de $GL_{d_i}$ diagonaux), d'une
repr{\'e}sentation temp{\'e}r{\'e}e de $M$ et d'un caract{\`e}re non-ramifi{\'e} de $M$ ``dans la
chambre de Weil positive'', et sont pris {\`a} $G_d$-conjugaison pr{\`e}s. Le
signe $\Ip{M}{G_d}$ est l'induction parabolique normalis{\'e}e le long du
parabolique triangulaire sup{\'e}rieur dont le Levi est $M$.
Comme on sait de plus que les repr{\'e}sentations temp{\'e}r{\'e}es de $G_d$ sont
soit induites, soit des s{\'e}ries discr{\`e}tes, on en d{\'e}duit
$$R(G_d)= \left(\bigoplus_{\pi\in \Disc{C}{G_d}} \ZM[\pi]\right) \bigoplus
R_I(G_d),$$
ce qui montre que $\o{R}(G_d)$ est libre sur $\ZM$ et qu'une base en est
donn{\'e}e par les images des s{\'e}ries discr{\`e}tes.
\item  Soit $M$ un \levi standard de $G_d$ et $\pi$ une repr{\'e}sentation
  admissible de  $G_d$. Alors pour toute repr{\'e}sentation admissible de
  $G_d$, un argument de d{\'e}formation attribu{\'e} {\`a} Kazhdan dans \cite{SS2} montre que
$$\la \pi,\ip{M}{G_d}{\sigma} \ra = \la \ip{M}{G_d}{\sigma},\pi \ra =0.$$
\item Soit $\pi,\pi'$ deux s{\'e}ries discr{\`e}tes de $G_d$, on a par \cite{Vigext}
$$ \la \pi,\pi'\ra =\delta_{\pi\pi'}.$$
  \end{enumerate}
\end{proof}

D'autre part, la formule des caract{\`e}res induits de Van Dijk
\cite{vdijk} montre  que l'application qui {\`a} un {\'e}l{\'e}ment $x\in R(G_d)$
associe la restriction de son caract{\`e}re-distribution aux {\'e}l{\'e}ments
elliptiques r{\'e}guliers se factorise {\`a} travers le quotient
$\o{R}(G_d)$. On d{\'e}duit alors du th{\'e}or{\`e}me \ref{JL1} le

\begin{coro} \label{JL2}
  L'application $R(\dd)\To{JL_d} R(G_d)$ induit un isomorphisme isom{\'e}trique
  $$R(\dd)\simto \o{R}(G_d)$$ caract{\'e}ris{\'e} par l'{\'e}galit{\'e} de caract{\`e}res du
  th{\'e}or{\`e}me \ref{JL1}.
%  correspondance de Jacquet-Langlands sur les classes de conjugaison
%  elliptiques r{\'e}guli{\`e}res se correspondant.
\end{coro}

En cons{\'e}quence, on a un morphisme dans l'autre sens $R(G_d)\To{}
\o{R}(G_d) \simto R(\dd)$ que nous noterons $LJ_d$. Cette notation est d\^ue \`a I. Badulescu qui a d\'efini dans \cite{Badu2} des r\'etractions similaires dans des situations plus g\'en\'erales. D'apr\`es \cite{Badu2}, ces r\'etractions n'envoient g\'en\'eralement pas irr\'eductibles sur irr\'eductibles, m\^eme au signe pr\`es. Dans le cas pr\'esent, le lemme suivant montre que l'image d'une irr\'eductible est soit nulle, soit une irr\'eductible au signe pr\`es.

\begin{lemme} \label{defell}
Pour une repr{\'e}sentation irr{\'e}ductible $\pi$ de $G_d$, les
  propri{\'e}t{\'e}s suivantes sont {\'e}quivalentes :
  \begin{enumerate}
  \item $\pi$ a le m{\^e}me support cuspidal qu'une s{\'e}rie discr{\`e}te.
  \item $\pi$ a un caract{\`e}re non nul sur les {\'e}l{\'e}ments elliptiques
    semi-simples r{\'e}guliers.
  \item L'image de $\pi$ dans $\o{R}(G_d)$ est non nulle.
  \end{enumerate}
Une repr{\'e}sentation satisfaisant ces propri{\'e}t{\'e}s 
sera dite {\em elliptique}. Son image par $JL_d$ co\"incide au signe pr\`es avec celle de la s\'erie discr\`ete de m\^eme support cuspidal. 
 \end{lemme}

Rappelons que le support cuspidal d'une repr{\'e}sentation
irr{\'e}ductible $\pi$ est l'unique classe de conjugaison de couples
$(M,\tau)$ form{\'e}s d'un \levi $M$ et d'une repr{\'e}sentation cuspidale
irr{\'e}ductible $\tau$ de $M$ qui appara{\^\i}t dans le module de Jacquet
normalis{\'e} de  $\pi$ le long d'un parabolique de Levi $M$.

\medskip

\begin{proof}
L'implication $ii) \Rightarrow iii)$ est une cons{\'e}quence de la formule
de Van Dijk \cite{vdijk} qui montre que le caract{\`e}re d'une induite
parabolique en un {\'e}l{\'e}ment semi-simple r{\'e}gulier elliptique est nul.

Pour voir l'implication $iii) \Rightarrow i)$,
% notons $[M,\tau]$ le
%support cuspidal de $\pi$ et
{\'e}crivons  $[\pi]=\sum_{(M,\sigma)} [\ip{M}{G_d}{\sigma}]$ (somme
d'induites de repr{\'e}sentations essentiellement temp{\'e}r{\'e}es) comme
nous le permet la classification par le quotient de Langlands.
On peut supposer que les supports cuspidaux de chaque $\sigma$ sont
contenus dans celui de $\pi$.
Comme les temp{\'e}r{\'e}es sont soit induites, soit des s{\'e}ries
discr{\`e}tes, on voit que si $\pi\notin R_I(G_d)$,
alors il y a une s{\'e}rie 
discr{\`e}te de m{\^e}me support cuspidal que $\pi$. Donc $iii)
\Rightarrow i)$.

Pour l'implication $i)\Rightarrow ii)$, on peut utiliser la
combinatoire de la classification de Zelevinski. Soit $\pi^{disc}$
l'unique s{\'e}rie discr{\`e}te de m{\^e}me support cuspidal que
$\pi$. Avec les notations de \cite{Zel}, on sait qu'il existe 
un (unique) segment 
$\Delta=[\tau,\tau']$ o{\`u} $\tau$ est une cuspidale irr{\'e}ductible de
$G_{d'}$ avec $d'$ diviseur de $d$, de sorte que 
\begin{enumerate}
\item L'ensemble partiellement ordonn{\'e} des multisegments de support
  $[\tau,\tau']$ a pour plus petit {\'e}l{\'e}ment $a_{min}=\{\Delta\}$ et
  plus grand {\'e}l{\'e}ment $a_{max}=\{\{\tau'\},\cdots, \{\tau\}\}$.
\item $\pi^{disc}=\la a_{max} \ra$ et toutes les repr{\'e}sentations de
  m{\^e}me support cuspidal sont de la forme $\la a \ra$ pour un
  multisegment $a$ de support $[\tau,\tau']$.
\item La repr{\'e}sentation irr{\'e}ductible $\la b \ra$ associ{\'e}e {\`a} $b$ appara{\^\i}t comme
  sous-quotient de la repr{\'e}sentation induite $\pi(a)$ associ{\'e}e {\`a} $a$
  \ssi\ $b\leq a$, et dans ce cas sa multiplicit{\'e} est $1$. 
\end{enumerate}
%Les repr{\'e}sentations de m{\^e}me
%support cuspidal que $\pi$ sont de la forme $Z(m)$ o{\`u} $m$ est un
%multisegment de support $\Delta$. Notons $I(m)$ l'induite associ{\'e}e. On
%sait que ces induites sont de multiplicit{\'e} $1$ et par le
%th{\'e}or{\`e}me \cite[7.1]{Zel} que $Z(m')$ appara{\^\i}t comme sous-quotient de
%$I(m)$ \ssi\ $m'\leq m$. 
Supposons maintenant que $\pi=\la a \ra$.
Pour $b< a$ notons $d(b,a)$ la longueur
d'une cha{\^\i}ne maximale $b<m_1<\cdots< a$ ({\em i.e} le nombre
d'op{\'e}rations ``{\'e}l{\'e}mentaires'' au sens de \cite[7]{Zel} pour passer de
$a$ {\`a} $b$), 
on obtient donc l'{\'e}galit{\'e} dans $R(G_d)$
$$ [\la a \ra] = \sum_{b\leq a} (-1)^{d(b,a)} [\pi(b)] $$
%Or l'ensemble des multisegments de support $\Delta$ a un plus grand
%{\'e}l{\'e}ment $m_{max}$ auquel correspond $Z(m_{max})=<\Delta>=\pi^{disc}$
%et un plus petit {\'e}l{\'e}ment  $m_{min}$ pour
Or, pour $b=a_{min}$, on a $\la a_{min}\ra=\pi(a_{min})$ (c'est la
repr{\'e}sentation de Speh associ{\'e}e {\`a} $\Delta$ et c'est l'image de la
s{\'e}rie discr{\`e}te $\pi^{disc}$ par l'involution de 
Zelevinski) et  pour $a_{min}\neq b$, la repr{\'e}sentation $\pi(b)$ est une
induite ``propre''. 
On obtient donc la congruence $[\pi]= \pm [\la  a_{min}\ra]
\hbox{ mod } R_I(G_d)$ dans $R(G_d)$. En l'appliquant aussi {\`a}
$\pi^{disc}$, on obtient $$[\pi]=\pm [\pi^{disc}] \hbox{ mod }
R_I(G_d).$$ 
La formule
de Van Dijk montre alors que les caract{\`e}res de $\pi$ et $\pi^{disc}$
 co{\"\i}ncident au signe pr{\`e}s sur les {\'e}l{\'e}ments elliptiques. Or, les
formules d'orthogonalit{\'e} pour les s{\'e}ries discr{\`e}tes montrent que le
caract{\`e}re d'une s{\'e}rie discr{\`e}te sur ces {\'e}l{\'e}ments est non nul.

% o{\`u} l'on a not{\'e} $\pi'$
%l'unique s{\'e}rie discr{\`e}te de m{\^e}me support cuspidal que
%$\pi$. N{\'e}anmoins nous ne donnons pas l'argument complet ici et
%donnerons un peu plus loin un argument utilisant les {\'e}quivalences de
%cat{\'e}gories de Bushnell-Kutzko.

%notons
%$R(G_d)_{(M,\tau)}$ le sous-groupe engendr{\'e} par les irr{\'e}ductibles de $G_d$
%dont le support cuspidal contient la paire $(M,\tau)$ o{\`u} $M$ est un
%\levi de $G_d$ et $\tau$ une repr{\'e}sentation irr{\'e}ductible cuspidale de
%$M$. 
%D'apr{\`e}s \cite{TPW}, 
%la d{\'e}compositon  $R(G_d)=\bigoplus_{(M,\tau)}R(G_d)_{(M,\tau)}$ selon le
%support cuspidal induit une d{\'e}composition
%$R_I(G_d)=\bigoplus_{(M,\tau)} (R_I(G_d)\cap R(G_d)_{(M,\tau)})$.
%En particulier, l'application

\end{proof}

Nous noterons $\Ell{C}{G_d} \subset
\Irr{C}{G_d}$ le sous-ensemble des classes de repr{\'e}sentations
elliptiques. On a bien-s{\^u}r
$$ \Cu{C}{G_d} \subset \Disc{C}{G_d} \subset \Ell{C}{G_d}. $$
Nous allons donner une classification de ces repr{\'e}sentations adapt{\'e}e
aux besoins de ce texte. Auparavant, nous devons faire quelques rappels de th\'eorie des repr\'esentations de groupes lin\'eaires $p$-adiques.

%\noindent{\em Notations} :
\alin{Rappels sur la th{\'e}orie de Bernstein \cite{bernstein}} \label{decbernstein} 
Si $G$ est un groupe r\'eductif $p$-adique, $\Mo{C}{G}$ d\'esigne la
cat{\'e}gorie ab{\'e}lienne de toutes les 
$C$-repr{\'e}sentations lisses de $G$.
Soit $(M,\tau)$ une paire Levi-cuspidale, d\'efinissons
%form{\'e}e d'un \levi de $G_d$ et d'une
%$C$-repr{\'e}sentation cuspidale irr{\'e}ductible de ce Levi et 
$\BG_{M,\tau}^{G}$
la sous-cat{\'e}gorie pleine de $\Mo{C}{G}$ form{\'e}e des objets dont tous
les sous-quotients irr{\'e}ductibles contiennent $(M,\tau\psi)$ dans leur
support cuspidal, pour un certain caract{\`e}re non-ramifi{\'e} $\psi$ de $M$. On
sait que la cat{\'e}gorie $\BG_{M,\tau}^{G}$ est un ``facteur direct
ind{\'e}composable'' (que nous appellerons bloc de Bernstein associ{\'e} {\`a}
$(M,\tau)$) de
$\Mo{C}{G}$ et qu'on a une d{\'e}composition
$$ \Mo{C}{G} \simeq \bigoplus_{(M,\tau)/\sim} \BG_{M,\tau}^{G} $$ 
o{\`u} $(M,\tau)\sim (M',\tau')$ \ssi\ il existe $g\in G$ et $\psi$
caract{\`e}re non ramifi{\'e} de $M$ tels que $M'=M^g$ et $\tau'=(\tau\psi)^g$
({\'e}quivalence ``inertielle''). En particulier, les idempotents centraux
primitifs de $\Mo{C}{G}$ sont en bijection avec les classes
inertielles de paires $(M,\tau)$.

Pour un produit de groupes lin\'eaires, lorsque $M = T$ est un tore
maximal et $\tau$ est un caract{\`e}re 
non-ramifi{\'e} de $T$, %$\UG_d=\BG_{M,\tau}^{G_d}$ et 
on appelle parfois le bloc associ\'e $\BG^{G}_{T,1}$
le ``bloc unipotent'' de $G$.

\alin{Rappels sur la th\'eorie de Bushnell-Kutzko}
Cette th\'eorie permet de d\'ecrire les blocs de Bernstein dans le cas o\`u $G$ est un groupe lin\'eaire. On s'int\'eresse ici au cas particulier suivant : on fixe un entier $n$ et
 une repr\'esentation supercuspidale irr\'eductible $\tau$ de $G_{n}=GL_{n}(K)$.
 Pout tout entier  $e\geq 1$, on voit le groupe produit  $(G_{n})^{e}$ 
 comme un Levi standard (diagonal) dans $G_{ne}$, d'o\`u une paire Levi-cuspidale $((G_n)^e,\tau^e)$ pour $(G_{n})^{e}$.
Si $K'$ d\'esigne un autre corps local, nous noterons aussi $G'_n:=GL_{n}(K')$.

\begin{fact} \label{equivcat}
La th\'eorie de Bushnell-Kutzko nous fournit une extension $K'$ de $K$, de degr\'e $n$ et degr\'e r\'esiduel le nombre $f_\tau$ de caract\`eres non ramifi\'es $\chi: K^\times \To{} C^\times$ tels que $(\chi\circ \hbox{det})\otimes\tau \simeq \tau$, et une famille d'\'equivalences de cat\'egories 
  $$ \alpha_\tau^{e} :\; \BG_{(G_n)^e,\tau^{e}}^{G_{ne}} \simto \BG_{(G'_1)^e,1^{e}}^{G'_e}$$
  % o\`u %$G'_e$ d\'esigne le tore des matrices diagonales de 
  %$G'_e:= GL_e(K_\tau)$, 
  satisfaisant les propri\'et\'es suivantes :
 \begin{enumerate}
\item Normalisation : $\alpha_\tau^1(\tau) = 1$.
        \item Compatibilit\'e \`a l'induction parabolique normalis\'ee\footnote{On rappelle que les paraboliques standards sont les paraboliques triangulaires sup\'erieurs,} : si $e=\sum_{i=1}^r e_i$ est une partition de $e$ alors 
        $$ \alpha_\tau^{e_1}\times \cdots \times \alpha_{\tau}^{e_r} \simeq \alpha_\tau^e \circ \Ip{(G_n)^{e_1} \times \cdots \times (G_n)^{e_r}}{G_{ne}} $$ 
   \item  Compatibilit\'e \`a la torsion : si $\chi : K^\times \To{} C^\times$ est un caract\`ere, alors pour toute $\pi \in \BG^{G_{ne}}_{(G_n)^e,\tau^{e}}$, on a
  $$\alpha_\tau^e \left((\chi\circ\hbox{det})\otimes \pi\right) \simeq (\chi\circ N_{K'|K}\circ\hbox{det})\otimes\alpha_\tau^e(\pi). $$
\item Compatibilit\'e aux caract\`eres centraux : pour toute repr\'esentation $\pi\in \BG^{G_{ne}}_{(G_n)^e,\tau^{e}}$ telle que $\pi(\varpi_K)$ soit un scalaire, $\alpha_\tau^e(\pi)(\varpi_K)$ est aussi un scalaire et on a $\alpha_\tau^e(\pi)(\varpi_K)=\pi(\varpi_K) \tau(\varpi_K)^{-1}$.

\end{enumerate}
\end{fact}

\begin{proof}
On utilise les notations de \cite[Thm (7.6.20)]{BK}. L'\'equivalence $\alpha_\tau^1$ est donn\'ee par ce qui est not\'e ${\mathfrak A}{\mathfrak d}(\Psi_1)^{-1}$ dans {\em loc. cit}, le choix de $\Psi_1$ \'etant fix\'e par notre propri\'et\'e i).
 Pour chaque $e>1$ on d\'efinit alors $\alpha_\tau^e:=  {\mathfrak A}{\mathfrak d}(\Psi)^{-1}$, 
  donn\'e par \cite[Cor (7.6.21)]{BK}. En fait, pour \^etre un peu plus pr\'ecis,
  l'\'equivalence de cat\'egories donn\'ee par {\em loc.cit} ne concerne que les
  repr\'esentations admissibles de longueur finie. Pour obtenir
  l'\'equivalence des blocs ``en entier'', il faut utiliser aussi
  \cite{BK1}.   La compatibilit\'e \`a la torsion par les caract\`eres se d\'eduit de \cite[(7.5.12)]{BK}. La  propri\'et\'e iv) est mentionn\'ee en \cite[(7.7.6)]{BK} pour les repr\'esentations admissibles. Elle se g\'en\'eralise aux autres repr\'esentations, en utilisant la discussion de \cite[(7.5.9)]{BK} par exemple.
La compatibilit\'e \`a l'induction parabolique est explicite 
  \cite[7.6.21]{BK} dans le cas ``minimal" $e=1+\cdots +1$. Dans le cas g\'en\'eral, elle est donn\'ee par  \cite{BK1}.
  \end{proof}

\begin{nota} \label{notationsrho}
  Soit $\rho\in \Irr{C}{\dd}$. Nous noterons
\begin{itemize}
%\item $(M_\rho,\tau_\rho)$ le support cuspidal de $JL_d(\rho)$ : c'est
%  une paire form{\'e}e d'un \levi de $G_d$ et d'une
%$C$-repr{\'e}sentation cuspidale irr{\'e}ductible de ce Levi (en bref, une
%paire ``Levi-cuspidale'') et qui est  d{\'e}finie {\`a}
%  $G_d$-conjugaison pr{\`e}s.
\item $d_\rho\in \NM^*$ l'unique diviseur de $d$ et $\pi_\rho$
  l'unique repr\'esentation cuspidale irr\'eductible de $G_{d/d_\rho}$
  tels que $JL_d(\rho)$ apparaisse dans l'induite parabolique standard
  normalis\'ee
  %$$\ip{(G_{d/d_\rho})^{d_\rho}}{G_d}{\pi_\rho^{d_\rho}} =
  $ {|det|^{\frac{1-d_\rho}{2}}\pi_\rho \times\cdots\times |det|^{\frac{d_\rho-1}{2}} \pi_\rho}.$ L'existence de
  $d_\rho$ et $\pi_\rho$ est assur{\'e}e par la classification des
  s{\'e}ries discr{\`e}tes de $G_d$ par Bernstein-Zelevinski,
  \cite[Thm 9.3]{Zel}.
\item $M_\rho$ le \levi\ standard $(G_{d/d_\rho})^{d_\rho}$ et
  $\vec\pi_\rho$ la repr\'esentation supercuspidale irr\'eductible $
  \pi_\rho |det|^{(1-d_\rho)/2} \otimes \cdots \otimes \pi_\rho
  |det|^{(d_\rho-1)/2}$ de $M_\rho$.  Ainsi la paire $(M_\rho,\vec\pi_\rho)$
  est un repr\'esentant du support cuspidal de $JL_d(\rho)$
\end{itemize}
\end{nota}

%Le fait que l'image de $JL_d(\rho)$ est la Steinberg est une cons\'equence de 
%  \cite[Thm 7.7.1]{BK} et du fait que le support cuspidal de cette image contient $(T'_{d_\rho},1)$ par d\'efinition de $\alpha_{\pi_\rho}$ et compatibilit\'e aux foncteurs paraboliques. 

D'apr\`es la caract\'erisation \ref{defell} iii) des repr\'esentations elliptiques et la propri\'et\'e \ref{equivcat} ii), 
l'\'equivalence de cat\'egories $\alpha_{\pi_\rho}^{d_\rho}$ induit une bijection entre l'ensemble
$\Ell{C}{\BG_{M_\rho,\pi_\rho}^{G_d}}$ des repr{\'e}sentations
elliptiques de $G_d$ dans le bloc  $\BG_{M_\rho,\vec\pi_\rho}^{G_d}$ et l'ensemble
$\Ell{C}{\BG^{G'_{d_\rho}}_{T_{d_\rho},1}}$ des repr{\'e}sentations elliptiques de
$G'_{d_\rho}$ dans son bloc unipotent. Notons que toute  famille d'\'equivalences de cat\'egories \ref{equivcat} satisfaisant les propri\'et\'es i) \`a iv) induit la m\^eme bijection. Par ailleurs, $\pi_\rho$ a \'et\'e choisie pour que cette bijection envoie $JL_d(\rho)$ sur la repr\'esentation de Steinberg de $G'_{d_\rho}$.

\alin{Classification des repr{\'e}sentations elliptiques} \label{classrepell}
%Nous utiliserons la classification de Rodier plut\^ot que celle de Zelevinski.
 Rappelons tout d'abord la classification des elliptiques du bloc unipotent de $G_d$ que nous avons utilis\'ee dans \cite[2.1.3]{Dat1}. Notons $S_d$ l'ensemble des racines simples du tore diagonal de $G_d$ dans le Borel sup\'erieur. Tout sous-ensemble $I\subseteq S_d$ d\'etermine un sous-groupe de Levi diagonal $M_I$ et un parabolique $P_I$ triangulaire sup\'erieur, dont l'induite droite $\ind{P_I}{G_d}{1}$ poss\`ede un unique quotient irr\'eductible que nous noterons $\pi_1^I$. Celui-ci est elliptique, et  l'application 
 $$ I\subseteq S_d \mapsto \pi_1^I \in \Ell{C}{\BG^{G_d}_{(G_1)^d,1^d}}$$ est une bijection.
En termes d'induites normalis\'ees, on a $\pi_1^I=   \hbox{Cosoc}\left(\Ip{M_{I}}{G_d}\left(\hbox{Soc}\left(\Ip{T_d}{M_{I}}\left(\vec\pi_1\right)\right)\right) \right)
$.

Pour $\rho \in \Irr{C}{\dd}$, on introduit l'ensemble $S_\rho$ des racines simples du centre de $M_\rho$ dans l'alg\`ebre de Lie du parabolique triangulaire sup\'erieur $P_\rho$ de Levi $M_\rho$, et pour $I\subseteq S_\rho$, on note $M_{\rho,I}$ le %parabolique contenant $P_\rho$ et dont la composante de Levi est le 
centralisateur du noyau commun des $\alpha\in I$. La repr\'esentation 
\ini
\begin{equation}\label{notaell}
 \pi^I_\rho := \hbox{Cosoc}\left(\Ip{M_{\rho,I}}{G_d}\left(\hbox{Soc}\left(\Ip{M_{\rho}}{M_{\rho,I}}\left(\vec\pi_\rho\right)\right)\right) \right)
\end{equation}
de $G_d$ est irr\'eductible et elliptique, par compatibilit\'e de l'\'equivalence $\alpha_{\pi_\rho}^{d_\rho}$ aux inductions et torsions, et l'application
 $$ I\subseteq S_\rho \mapsto \pi^I_\rho \in \Ell{C}{\BG^{G_d}_{M_\rho,\pi_\rho^{d_\rho}}}$$ est une bijection.
Notons que  $JL_d(\rho)=\pi_\rho^\emptyset$. 

%On a donc obtenu la bijection $(\rho \in \Irr{C}{\dd},I\subseteq S_{\rho}) \mapsto
%\pi_\rho^I\in \hbox{Ell}_C(G_d)$ qui est la classification annonc\'ee des repr{\'e}sentations elliptiques de $G_d$. 

\begin{conv} \label{conv}
Bien que cela paraisse moins intrins\`eque, nous devons num\'eroter $S_\rho$, {\em i.e.} choisir une bijection $S_\rho\simto \{1,\cdots,d_\rho-1\}$ pour la suite. Nous choisissons de num\'eroter les racines simples des paraboliques triangulaires sup\'erieurs de haut en bas. 
\end{conv}
En d\'ecorant de signes $'$ les objets relatifs au corps $K'$ de l'\'equivalence $\alpha_{\pi_\rho}^{d_\rho}$ de \ref{equivcat}, 
la convention ci-dessus permet en particulier d'identifier $S_\rho$ et $S'_{d_\rho}$, et on a pour tout $I\subseteq S_\rho$ 
$$ \alpha_{\pi_\rho}^{d_\rho}(\pi_\rho^I)\simeq {\pi'}_1^I.  $$

\begin{rema}
Pour tous $\rho\in \Irr{C}{\dd}$ et $I\subseteq S_\rho$, on a dans $\o{R}(G_d)$  l'{\'e}galit{\'e}
$[\pi_\rho^I]=(-1)^{|I|}[\pi_\rho^\emptyset]$. Par cons{\'e}quent on a
aussi
$LJ_d[\pi^I_\rho] = (-1)^{|I|} 
[\rho]$ dans $R(\dd)$.
\end{rema}
\begin{proof}
Gr{\^a}ce aux \'equivalences de cat\'egories \ref{equivcat} on est ramen{\'e} au bloc unipotent, c'est-\`a dire au cas $\rho=1$.
Dans ce cas, le lemme X.4.6 de \cite{BW} montre que
%on remarque que le
%complexe acyclique de \cite[Prop 6.13]{SS1} ou \cite[{BW}] montre que
dans $R(G_{d})$, on a
$$ [\pi^{I}_1] = \sum_{I\subseteq J \subseteq S_{d}} (-1)^{|J\setminus I|}  
   [\ind{P_I}{G_{d}}{1}]. $$
En appliquant ceci {\`a} $I$ et $\emptyset$, on obtient $[\pi_1^{I}]
= (-1)^{d-1-|I|}[1_{G_{d}}] = 
(-1)^{|I|}[\hbox{St}_{G_{d}}]$.

\end{proof}

\begin{rema} \label{dualiteell}
La duale de $\pi_\rho^I$ est $\pi_{\rho^\vee}^{\o{I}}$, o\`u l'on identifie $S_{\rho}=S_{{\rho^\vee}}$ et $\o{I}$ est l'image de $I$ par l'involution $i\mapsto d_\rho-i$ de $S_{\rho}$ num\'erot\'e comme en \ref{conv}.
\end{rema}

\begin{proof}
Par compatibilit\'e de la correspondance de Jacquet-Langlands avec le passage \`a la contragr\'ediente, on a $JL_d(\rho^\vee)=JL_d(\rho)^\vee$, et donc
$M_{\rho^\vee}=M_\rho$ et $\vec\tau_{\rho^\vee}\simeq (\vec\pi_\rho)^\vee$. 
En dualisant la d\'efinition de $\pi^I_\rho$, on obtient
$$ (\pi^I_\rho)^\vee = \hbox{Soc}\left(\Ip{M_{\rho,I}}{G_d}\left(\hbox{Cosoc}\left(\Ip{M_{\rho}}{M_{\rho,I}}\left((\vec\pi_\rho)^\vee\right)\right)\right) \right).
$$
Appliquant l'\'equivalence $\alpha_{\pi_{\rho^\vee}}^{d_\rho}$ de \ref{equivcat} et  identifiant $S_\rho$ \`a  $S'_{d_\rho}$ gr\^ace \`a \ref{conv}, on obtient
$\alpha_{\tau_{\rho^\vee}}^{d_\rho}\left((\pi^I_\rho)^\vee\right) = ({\pi'}_1^{I})^\vee $
qui d'apr\`es \cite[2.3.3 iv)]{Dat1} n'est autre que ${\pi'}_1^{\o{I}}$ (le $'$ renvoie au corps $K'$ de \ref{equivcat}).
\end{proof}

\alin{Extensions entre repr{\'e}sentations elliptiques} 
Nous identifierons les centres de $G_d$ et $\dd$ \`a $K^\times$ par les plongements canoniques. En particulier le caract\`ere central d'une repr\'esentation de $G_d$ ou $\dd$ sera vu comme un caract\`ere de $K^\times$.
Avec cette convention les caract\`eres centraux de $\rho\in\Irr{\o\QM_l}{\dd}$ et $JL_d(\rho)\in
\Irr{\o\QM_l}{G_d}$ sont \'egaux. Dans la proposition suivante, la notation $\hbox{Ext}_{G_d,\omega}$ d\'esigne les groupes d'extensions calcul\'es dans la cat\'egorie des repr\'esentations lisses de $G_d$ de caract\`ere central $\omega$.

\begin{prop} \label{theoextell}
Soient $\rho \in \Irr{\o\QM_l}{\dd}$ de caract\`ere central $\omega$ et $I\subseteq S_{\rho}$.
\begin{enumerate}
  \item Pour tout $i\in\NM$ et toute $\pi\in \Irr{\omega}{G_d}$, on a en notant $\delta(I,I')=|I\cup I'|-|I\cap I'|$
    $$  \cas{\ext{i}{\pi^I_{\rho}}{\pi}{G_d,\omega}}{\o\QM_l}{\pi\simeq \pi_\rho^{I'} \hbox{ et }i=\delta(I,I')}{0}{\hbox{le contraire}}. $$
\item Soient $I,J,K$ trois sous-ensembles de $S_{\rho}$ tels que
  $\delta(I,J)+\delta(J,K)=\delta(I,K)$, alors le cup-produit
$$  \cup \;:\;\;
\ext{\delta(I,J)}{\pi^I_{\rho}}{\pi^J_{\rho}}{G_d,\omega} \otimes_{\o\QM_l}
\ext{\delta(J,K)}{\pi^J_{\rho}}{\pi^K_\rho}{G_d,\omega} \To{} 
\ext{\delta(I,K)}{\pi^I_{\rho}}{\pi^K_\rho}{G_d,\omega}
$$ est un isomorphisme.
\end{enumerate}
\end{prop}

\begin{proof}
 Rappelons  que si $\pi,\pi'$ sont deux repr\'esentations
irr\'eductibles de m\^eme caract\`ere central $\omega$ mais de supports
cuspidaux distincts, alors $\ext{i}{\pi}{\pi'}{G,\omega}=0$ pour tout
$i\in\NM$, voir par exemple \cite[6.1]{Vigext}. Ceci montre la nullit\'e des Ext lorsque $\pi$ n'est pas elliptique de type $\rho$.

Pour le reste du th\'eor\`eme, nous allons nous ramener au
th\'eor\`eme 1.3 de \cite{Dat1} gr\^ace \`a l'\'equivalence
$\alpha_{\pi_\rho}^{d_\rho}$. Pour cela,
remarquons que par les propri\'et\'es i) et iv) de \ref{equivcat},
 celle-ci induit une \'equivalence entre la sous-cat\'egorie pleine des objets de $\BG^{G_d}_{M_\rho,\vec\pi_\rho}$ o\`u $\varpi$ agit par le scalaire $\omega(\varpi)$ et la sous-cat\'egorie pleine des objets
de $\BG^{G'_{d_\rho}}_{(G'_1)^{d_\rho},\vec{1}}$ o\`u $\varpi$ agit trivialement. % (rappelons que $G'_{d_\rho}=GL_{d_\rho}(K_\rho)$ avec les notations de \ref{equivcat}). 
Comme la premi\`ere contient la sous-cat\'egorie pleine des objets de caract\`ere central $\omega$ comme {\em facteur direct} (par compacit\'e de $\OC_K^\times$) et la seconde contient la sous-cat\'egorie des objets de caract\`ere central trivial comme facteur direct aussi (par compacit\'e de $\OC_{K'}^\times$), et que par cons\'equent de chaque c\^ot\'e les Ext sont les m\^emes dans l'une ou l'autre sous-cat\'egorie pleines en question,
il s'ensuit que 
$$\ext{*}{\pi^I_\rho}{\pi^J_\rho}{G_d,\omega} = \ext{*}{{\pi'}^{I}_{1}}{{\pi'}^{J}_{1}}{PG'_{d_\rho}}$$
o\`u l'on identifie $S_\rho$ et $S'_{d_\rho}$ par \ref{conv} et o\`u le $'$ renvoie au corps $K'$ de \ref{equivcat}. Il ne reste plus qu'\`a appliquer \cite[1.3]{Dat1}
%Mais puisque ${K'}^\times/\varpi^\ZM$ est compact, les objets projectifs de $\Mo{\o\QM_l}{PG'_{d_\rho}}$ sont encore projectifs lorsqu'ils sont vus comme objets de $\Mo{\o\QM_l}{G'_{d_\rho}/\varpi^\ZM}$ de sorte que 
%$$\ext{*}{{\pi'}^{I}_{1}}{{\pi'}^{J}_{1}}{G'_{d_\rho}/\varpi^\ZM}
%=\ext{*}{{\pi'}^{I}_{1}}{{\pi'}^{J}_{1}}{PG'_{d_\rho}}.$$
%On est donc ramen\'e au th\'eor\`eme 1.3 de \cite{Dat1}.
 \end{proof}

\begin{rema} \label{remaextell}
Comme le montre la preuve ci-dessus, on peut dans la proposition pr\'ec\'edente remplacer les $\hbox{Ext}_{G_d,\omega}$ par les $\hbox{Ext}$ calcul\'es dans la sous-cat\'egorie pleine des objets de $\Mo{\o\QM_l}{G_d}$ o\`u $\varpi$ agit par le scalaire $\omega(\varpi)$.
\end{rema}

\subsection{Correspondance de Langlands locale} \label{corLan}
\def\rep{\hbox{Rep}}

Soit $K^{nr}$ la sous-extension non ramifi{\'e}e maximale de $K$ dans
$K^{ca}$ et $\wh{K^{nr}}$ sa compl{\'e}tion. On note toujours $I_K :=
\gal(K^{ca}/K^{nr})\subset W_K$ le groupe d'inertie de $K$.

\alin{Formulation {\`a} la Weil-Deligne}
Rappelons qu'une  ``représentation de Weil-Deligne" $\sigma$ de $W_K$ à valeurs dans le corps $C$ est un triplet $(\sigma^{ss},N_\sigma,V_\sigma)$ où 
\begin{itemize}
\item $V_\sigma$ est un espace vectoriel de dimension finie sur $C$,  
\item $\sigma^{ss}: W_K\To{} GL(V_\sigma)$ est une représentation continue (pour la topologie discrète de $C$) et semi-simple de $W_K$.
\item $N_\sigma\in \endo{C}{V_\sigma}$ est un endomorphisme nilpotent de $V$ tel que pour tout $w\in W_K$, on a 
\ini\begin{equation}\label{monod}
\sigma^{ss}(w)N_\sigma\sigma^{ss}(w)^{-1} = |w|N_\sigma,
\end{equation}
o\`u $|-|$ d\'esigne le caract\`ere non ramifi\'e de $W_K$ qui envoie les Frobenius arithm\'etiques sur l'ordre $q$ du corps r\'esiduel de $K$.
\end{itemize}
%Notons $WD_K:=W_K\ltimes \GM_a$ le groupe de Weil-Deligne de $W_K$ et
%$\rep^d_C(WD_K)$ l'ensemble des classes d'{\'e}quivalences de
%$C$-repr{\'e}sentations $\sigma$  de $WD_K$ de dimension $d$ {\em
%  admissibles}, {\em i.e.} v{\'e}rifiant
%\begin{itemize}
%\item $\sigma(\phi)$ est semi-simple.
%\item $\sigma_{|\GM_a}$ est algebrique.
%\item $\sigma_{|I_K}$ a un noyau d'indice fini.
%\end{itemize}
%$
Notons $\rep^d_C(WD_K)$ l'ensemble des classes d'{\'e}quivalences des
 $C$-repr{\'e}sentations de Weil-Deligne de dimension $d$.
La correspondance de Langlands sur $K$ est une famille de bijections
$(\sigma_d)_{d\in\NM^*}$ 
$$\sigma_d :\;\;\Irr{\CM}{GL_d(K)} \simto \rep^d_\CM(WD_K)$$ 
 qui v{\'e}rifie un certain nombre de propri{\'e}t{\'e}s
suffisantes pour la rendre unique
(compatibilit{\'e} avec la th{\'e}orie du
corps de classe qui donne aussi le cas $d=1$,  pr{\'e}servation de certains
invariants de nature arithm{\'e}tico-analytique (facteurs $L$ et
$\epsilon$ de paires), etc...). Nous renvoyons {\`a} \cite{HeLang},\cite{HeSMF} pour l'{\'e}nonc{\'e}
pr{\'e}cis de ces propri{\'e}t{\'e}s caract{\'e}ristiques que nous n'utiliserons pas
en totalit{\'e}. Nous utiliserons sans commentaires particuliers la compatibilité à la contragrédiente et la compatibilité à la torsion par les caractères qui dans le cas des caractéres non ramifiés s'exprime formellement ainsi :
$$\forall \pi\in\Irr{\CM}{G_d},\forall a\in\CM^\times,\;\; \sigma_d((|\hbox{det}|^a\pi)=|-|^a \sigma_d(\pi) $$
(en particulier, on a normalisé le corps de classes de manière à ce que les uniformisantes correspondent aux Frobenius géométriques).
Rappelons par ailleurs que 
$$ \sigma_d(\Cu{\CM}{G_d}) = \{\sigma \in \rep^d_\CM(WD_K) \; \hbox{
  irr{\'e}ductibles}\}. $$
En fait, la correspondance est  d{\'e}termin{\'e}e par sa restriction aux
repr{\'e}sentations cuspidales de $G_d$ et irr{\'e}ductibles de $W_K$, la
classification de Zelevinski \cite{Zel} et la classification des
repr{\'e}sentations de Weil-Deligne en fonction des irr{\'e}ductibles de
$W_K$.

La g{\'e}om{\'e}trie alg\'ebrique r{\'e}alise plut\^ot une variante de la
correspondance de Langlands, appel{\'e}e parfois correspondance de Hecke
et qui se d{\'e}duit de celle de Langlands par une simple torsion ``{\`a} la
Tate'' par le caract{\`e}re $w \mapsto |w|^{{1-d}\over 2}$. Il est sous-entendu ici que l'on choisit toujours (si besoin) la racine carrée positive de $q$ dans $\CM$ pour définir une puissance demi-entière de $|-|$.
Pour les représentations cuspidales, cette variante se trouve {\^e}tre compatible 
%Comme nous avons besoin de changer de corps de coefficients
%(alg{\'e}briquement clos de caract{\'e}ristique nulle), mentionnons que ces
%bijections sont compatibles 
aux actions des automorphismes de $\CM$ des deux c{\^o}t{\'e}s, cf \cite[(7.4)]{HeTNB}. 
L'extension aux  représentations non cuspidales respecte cette compatibilit\'e, 
comme nous allons le 
%cela ne semble plus vrai.
préciser maintenant en suivant \cite[(7.4)]{HeTNB}.
Si $\tau$ est un automorphisme de $\CM$ et $\sigma=(\sigma^{ss},N,V)$ une représentation de Weil-Deligne, notons $\sigma^\tau:=(\sigma^{ss}\otimes 1, N\otimes 1, V\otimes_{\CM,\tau}\CM)$.
Cela définit une action de $\aut{}{\CM}$ sur 
$\rep^d_\CM(WD_K)$ et de même on en définit une autre sur $\Irr{\CM}{G_d}$.

\begin{fact}
Pour tout automorphisme $\tau$ du corps $\CM$ et toute représentation $\pi\in\Irr{\CM}{G_d}$, on a
$$ |-|^{(1-d)/2}\sigma_d(\pi^\tau) = (|-|^{(1-d)/2}\sigma_d(\pi))^\tau .$$
\end{fact}
\begin{proof}
Le cas le plus difficile est celui o\`u $\pi$ est cuspidale ; il est prouvé dans \cite[(7.4)]{HeTNB}
% que pour tout automorphisme $\tau$ du corps $\CM$, on a 
%$$ |-|^\frac{1-d}{2}\otimes \sigma_d(\tau(\pi)) = \tau(|-|^\frac{1-d}{2}\otimes \sigma_d(\pi)) %.$$
(et cela repose sur des compatibilités avec des cas de correspondance globale).
Le cas g\'en\'eral d\'ecoule certainement de la version rationnelle de la classification de Langlands de \cite[Prop. 3.2]{CloAA}. Pour le confort du lecteur nous donnons une preuve compl\`ete.

Modifions un peu l'\'enonc\'e en introduisant
la notation $\varepsilon_\tau = \tau(\sqrt q)/\sqrt q \in \{\pm 1\}$ et les caract\`eres $$\varepsilon_\tau:\;\; g\in G_d \mapsto \varepsilon_\tau^{\hbox{val}\circ\hbox{det}(g)}  
\;\;\hbox{et}\;\; \varepsilon_\tau:\;\; w\in W_K \mapsto \varepsilon_\tau^{\hbox{log}_q|w|}.$$
L'\'enonc\'e que l'on veut prouver se reformule en
\ini\begin{equation}\label{refor}
 \sigma_d(\pi^\tau)=\varepsilon_\tau^{d-1}\sigma_d(\pi)^\tau. 
 \end{equation}
\'Etant donn\'ees $\pi_1,\cdots,\pi_r$ des repr\'esentations de $G_{d_1},\cdots,G_{d_r}$, notons $\pi_1\times \cdots \times \pi_r$ la repr\'esentation de $G_{d_1+\cdots +d_r}$ induite parabolique standard normalis\'ee. Alors si $d=d_1+\cdots +d_r$, on calcule 
$$ \left(\pi_1\times\cdots \times \pi_r\right)^\tau = \varepsilon_\tau^{d-d_1}\pi_1^\tau\times \cdots \times \varepsilon_\tau^{d-d_r}\pi_r^\tau. $$
Supposons maintenant que cette induite est irr\'eductible et que (\ref{refor}) est connu pour chaque $\pi_i$. On sait alors que
$$ \sigma_d(\pi_1\times\cdots \times \pi_r)=\sigma_{d_1}(\pi_1)\oplus\cdots \oplus\sigma_{d_r}(\pi_r) $$
et on en tire imm\'ediatement l'\'enonc\'e pour l'induite. D'apr\`es la classification de Zelevinski, cela nous ram\`ene au cas o\`u $\pi$ est l'unique quotient irr\'eductible de 
$$\hbox{St}_{k_1}(|det|^{n_1}\pi_g)\times \cdots \times \hbox{St}_{k_r}(|det|^{n_r}\pi_g)$$
avec $n_1\geq\cdots\geq n_r \in\ZM$ et $\pi_g\in\Cu{\CM}{G_g}$. On a not\'e ici 
$\hbox{St}_{k}{\pi_g}$ l'unique quotient irr\'eductible de $\pi_g|det|^{(1-k)/2}\times \cdots \times \pi_g|det|^{(k-1)/2}$ ($k$ facteurs, on passe au suivant en multipliant par $|det|$).
Remarquons que $(\hbox{St}_{k}{\pi_g})^\tau=\hbox{St}_{k}{(\varepsilon_\tau^{kg-g-k+1} \pi_g^\tau)}$.
Par ailleurs, d'apr\`es \cite[2.7]{HeSMF}, on a $\sigma_{kg}(\hbox{St}_{k}{\pi_g})=\sigma_g(\pi_g)\otimes\sigma_k(\hbox{St}_{G_k})$.
On calcule donc 
\begin{eqnarray*}
\sigma_{kg}((\hbox{St}_{k}{\pi_g})^\tau)=\sigma_{kg}(\hbox{St}_{k}{(\varepsilon_\tau^{kg-g-k+1} \pi_g^\tau)}) & = & \varepsilon_\tau^{kg-g-k+1}\sigma_g(\pi_g^\tau)\otimes \sigma_k(\hbox{St}_{G_k}) \\
& = & \varepsilon_\tau^{kg-g-k+1}\varepsilon_\tau^{1-g}\varepsilon_\tau^{1-k}
(\sigma_g(\pi_g)\otimes\sigma_k(\hbox{St}_{G_k}))^\tau \\
& = & \varepsilon_\tau^{1-kg}\sigma_{kg}(\hbox{St}_{k}{\pi_g})^\tau
\end{eqnarray*}
compte tenu de ce que l'on sait d\'eja pour $\pi_g$ et $\hbox{St}_{G_k}$. Ainsi (\ref{refor}) est v\'erifi\'e pour les repr\'esentations du type $\hbox{St}_k(\pi_g)$. 
Revenons \`a notre repr\'esentation irr\'eductible $\pi$ et rappelons que selon \cite[2.9]{HeSMF}, il lui est associ\'e la repr\'esentation galoisienne
$$ \sigma_d(\pi)=\sigma_{k_1g}(\hbox{St}_{k_1}(|det|^{n_1}\pi_g))\oplus \cdots \oplus \sigma_{k_rg}\hbox{St}_{k_r}(|det|^{n_r}\pi_g)). $$
Remarquons par ailleurs que $\pi^\tau$ est l'unique quotient irr\'eductible de la repr\'esentation
$$\varepsilon_\tau^{d-k_1g}(\hbox{St}_{k_1}(|det|^{n_1}\pi_g))^\tau\times \cdots \times \varepsilon_\tau^{d-k_rg}(\hbox{St}_{k_r}(|det|^{n_r}\pi_g))^\tau$$  
qui n'est autre que la repr\'esentation
$$\hbox{St}_{k_1}\left( \varepsilon_\tau^{d-k_1g}\varepsilon_\tau^{k_1g-g-k_1-1} |det|^{n_1}\pi_g^\tau\right)\times\cdots\times
\hbox{St}_{k_1}\left( \varepsilon_\tau^{d-k_1g}\varepsilon_\tau^{k_1g-g-k_1-1} |det|^{n_r}\pi_g^\tau\right).
$$
Il lui est donc associ\'e la repr\'esentation
$$\sigma_d(\pi^\tau)= \sigma_{k_1g}\left(  \hbox{St}_{k_1}(\varepsilon_\tau^{d-g-k_1-1}|det|^{n_1}\pi_g^\tau)\right)\oplus\cdots\oplus \sigma_{k_rg}\left(  \hbox{St}_{k_r}(\varepsilon_\tau^{d-g-k_r-1}|det|^{n_r}\pi_g^\tau)\right)
$$
qui n'est autre que
$$\varepsilon_\tau^{d-k_1g}\sigma_{k_1g}\left((\hbox{St}_{k_1}(|det|^{n_1}\pi_g))^\tau\right)\oplus \cdots \oplus \varepsilon_\tau^{d-k_rg}\sigma_{k_rg}\left((\hbox{St}_{k_r}(|det|^{n_r}\pi_g))^\tau\right).$$  
Ainsi, par le cas d\'eja trait\'e des repr\'esentations $\hbox{St}$, on en d\'eduit (\ref{refor}) pour $\pi$.

\end{proof}

La compatibilit{\'e} aux automorphismes de $\CM$  permet de transposer sans ambigu{\"\i}t{\'e}
la correspondance de Langlands tordue \`a la Hecke au corps abstrait $C$. 
 Pour obtenir une correspondance de Langlands sur $C$,
il faut alors choisir 
(si besoin) une racine du cardinal du corps r{\'e}siduel de $K$.

\alin{Formulation continue $l$-adique}  \label{formcontl}Lorsque $C=\o\QM_l$,
nous avons besoin d'une formulation en termes de repr{\'e}sentations
continues $l$-adiques de $W_K$ plut{\^o}t que de repr{\'e}sentations de Weil-Deligne.
Voici  bri\`evement le lien entre les deux formulations expliqu\'e par Deligne
dans \cite[8]{DelAntwerp}.
Rappelons  que $W_K$ est muni de la topologie d{\'e}finie par la
topologie profinie de $I_K$ et la topologie discr{\`e}te de $W_K/I_K\simeq
\ZM$. On note alors  $\rep^d_l(W_K)$ l'ensemble des classes
d'{\'e}quivalences de $\o\QM_l$-repr{\'e}sentations Frobenius-semisimples 
%qui sont obtenues par 
%extension des scalaires de repr{\'e}sentations 
et {\em continues} de
$W_K$. 

Notons $\ZM_l(1):=\limproj\mu_{l^n}$ et $t_l : I_K\To{} \ZM_l(1)$ le $l$-quotient de l'inertie mod\'er\'ee, donn\'e par l'action sur les $l^n$-\`emes racines d'une uniformisante.
\`A tout prog\'en\'erateur $\mu$ de $\ZM_l(1)$, est donc associ\'e un morphisme surjectif
$ t_\mu:\;\; I_K\To{t_l} \ZM_l(1)\To{\mu^*} \ZM_l.$
Soit $(\sigma^{ss},N_\sigma,V_\sigma)$ une $\o\QM_l$-repr{\'e}sentation de Weil-Deligne
comme dans le paragraphe pr{\'e}c{\'e}dent et soit $\phi$
un rel{\`e}vement de Frobenius g{\'e}om{\'e}trique dans $W_K$. L'équation \ref{monod} montre 
 que la formule
\ini\begin{equation} \label{WD}
\sigma^\phi(w):= \sigma^{ss}(w) \hbox{exp}(N_\sigma t_\mu(i_\phi(w))), \;\;\hbox{ o{\`u}
  } w=\phi^{\nu(w)} i_\phi(w) \in \phi^\ZM\ltimes I_K
  \end{equation}
d{\'e}finit une repr{\'e}sentation continue et Frobenius-semisimple de $W_K$
sur l'espace  $V_{\sigma}$.
 Le th{\'e}or{\`e}me ``de la
monodromie $l$-adique''  de Grothendieck montre que l'application 
$\rep^d_{\o\QM_l}(WD_K) \To{} \rep^d_l(W_K)$ ainsi obtenue est une bijection.
Deligne a montr\'e qu'elle ne d\'epend pas des choix de $\phi$ et $\mu$.
L'op{\'e}rateur nilpotent $N_\sigma \in \endo{\o\QM_l}{V_{\sigma}}$ est
appel{\'e} ``monodromie'' et contr{\^o}le le d{\'e}faut de semi-simplicit{\'e} de la
repr{\'e}sentation $\sigma$. En particulier la semisimplifi{\'e}e de
$\sigma$ n'est autre que $\sigma^{ss}$.

%bijection avec l'ensemble $\rep^d_{\o\QM_l}(WD_K)$ d{\'e}fini
%ci-dessus. Plus pr{\'e}cis{\'e}ment, si $(\sigma,V)\in\rep^n_k(W_K)$, il
%affirme l'existence d'un endomorphisme nilpotent $N$ de $V$ tel que en
%posant pour tout $w=\phi^k.i$, $i\in I_K$
%$$ \sigma'(w):=\sigma(w) \exp(-Nt_l(i))$$ o{\`u} $t_l : I \simto \ZM_l$

La correspondance de Langlands locale, tordue \`a la Hecke, fournit donc une correspondance entre repr\'esentations irr\'eductibles de $G_d$ et repr\'esentations $l$-adiques et  c'est cette derni{\`e}re que la
g{\'e}om{\'e}trie ({\em i.e.} la cohomologie $l$-adique) peut pr{\'e}tendre r{\'e}aliser.
Moyennant le
choix d'une racine du cardinal du corps r{\'e}siduel dans $\o\QM_l$, on obtient la ``vraie"
correspondance, que nous noterons toujours
$$\pi \in \Irr{\o\QM_l}{G_d} \mapsto
\sigma_d(\pi) \in \rep^d_l(W_K).$$

\alin{Repr{\'e}sentations elliptiques} \label{corlanell}
Pour expliquer ce qu'il advient des repr{\'e}sentations elliptiques de
$G_d$ {\`a} travers la correspondance de Langlands, rappelons  que
$$ \sigma_d(\Disc{\o\QM_l}{G_d}) = \{\sigma \in \rep^d_l(W_K) \; \hbox{
  ind{\'e}composables}\}. $$
%Ainsi,  la composition  $ \sigma_d\circ JL_d :\; \Irr{\o\QM_l}{\dd}\To{} {\rep^d_l(W_K)}$
%induit une bijection de $\Irr{\o\QM_l}{\dd}$ sur 
%l'ensemble des ind{\'e}composables de $\rep^d_l(W_K)$, qui constitue un exemple  simple de fonctorialit{\'e} de Langlands.
\def\I#1{\{1,\cdots,#1\}}
%Maintenant, c
Comme le support cuspidal d'une  repr{\'e}sentation elliptique est 
aussi celui d'une  s{\'e}rie discr{\`e}te,
la compatibilit{\'e} de  la correspondance de Langlands {\`a} l'induction
parabolique donne la caract{\'e}risation
$$ \sigma_d(\Ell{\o\QM_l}{G_d}) = \{\sigma \in \rep^d_l(W_K), \exists \sigma'  \hbox{
  ind{\'e}composable telle que }  \sigma^{ss}={\sigma'}^{ss}\}. $$
Ainsi pour toute repr\'esentation elliptique $\pi$ de type $\rho$, on a
 $$\sigma_d(\pi)^{ss}= \sigma_{d/d_\rho}(\pi_\rho)\otimes \tau^{ss}_{d_\rho} |-|^{\frac{1-d_\rho}{2}} $$
 o\`u $\tau^{ss}_{d_\rho} := \o\QM_l \oplus \o\QM_l(1) \oplus \cdots \oplus \o\QM_l(d_\rho-1)$.
Pour tout sous-ensemble $I \subset \{1,\cdots, d_\rho-1\}$, introduisons la repr\'esentation $l$-adique $\tau_{d_\rho}^{I,\phi}$ associ\'ee \`a la repr\'esentation de Weil-Deligne $(\tau_{d_\rho}^{ss},N_I)$ par la formule \ref{WD}, et o\`u $N_I$ est l'op\'erateur de monodromie $\tau^{ss}_{d_\rho} \To{} \tau^{ss}_{d_\rho}(-1)$
 donn\'e par la matrice $\sum_{i\in I^c} E_{i,i-1}$, o\`u $I^c$ est le compl\'ementaire de $I$ (voir  \cite[4.1.1]{Dat1} pour plus de d\'etails). Nous noterons simplement $\tau_{d_\rho}^I$ la classe d'isomorphisme de cette repr\'esentation.

\begin{lemme} \label{corell}
  Soit $\rho \in \Irr{\o\QM_l}{D}$ et 
$I\subseteq S_{\rho}$. Utilisant les notations \ref{notationsrho} et la convention \ref{conv}, on a 
$$ \sigma_d(\pi^I_\rho) \simeq \sigma_{d/d_\rho}(\pi_\rho)\otimes
\sigma_{d_\rho}(\pi_1^{I}) \simeq 
\sigma_{d/d_\rho}(\pi_\rho)\otimes \tau_{d_\rho}^{I} \otimes |.|^{\frac{1-d_\rho}{2}}.$$
%o\`u l'on identifie $S_\rho$ et $S_{d_\rho}$ de mani\`ere \`a respecter l'ordre des blocs diagonaux.
\end{lemme}
\begin{proof}
Dans la preuve de \cite[4.1.2]{Dat1}, nous exhibons des entiers $d_1, \cdots, d_r$ de somme $d_\rho$, et  $n_1\geq \cdots \geq n_r $ tels que la repr\'esentation ${\pi}_1^I$ soit l'unique quotient de l'induite
$$ |det|_{K}^{n_1}\hbox{St}_{G_{d_1}}\times \cdots \times |det|_{K}^{n_r}\hbox{St}_{G_{d_r}}. $$
Ces entiers sont  ind\'ependants du corps $K$, de sorte que via l'\'equivalence de cat\'egories \ref{equivcat} et la convention \ref{conv}, on en d\'eduit que $\pi_\rho^I$ est l'unique quotient de l'induite 
$$ |det|^{n_1}\hbox{St}_{d_1}(\pi_\rho)\times\cdots\times |det|^{n_r}\hbox{St}_{d_r}(\pi_\rho), $$
ce qui nous fournit les param\`etres de Langlands de $\pi_\rho^I$. Alors d'apr\`es
\cite[2.9]{HeSMF}, on a
\begin{eqnarray*}
\sigma_d(\pi_\rho^I) & =  & |-|^{n_1}\sigma_{dd_1/d_\rho}(\hbox{St}_{d_1}(\pi_\rho))\oplus\cdots\oplus |-|^{n_r}\sigma_{dd_r/d_\rho}(\hbox{St}_{d_r}(\pi_\rho)) \\
 & = & |-|^{n_1} \sigma_{d/d_\rho}(\pi_\rho)\otimes \sigma_{d_1}(\hbox{St}_{G_{d_1}}) \oplus\cdots\oplus |-|^{n_r} \sigma_{d/d_\rho}(\pi_\rho)\otimes \sigma_{d_r}(\hbox{St}_{G_{d_r}}) \\
 & = & \sigma_{d/d_\rho}(\pi_\rho)\otimes \sigma_{d_\rho}(\pi_1^{I})
 \end{eqnarray*}
d'o\`u la premi\`ere \'egalit\'e de l'\'enonc\'e. La seconde d\'ecoule de \cite[4.1.2]{Dat1}.
\end{proof}

\section{Espaces modulaires de Drinfeld} \label{drinfeld}

Dans cette section, nous pr\'ecisons %rappelons  bri{\`e}vement la d{\'e}finition de 
la version des espaces de
modules de groupes formels que nous utilisons pour les \'enonc\'es cohomologiques principaux.
%consid{\'e}rons et la d{\'e}finition usuelle
%de leur cohomologie. 
Pour les sp{\'e}cialistes, il suffira de dire qu'on ne fixe pas la hauteur des quasi-isog{\'e}nies qui
rigidifient les probl{\`e}mes de modules (suivant en cela Rapoport-Zink). La cohomologie des espaces obtenus n'est donc pas admissible mais seulement admissible modulo le centre.
% et on
%{\em ne fixe pas} de caract{\`e}re central, du moins pour l'{\'e}nonc{\'e} du
%th{\'e}or{\`e}me principal.
On d{\'e}finit  ensuite les complexes de cohomologie $R\Gamma_c$ qui jouent
{\'e}videmment un r{\^o}le central dans ce texte.
On fixe dor\'enavant un entier $d\geq 1$, et on note $G:=GL_d(K)$ et $D$ l'alg\`ebre \`a division centrale sur $K$ d'invariant $1/d$, pr\'ec\'edemment not\'ee $D_d$.

\def\dd{{D^\times}}

\def\od{\OC_{D}}

%On notera $\OC$ l'anneau des entiers de $K$, $\od$ celui de $D$ 
%et $\wh\OC^{nr}$ celui de la compl{\'e}tion de la
%sous-extension maximale non ramifi{\'e}e $K^{nr}$ de $K$ dans $K^{ca}$.
%Le corps r{\'e}siduel de $K$ est $k$ et celui de $K^{nr}$ est une cl{\^o}ture
%alg{\'e}brique $k^{ca}$ de $k$. Enfin $\varpi$ d{\'e}signe une uniformisante de
%$\OC$ et $\wh\OC^{nr}$.
%Les espaces qui nous int{\'e}ressent sont des espaces de modules de
%$\OC$-modules formels. 

\subsection{La tour de Drinfeld} \label{tourDr}

La d\'efinition du complexe de cohomologie repose sur la description formelle suivante :

\ali \label{resumeDr}{\em  La tour de Drinfeld est un syst\`eme 
$\;\; \cdots \MC_{Dr,n} \To{\pi_{n,n-1}} \cdots \To{\pi_{1,0}} \MC_{Dr,0} \To{\xi_{Dr}} \Omega_K^{d-1} $, o\`u 

\begin{enumerate}
        \item Les $\MC_{Dr,n}$, $n\in \NM$ sont des $\knr$-espaces analytiques munis 
        \begin{enumerate}
        \item d'une action   de $G$ {\em continue} au sens de \cite[par. 6]{Bic3}.
        \item d'une action de $D^\times$ dont la restriction \`a $\OC_{D}^\times$  fait de $\pi_{n,0}:=\pi_{1,0}\circ\cdots \circ \pi_{n,n-1}$ un rev\^etement \'etale galoisien de groupe $\OC_{D}^\times/(1+ \varpi^n\OC_{D})$.
        \item d'une donn\'ee de descente \`a la Weil pour l'extension $\knr|K$ (au sens de \cite[3.45]{RZ}).
        \end{enumerate}
Ces actions commutent entre elles et aux  $\pi_{n,n'}$, et le sous-groupe $K^\times_{diag}$ de $G\times \dd$ agit trivialement.
\item l'augmentation $\xi_{Dr}$ est un morphisme d'espaces analytiques ``au-dessus de $K$" ({\em cf} \cite[p. 30]{Bic2}) qui, apr\`es extension des scalaires, devient localement isomorphique, compatible aux donn\'ees de descentes, et \'equivariant si l'on munit  l'espace sym\'etrique de Drinfeld $\Omega_K^{d-1}$ de l'action naturelle de $G$ et triviale de $\dd$.  
\end{enumerate}
}

\ali  Pour le confort du lecteur nous allons bri\`evement rappeler la d\'efinition de ces objets. On reprend les notations de l'introduction, en abr\`egeant $\OC:=\OC_K$ et en fixant une uniformisante $\varpi$ de $\OC$.

Si $B$ est une $\OC$-alg{\`e}bre, un $\OC$-module
formel sur $B$ est un groupe formel muni d'une action de $\OC$
relevant l'action naturelle sur l'alg{\`e}bre de Lie. Un $\od$-module
formel sur $B$ est un $\OC$-module formel muni d'une action de $\od$
qui {\'e}tend celle de $\OC$. Notons $\OC_d \subset \od$ l'anneau des
entiers d'une sous-extension non ramifi{\'e}e maximale de $K$ dans $D$.
 Suivant Drinfeld, un $\od$-module formel
sur $B$ est dit sp{\'e}cial si son alg{\`e}bre de Lie est localement libre de rang
$1$ comme $\OC_d\otimes_\OC B$-module.

Une isog{\'e}nie de $\OC$-modules formels est une isog{\'e}nie des groupes
formels sous-jacents compatible aux actions de $\OC$. Sa
($\OC$-)hauteur est le quotient (entier) de la hauteur au sens des
groupes par le degr{\'e} de $k$ sur $\FM_p$. La hauteur d'un
$\OC$-module formel $X$ est la hauteur de l'isog{\'e}nie $X(\varpi)$.
Une quasi-isog{\'e}nie $X\To{} Y$ de $\OC$-modules formels est un
{\'e}l{\'e}ment de $\hom{X}{Y}{\OC-mf}\otimes_\OC K$ qui admet un inverse
dans $\hom{Y}{X}{\OC-mf}\otimes_\OC K$. On montre que c'est en fait une isog{\'e}nie {\`a}
multiplication par une puissance de $X(\varpi)$ pr{\`e}s, ce qui permet d'en
d{\'e}finir la hauteur.

\alin{L'espace $\MC_{Dr,0}$} \label{modulesDr}
\def\nilp{\hbox{Nilp}}
La d\'efinition de la tour de Drinfeld
 repose sur l'existence d'un $\od$-module formel $\XM$ sp{\'e}cial (donc de
dimension $d$) de hauteur $d^2$ sur $k^{ca}$, qui est unique {\`a} isog{\'e}nie
pr{\`e}s, \cite[Prop II.5.2]{BouCar}, \cite[3.60]{RZ}.  

Soit $\nilp$ la cat{\'e}gorie des $\wh\OC^{nr}$-alg{\`e}bres o{\`u} l'image
de $\varpi$ est nilpotente.
% Une {\em d{\'e}formation par quasi-isog{\'e}nie}
%de $\XM$ sur une telle alg{\`e}bre $B$ est un 
 On consid{\`e}re le foncteur $\wt{G} : \nilp
\To{} \hbox {Ens}$ qui {\`a} $B$ associe l'ensemble des classes
d'isomorphisme de %d{\'e}formations par quasi-isog{\'e}nie de $\XM$ sur $B$. 
couples $(X,\rho)$ avec $X$
un $\od$-module formel sur $B$ et $\rho : \XM \otimes_{k^{ca}} (B/\varpi B)
\To{} X \otimes_B (B/\varpi B)$ une quasi-isog{\'e}nie.
On
a une d{\'e}composition {\'e}vidente $\wt{G} = \bigsqcup_{h\in \ZM} G^{(h)}$ o{\`u}
$G^{(h)}$ classifie les classes de couples $(X,\rho)$ avec $\rho$ de
hauteur $dh$. Chaque $G^{(h)}$ est non-canoniquement isomorphe {\`a}
$G^{(0)}$ et $G^{(0)}$ est le foncteur $G$ de Drinfeld
\cite{Drincov}. On sait que $G$ est repr{\'e}sentable ({\em cf}
\cite{Genestier}, \cite{RZ}) par un
sch{\'e}ma formel localement de type fini sur $\wh\OC^{nr}$. On note
$\wh{\MC}_{Dr,0}^{(0)}$ ce sch{\'e}ma formel. De m{\^e}me on note $\wh{\MC}_{Dr,0}^{(h)}$ le
sch{\'e}ma formel repr{\'e}sentant $G^{(h)}$ et $\wh{\MC}_{Dr,0}$ celui qui repr{\'e}sente
$\wt{G}$. On a donc non-canoniquement $\wh{\MC}_{Dr,0}\simeq
\wh{\MC}_{Dr,0}^{(0)}\times \ZM$. 
Enfin on note sans $\;\wh{}\;$ les fibres g{\'e}n{\'e}riques au sens de
Raynaud-Berkovich de ces espaces : ce sont donc des $\knr$-espaces
analytiques au sens de \cite{Bic2}.

%Fixons un entier $d\geq
%1$. Notons $\Omega$ le ``demi-plan sup{\'e}rieur'' de Drinfeld de
%dimension $d-1$ et $(\Sigma_n)_n$ la tour de rev{\^e}tements {\'e}tales d{\'e}finie
%dans \cite{Drincov}. En particulier $\Sigma_0=\Omega \wh\otimes_K \knr$.
%Ces espaces ont {\'e}t{\'e} d{\'e}finis comme des espaces analytiques rigides (sur
%$K$ pour le demi-plan et sur $\wh{K^{nr}}$ pour ses rev{\^e}tements) mais
%nous les verrons g{\'e}n{\'e}ralement comme des espaces analytiques au
%sens de Berkovich, pour leur appliquer sa th{\'e}orie de la cohomologie
%{\'e}tale {\'e}quivariante.

%Suivant \cite{HaCusp} et \cite{Hausb} on a besoin d'une variante de ces espaces : pour
%cela on fixe un g{\'e}n{\'e}rateur topologique $\phi$ de $\gal(\wh{K^{nr}}/K)$
%et on consid{\`e}re les $\knr$-espaces analytiques suivants (moralement
%des changements de base de ``restriction des scalaires {\`a} la Weil''
%pour l'extension $\knr/K$) :
%$$\breve{\Sigma}_n := \bigsqcup_{n\in \ZM}
%\Sigma_n\otimes_{\knr,\phi^n} \knr.$$
%Enfin nous abr{\`e}gerons $\breve\Sigma_n^{ca}:= \breve\Sigma_n
%\wh\otimes_\knr \ka$.

\alin{Structures de niveau}
Notons $({X}_u,\rho_u)$ l'objet universel au-dessus de $\wh\MC_{Dr,0}$. Le noyau
${X}_u[\varpi^n]$ de
la multiplication par $\varpi^n$ dans ${X}_u$ est un sch{\'e}ma formel en groupes
plat  fini de rang  $p^{nd^2}$ au-dessus de $\wh\MC_{Dr,0}$, et qui
est {\'e}tale en fibre g{\'e}n{\'e}rique. Plus 
pr{\'e}cis{\'e}ment, sa fibre g{\'e}n{\'e}rique est
localement pour la topologie {\'e}tale isomorphe {\`a} $\mdro \times
\od/\varpi^{n}\od$. Le
$(\od/\varpi^n\od)^\times$-torseur sur $\mdro$
$$ \underline{\hbox{Isom}_{\od}}\left((\varpi^{-n}\od/\od)_{\mdro},{X}_u[\varpi^n]\right)$$
est  donc repr{\'e}sent{\'e} par un rev{\^e}tement {\'e}tale de $\mdro$, galoisien  de
groupe $\od^\times /(1+\varpi^n\od)$, qui est le $\mdrn$ de \ref{resumeDr}. %\footnote{Cette d{\'e}finition des rev{\^e}tements
%  est emprunt{\'e}e {\`a} \cite{Genestier}}. 
Pour
$m\leq n$, l'inclusion $\varpi^{-m}\od/\od \subset \varpi^{-n}\od/\od$
induit le
%a multiplication par $\varpi^m$ sur $X[\varpi^n]$ induit un
morphisme $\mdrn\To{\pi_{n,m}} \MC_{Dr,m}$ de \ref{resumeDr}. %On obtient ainsi une ``tour'' de
%rev{\^e}tements Galoisiens de $\mdro$ dont le pro-groupe de Galois est
%$\od^\times$. 

\alin{Actions des groupes}  \label{actionDr}
On sait  que le groupe des
quasi-isog{\'e}nies du $\od$-module formel 
$\XM$ s'identifie {\`a} $G=GL_d(K)$, \cite[3.60]{RZ} ou
\cite[II.5.2]{BouCar}. On r{\'e}cup{\`e}re donc une action ({\`a} gauche) de $G$ 
sur le foncteur $\wt{G}$ qui envoie le couple $(X,\rho)$ sur le couple
$(X,\rho \circ ({^tg})_\XM)$ o{\`u} $g_\XM$ d{\'e}signe la
quasi-isog{\'e}nie de $\XM$ associ{\'e}e {\`a} $g$ {\em et ${^tg}$ d\'esigne la transpos\'ee de $g$}\footnote{Cette normalisation de l'action permet 
%Nous nous \'ecartons ici des conventions usuelles. Il y a trois raisons \`a cela : d'une part, 
de rendre l'\'enonc\'e final \ref{main} plus joli car exempt de contragr\'ediente.
% d'autre part, cela rendra le morphisme des p\'eriodes de Drinfeld compatible avec l'action sur $\Omega$ de la partie \ref{dp}, et enfin, 
 De plus le th\'eor\`eme de Faltings fait intervenir un passage \`a la transpos\'ee.}. Cette action en induit
une sur $\wh{\MC}_{Dr,0}$ par Yoneda, et l'objet universel est \'equivariant, {\em i.e.} muni
d'isomorphismes canoniques 
$({X}_u,\rho_u\circ g_\XM^{-1}) \simto g^*({X}_u,\rho_u) $. 
%En particulier ${X}_u$ est
%$G$-{\'e}quivariant au-dessus de $\wh\MC_{Dr,0}$.
Par cons{\'e}quent tous les $\mdrn$ sont
munis d'une action de $G$ %(explicitement d{\'e}finie par
%$$ X_u[\varpi^n] \simto g^*X_u[\varpi^n] = X_u[\varpi^n]\times_{\wh\MC_{Dr,0},
%  g} \wh\MC_{Dr,0} \To{proj_1} X_u[\varpi^n] )$$
et les morphismes de transition sont
$G$-{\'e}quivariants. D'apr{\`e}s Berkovich, l'action de $G$ sur ces $\knr$-espaces
analytiques est continue au sens de \cite[par. 6]{Bic3} pour la topologie
naturelle de $G$.

\medskip

On d{\'e}finit maintenant l'action de $\dd$ sur $\MC_{Dr,0}$. Pour $d\in
\dd$ et $X$ un $\od$-module formel sur $B$, notons ${^dX}$ le
$\od$-module formel dont le $\OC$-module formel sous-jacent est
encore $X$ mais dont l'action de $\od$ est donn{\'e}e par
${^dX}(x):=X(d^{-1}xd)$, de sorte que $X(d^{-1})$ est une quasi-isog{\'e}nie
$X\To{} {^dX}$. Pour un couple $(X,\rho)$ on pose
$d(X,\rho):=({^dX}, \rho \circ \XM(d^{-1}))$. On obtient ainsi
une action {\`a} gauche de $\dd$  sur le foncteur
$\wt{G}$, et donc 
sur $\wh\MC_{Dr,0}$, triviale sur $\OC_{D}^\times$. L'objet universel est alors muni d'isomorphismes
canoniques $d^*(X_u,\rho_u) \simto ({^dX_u},\rho_u\circ \XM(d^{-1}))$ et cela permet de prolonger l'action naturelle de $\OC_{D}^\times$ sur les $\MC_{Dr,n}$ \`a $\dd$, de mani\`ere compatible aux morphismes de transition.

%On d{\'e}finit  la structure $\dd$-{\'e}quivariante  sur
%$\mdrn$ au-dessus de $\mdro$ par la suite d'isomorphismes
%\begin{eqnarray*}
%    d^*  \underline{\hbox{Isom}}_{\od}\left(\varpi^{-n}\od/\od,
%X_u[\varpi^n]\right) & \simeq &   \underline{\hbox{Isom}}_{\od}\left(\varpi^{-n}\od/\od,
%d^*X_u[\varpi^n]\right) \\
%& \simeq &  \underline{\hbox{Isom}}_{\od}\left(\varpi^{-n}\od/\od,
%X_u^d[\varpi^n]\right) \\
%& \simto &  \underline{\hbox{Isom}}_{\od}\left(\varpi^{-n}\od/\od,
%X_u[\varpi^n]\right)
%\end{eqnarray*}
%o{\`u} les deux premiers isomorphismes sont canoniques et le troisi{\`e}me est
%donn{\'e} par $\phi \mapsto \phi\circ\,\hbox{ad}(d^{-1})$, en notant
%$\hbox{ad}(d^{-1})$ l'application $a\in \varpi^{-n}\od/\od \mapsto
%d^{-1}ad \in \varpi^{-n}\od/\od$.   

%Comme pr{\'e}c{\'e}demment pour l'action de $G$, ces structures
%$\dd$-{\'e}quivariantes induisent (et 
%{\'e}quivalent {\`a}) des actions de $\dd$ sur les $\mdrn$ compatibles avec les
%morphismes de transition.
%Notons que si
%$d\in \od^\times$, l'action de $d$ est triviale sur $\mdro$ et co{\"\i}ncide
%sur $\mdrn$ avec l'action naturelle sur le torseur qui d{\'e}finit
%$\mdrn$. En particulier l'action de $\dd$ se factorise sur chaque $\mdrn$
%par un quotient discret.

\medskip

Passons {\`a} l'action de Galois et fixons un g{\'e}n{\'e}rateur topologique $\phi$
de $\gal(K^{nr}/K)=\gal(k^{ca}/k)$. Rappelons qu'une donn\'ee de descente \`a la Weil \cite[3.45]{RZ} est une
 %on peut munir les $\knr$-espaces
%analytiques $\mdrn$ d{\'e}finis ci-dessus de 
structure $\phi^\ZM$-{\'e}quivariante au-dessus de $\knr$.
%\footnote{Remarquons que la donn{\'e}e d'une structure $\phi^\ZM$-{\'e}quivariante sur un $\knr$-espace  analytique  {\'e}quivaut {\`a} celle d'une action de $\phi^\ZM$ sur l'``espace analytique au-dessus de $K$'' sous-jacent ({\em cf}    \cite{Bic2} pour la terminologie), compatible {\`a} l'action de    $\phi^\ZM$ sur $\MC(\knr)$} 
%(appel{\'e}es ``donn{\'e}es de descente {\`a} la Weil'' dans
%\cite[3.45]{RZ}). 
%Commen{\c c}ons par $\wh\MC_{Dr,0}$ : 
Ceci \'etant,
si  $(X,\rho)$ est un
$\od$-module formel rigidifi{\'e}
au-dessus de $B$, alors on en d{\'e}finit un autre $(\phi^*X,{^\phi\rho})$
au-dessus de $\phi^*B :=B\otimes_{\OC^{nr},\phi} \OC^{nr}$ en posant
$\phi^*X:=X\otimes_{B} \phi^*B$ et 
$$^\phi\rho:\; \XM \otimes_{k^{ca}} \phi^*(B/\varpi) \To{Frob\otimes \id}
\phi^*\XM \otimes_{k^{ca}} \phi^*(B/\varpi) \To{\phi^*\rho} \phi^*X
\otimes_{\phi^* B} \phi^* (B/\varpi).$$
On obtient ainsi un morphisme de foncteurs $\wt{G} \To{} \phi_*\wt{G}$
dont on voit imm{\'e}diatement qu'il est inversible, ce qui nous donne une
structure $\phi^\ZM$-{\'e}quivariante sur $\wt{G}$, et donc sur $\wh\MC_{Dr,0}$, au-dessus 
de $\hbox{Spf}(\wh\OC^{nr})$. Nommons momentan{\'e}ment $\tau_\phi:
\phi^*\wh\MC_{Dr,0}\simto \wh\MC_{Dr,0}$ le $\hbox{Spf}(\wh\OC^{nr})$-isomorphisme 
associ{\'e}  {\`a} cette $\phi^\ZM$-structure {\'e}quivariante (adjoint du pr{\'e}c{\'e}dent), on obtient pour
l'objet universel l'existence d'un isomorphisme canonique
$\sigma_\phi:\;\;(\phi^*(X_u),{^\phi\rho_u})\simto
\tau_\phi^*(X_u,\rho_u) $  au-dessus de $\phi^*\wh\MC_{Dr,0}$, d'o{\`u}, par
composition, une structure $\phi^\ZM$-{\'e}quivariante
$$ \phi^*(X_u) \To{\sigma_\phi} \tau_\phi^*X_u =
X_u\times_{\wh\MC_{Dr,0}} \phi^*\wh\MC_{Dr,0} \To{proj_1} X_u $$
au-dessus de $\hbox{Spf}(\wh\OC^{nr})$.
En passant aux fibres g{\'e}n{\'e}riques, on obtient sur les $\mdrn$ 
les donn\'ees de descente annonc\'ees dans \ref{resumeDr} i) (c).
%les structures
%$\phi^\ZM$-{\'e}quivariantes au-dessus de $\MC(\knr)$ (notation de
%\cite{Bic2}) cherch{\'e}es. Par construction, elles sont compatibles aux
%morphismes de transition. Remarquons que l'action de $\phi$ sur
%$\mdrn$ envoie $\mdrn^{(h)}$ sur $\mdrn^{(h+1)}$.

\alin{Sur $\ka$}
{\'E}tendons maintenant les scalaires {\`a} $\ka$, en suivant la d\'efinition de \cite[1.4]{Bic2}.
 On obtient une tour de $\ka$-espaces analytiques 
 $$\mdrn^{{ca}}:=\mdrn
\wh\otimes_{\knr} \ka, \;\; n\in \NM $$
que nous avons not\'ee simplement $\MC^{{ca}}_{Dr}$ dans l'introduction. 
 La structure
$I_K$-{\'e}quivariante sur $\mdrn^{{ca}}$ au-dessus de $\ka$ induite par ce changement
de base et la structure $\phi^\ZM$-{\'e}quivariante pr{\'e}c{\'e}dente se
``recollent''  en une structure $W_K$-{\'e}quivariante au-dessus de
$\MC(\ka)$. De mani{\`e}re {\'e}quivalente, on obtient une action de
$W_K$ {\em \`a droite} sur l'``espace analytique   
au-dessus de $K$'' $\mdrn^{{ca}}$, compatible {\`a} l'action naturelle de
$W_K$ sur $\ka$.  
Ainsi la tour $\MC^{ca}_{Dr}$ est munie d'une action de $G\times \dd \times W_K$, et on v\'erifie sur les d\'efinitions que celle ci est triviale sur $K^\times_{diag} \times 1$. 

%\medskip 

%Les actions de $G, \dd$ et $W_K$ qu'on vient de d{\'e}finir sur les espaces
%analytiques $\mdrn^{{ca}}$ commutent. On v{\'e}rifie facilement que le
%stabilisateur de $\MC^{d/K,(h)}_{Dr,n}$ 
%est $(G\times \dd\times W_K)^0$, avec la notation introduite en \ref{defun}.
%Remarquons pour terminer que l'action de $G\times \dd\times W_K$ est triviale sur
%le sous-groupe ${K^\times_{diag}}$ de $G\times \dd$ 
%d\'efini par 
%\ini
%\begin{equation} \label{defDelta}
%{K^\times_{diag}}=\hbox{Im}\left(\application{}{K^\times}{G\times\dd}{z}{(z,z)}\right)
%\end{equation}

\alin{Le morphisme de p{\'e}riodes} \label{perDr}
Soit $\Omega^{d-1}_K$ 
l'espace sym{\'e}trique
de Drinfeld, d\'efini comme le compl\'e\-mentaire des hyperplans rationnels dans $\PM^{d-1}_K$, et $\Omega^{d-1,nr}_K$ % := \Omega^{d-1}_K\wh\otimes_K \knr$ 
son changement de base {\`a} $\knr$. Drinfeld a construit dans
\cite{Drincov} un isomorphisme $\xi^{(0)}:\;\;\wh\MC_{Dr,0}^{(0)}\To{}
\wh{\Omega}^{d-1,nr}_K$ o{\`u} $\wh{\Omega}^{d-1,nr}_K$ d{\'e}signe un certain
mod{\`e}le formel de ${\Omega}^{d-1,nr}_K$ d{\'e}fini auparavant par
Deligne.
Nous n'avons besoin ici que de la fibre g{\'e}n{\'e}rique de ce
morphisme $\xi$, mais nous avons par contre besoin  de l'{\'e}tendre {\`a} tout
$\MC_{Dr,0}$. Pour cela nous utilisons l'isomorphisme de $\knr$-espaces
analytiques $\tau_{\phi^h}:\;\;(\phi^h)^*\mdro \To{} \mdro$ qui induit
un isomorphime $(\phi^h)^*\MC_{Dr,0}^{(0)} \simto \MC_{Dr,0}^{(h)}$
et nous d{\'e}finissons $\xi^{(h)}$ par le diagramme commutatif
$$ \xymatrix{ (\phi^h)^* \mdro^{(0)} \ar[d]_{(\phi^h)^*\xi^{(0)}} \ar[r]^\sim & \mdro^{(h)}
  \ar@{-->}[d]^{\xi^{(h)}} \\ (\phi^h)^*\Omega^{d-1,nr}_K \ar[r]^\sim
  & \Omega^{d-1,nr}_K }$$
o{\`u} la fl{\`e}che du bas est la donn{\'e} de descente naturelle sur
$\Omega^{d-1,nr}_K$. Enfin on pose
$\xi := \bigsqcup_{h\in\ZM} \xi^{(h)} : \MC_{Dr,0} \To{} \Omega^{d-1,nr}_K$.
Une construction directe (et plus ``simple'' que celle de Drinfeld) de
ce $\xi$ est propos{\'e}e dans le 
chapitre 5 de Rapoport-Zink \cite{RZ}, mais n'est r{\'e}dig{\'e}e qu'en
in{\'e}gales caract{\'e}ristiques (techniques cristallines).

%\begin{fact} Le morphisme
%$$ \xi:\;\; \mdro \To{} \Omega^{d-1,nr}_K $$  est
%$G\times \dd\times W_K$ {\'e}quivariant pour l'action naturelle de $G$ sur
%$\Omega^{d-1,nr}_K$ (qui se factorise donc par $PGL_d(K)$), l'action
%triviale de $\dd$ et l'action naturelle de $W_K$ (celle donn{\'e}e par
%l'extension des scalaires).
%\end{fact}

%\begin{proof}
Les assertions de \ref{resumeDr} ii) sont contenues dans les r\'ef\'erences usuelles \cite{Drincov}, \cite{BouCar}, \cite{Genestier} et \cite{RZ}. Commentons seulement la compatibilit\'e \`a l'action de $G$ : dans les r\'ef\'erences pr\'ec\'edentes  l'action de $G$ sur les $\MC_{Dr,n}$ est l'action ``naturelle" dont la notre se d\'eduit par $g\mapsto {^tg}^{-1}$. Le morphisme de p\'eriode y est \'equivariant pour l'action naturelle sur le mod\`ele ``dual" de $\Omega^{d-1}_K$, \`a savoir l'ensemble des hyperplans de $K^d$ ne contenant pas de droites rationnelles. C'est pourquoi, avec notre normalisation diff\'erente et notre d\'efinition de $\Omega_K^{d-1}$, le morphisme de p\'eriodes reste $G$-\'equivariant.

%\end{proof}

\subsection{La tour de Lubin-Tate} \label{tourLT}

\ali \label{resumeLT} {\em C'est un syst\`eme
$\;\; \cdots \MC_{LT,n} \To{\pi_{n,n-1}} \cdots \To{\pi_{1,0}} \MC_{LT,0} \To{\xi_{LT}} S_K^{1/d} $, o\`u 
\begin{enumerate}
        \item Les $\MC_{LT,n}$, $n\in \NM$ sont des $\knr$-espaces analytiques munis 
        \begin{enumerate}
        \item d'une action {\em continue} de $\dd$,
        \item d'un syst\`eme de morphismes \'etales finis  $g^{n'|n} : \MC_{LT,n} \To{} \MC_{LT,n'}$
d\'efinis pour tous $g,n,n'$ satisfaisant  $g M_d(\OC_K)g^{-1} \subset \varpi^{n'-n}M_d(\OC_K)$ et tels que 
\begin{enumerate}
\item $\forall g,h \in G$, $g^{n''|n'}h^{n'|n}= (gh)^{n''|n}$ lorsque tous sont bien d\'efinis.
\item $\forall n,n'$, $\pi_{n,n'}=1^{n'|n}$.
\item $\forall n$, les $g^{n|n}$ font de $\pi_{n,0}$  un rev\^etement \'etale galoisien de groupe $GL_d(\OC_K)/(1+\varpi^nM_d(\OC_K))$.
\end{enumerate}         
        \item d'une donn\'ee de descente \`a la Weil pour l'extension $\knr|K$.
        \end{enumerate}
Ces donn\'ees sont compatibles entre elles et l'action de $K^\times_{diag}\subset \dd\times G$ est triviale sur chaque $\MC_{LT,n}$.
%\item les $\pi_{n,n'}$ sont des morphismes de $\knr$-espaces analytiques compatibles \`a ces structures (et \'etales finis par le point pr\'ec\'edent).
\item l'augmentation $\xi_{LT}$ est un morphisme d'espaces analytiques au-dessus de $K$, qui apr\`es extension des scalaires, devient \'etale surjectif, compatible aux donn\'ees de descente, et \'equivariant si l'on munit la vari\'et\'e de Severi-Brauer $S_K^{1/d}$ d'invariant $1/d$ de l'action naturelle de $\dd$ et  triviale de $GL_d(K)$. On a donc $\xi_{LT}\circ g^{0|n} = \xi_{LT} \circ 1^{0|n}$ pour tous $g,n$ tels que $gM_d(\OC_K)g^{-1} \subset \varpi^{-n}M_d(\OC_K)$.  
\end{enumerate}
}

Nous donnons maintenant un peu de substance \`a cette enveloppe formelle en rappelant la d\'efinition de ces objets. Elle repose sur l'existence d'un $\OC_K$-groupe formel $\XM$ sur $k^{ca}$ de dimension $1$ et
hauteur $d$, au sens de \cite{Drinell} ou \cite{HG} et qui est unique {\`a} isomorphisme pr{\`e}s, \cite[prop
1.6-7]{Drinell}.  
%De plus sa
%$\OC_K$-alg{\`e}bre d'endomorphismes s'identifie {\`a} $\od$, l'anneau des
%entiers de la $K$-alg{\`e}bre {\`a} division centrale d'invariant $1/d$.
\alin{D{\'e}formations de $\XM$} \label{modulesLT}
On peut comme dans le livre
de Rapoport-Zink exprimer le probl{\`e}me de modules sur la cat{\'e}gorie
$\nilp$ et de mani{\`e}re semblable au cas de la tour de Drinfeld. Nous
donnons malgr{\'e} tout la d{\'e}finition historique  en termes de d{\'e}formations, car c'est la seule {\'e}crite en {\'e}gales caract{\'e}ristiques.
 Soit $\CC$ la
cat{\'e}gorie des $\wh\OC^{nr}$-alg{\`e}bres locales de corps r{\'e}siduel
$\o\FM_p$ et compl{\`e}tes pour leur 
topologie adique. Une {\em d{\'e}formation par quasi-isog{\'e}nie} de $\XM$ sur une
telle alg{\`e}bre $R$ est une paire $(X,\rho)$ form{\'e}e d'un $\OC_K$-module
formel $X$ sur $R$ et d'une quasi-isog{\'e}nie $\XM\To{\rho}X\otimes_R
R/\MG_R$.  Notons $\hbox{Def}$ le foncteur de   
$\CC$ dans les ensembles qui 
{\`a} $R$ associe l'ensemble des classes d'isomorphisme de d{\'e}formations
par quasi-isog{\'e}nie $(X,\rho)$ de $\XM$ sur $R$. Ce foncteur est une
r{\'e}union disjointe de sous-foncteurs $\hbox{Def}^{(h)}$ classifiant les
couples $(X,\rho)$ avec $\rho$ de hauteur $h$. Tous les $\hbox{Def}^{(h)}$
sont non-canoniquement isomorphes {\`a} $\hbox{Def}^{(0)}$. De plus, comme une
quasi-isog{\'e}nie de hauteur nulle entre deux $\OC_K$-modules formels de
dimension $1$ sur $\o\FM_p$ est un isomorphisme (c'est une
cons{\'e}quence de \cite{Drinell},props 1.6-2 et 1.7), on voit que
$\hbox{Def}^{(0)}$ est le foncteur de d{\'e}formations {\'e}tudi{\'e}
par Drinfeld.  
D'apr{\`e}s \cite[4.2]{Drinell} on sait que $\hbox{Def}^{(0)}$ est
repr{\'e}sentable 
par l'alg{\`e}bre $R^{(0)}:=\wh\OC^{nr}[[T_1,\cdots,T_{d-1}]]$  des s{\'e}ries
formelles {\`a} $d-1$ variables sur $\wh\OC^{nr}$. Nous noterons alors
$\MC_{LT,0}^{(0)}$ le $\knr$-espace analytique de Berkovich
associ{\'e} au sch{\'e}ma 
formel $\wh\MC_{LT,0}^{(0)}:=\hbox{Spf}(R^{(0)})$ (suivant la proc{\'e}dure de
Raynaud-Berthelot d{\'e}crite dans \cite{Bic4}) : c'est la boule
unit{\'e} ouverte de dimension $d-1$. 
Il s'ensuit que les foncteurs $\hbox{Def}^{(h)}$ pour $h\in \NM$ sont
aussi repr{\'e}sentables et nous noterons $\MC_{LT,0}^{(h)}$ les espaces
analytiques associ{\'e}s. Nous posons enfin
$\wh\MC_{LT,0}:=\bigsqcup_{h\in\NM} \wh\MC_{LT,0}^{(h)}$ et notons sans
$\;\wh{}\; $ l'espace analytique associ{\'e}.

\alin{Structures de niveau} On peut suivre  la m{\^e}me proc{\'e}dure que dans
le cas $Dr$. En notant encore $(X_u,\rho_u)$ l'objet universel sur
$\wh\MC_{LT,0}$, on d{\'e}finit $\mltn$ comme le rev{\^e}tement {\'e}tale
galoisien de
groupe $GL_d(\OC/\varpi^n\OC)$ de $\mlto$ repr{\'e}sentant le torseur
$$ \underline{\hbox{Isom}}((\varpi^{-n}\OC/\OC)^d,X_u[\varpi^n]).$$
Il se trouve que dans ce cas $LT$ on peut faire mieux en interpr{\'e}tant 
$\MC_{LT,n}$ comme la fibre g{\'e}n{\'e}rique d'un probl{\`e}me de modules
classifiant les ``structures de niveau de Drinfeld'' : cela  a
l'avantage de simplifier la description de l'action de
$G$\footnote{On peut n{\'e}anmoins se passer d'une telle
  interpr{\'e}tation modulaire comme il est esquiss{\'e}
  dans \cite[5.34]{RZ} et expliqu{\'e} dans \cite[2.3.8.3]{Fargues}}.
Soit  $\Lambda \subset K^d$ un $\OC$-r{\'e}seau.
Si $R\in \CC$ et $(X,\rho)$ est une d{\'e}formation par quasi-isog{\'e}nie
de $\XM$ au-dessus
de $R$, alors une $\Lambda$-structure de niveau $n$ de Drinfeld sur $X$ est un
morphisme de $\OC$-modules $\psi$ de  $\varpi^{-n}\Lambda/\Lambda$ vers l'id{\'e}al
maximal $\MG_R$ muni de la structure de $\OC$-module d{\'e}finie par $X$
et tel que
$$ \prod_{x\in \varpi^{-1}\Lambda/\Lambda} (T-\psi(x)) \;\;\hbox{ divise }\;\;
  X_\varpi(T) $$
o{\`u} $X_\varpi(T)$ d{\'e}signe la s{\'e}rie formelle donnant l'action de
$\varpi$ sur $X$.
Les foncteurs $\hbox{Def}^{(h)}_{\Lambda,n}$ classifiant les triplets $(X,\rho,\psi)$ {\`a}
isomorphisme pr{\`e}s sont repr{\'e}sentables \cite[4.3]{Drinell} par des
alg{\`e}bres $R^{(h)}_{\Lambda,n}$ finies et 
plates sur les $R^{(h)}$ et nous noterons $\wh\MC_{LT,\Lambda,n}:=\bigsqcup_{h\in
  H} \hbox{Spf} R^{(h)}_{\Lambda,n}$. Lorsque $\Lambda$ est le
r{\'e}seau ``canonique'' $\OC^d\subset K^d$, la fibre g{\'e}n{\'e}rique de
$\wh\MC_{LT,\Lambda,n}$ s'identifie canoniquement {\`a} $\MC_{LT,n}$. 
{\'E}videmment pour $n$ fix{\'e}, tous les $\wh\MC_{LT,\Lambda,n}$ sont
(non-canoniquement) isomorphes.

%Alternativement, on peut aussi incorporer des structures de
%niveau $n\in \NM$ dans le probl{\`e}me de d{\'e}formations $\hbox{Def}^{(0)}$
%pr{\'e}c{\'e}dent en utilisant les 
%``bases de Drinfeld''. On obtient alors un foncteur repr{\'e}sentable par un
%anneau $R_{K,d}^n$ plat et fini sur $R_{K,d}$ et dont la fibre
%g{\'e}n{\'e}rique du sch{\'e}ma formel associ{\'e} n'est autre que $\MC_{LT,n}^{(0)}$.

%Si maintenant $K$ est un sous-groupe ouvert compact de $GL_d(\OC)$, on
%peut trouver $n>0$ tel que $K\supset 1+\varpi^nM_d(\OC)$ et on pose alors
%$\MC_{LT,K}:=\MC_{LT,n}/(K/(1+\varpi^nM_d(\OC))) $ o{\`u} l'on voit
%$K/(1+\varpi^nM_d(\OC))$ comme un sous-groupe de $GL_d(\OC/\varpi^n\OC)$. On
%v{\'e}rifie alors que l'espace analytique $\MC_{LT,K}$ ne d{\'e}pend pas, {\`a}
%isomorphisme canonique pr{\`e}s, du choix de $n$.

\alin{Action des groupes} \label{actionLT}
On sait que le groupe des quasi-isog{\'e}nies du
$\OC$-module formel $\XM$ s'identifie {\`a} $\dd$. Ceci nous permet de
d{\'e}finir une action \`a gauche de $\dd$ sur $\mlto$, puis sur les
$\mltn$
d'une mani{\`e}re exactement analogue {\`a}  celle par laquelle on a d{\'e}fini
l'action de $G$ sur $\mdro$ et les $\mdrn$. Explicitement, on envoie
un triplet $(X,\rho,\psi)$ sur le triplet $(X, \rho\circ d_{\XM}^{-1},\psi)$.

\medskip

L'action de $G$ est plus d{\'e}licate. Elle est soigneusement d{\'e}finie
dans \cite[II.2]{HaTay} 
(o{\`u} les d{\'e}formations sont cependant par isomorphismes) et
\cite{Strauch} (qui traite les d{\'e}formations par
quasi-isog{\'e}nies). On rappelle ici bri{\`e}vement cette
d{\'e}finition, en l'exposant un peu diff\'eremment des r\'ef\'erences
usuelles. Pour deux r\'eseaux $\Lambda, \Lambda'$ dans $K^d$, nous
noterons $d(\Lambda,\Lambda')$ la distance combinatoire entre les
points de l'immeuble de $PGL_d(K)$ associ\'es. C'est le minimum de la
somme $r+k$ pour tous les couples d'entiers
% Supposons 
%donn{\'e}s deux r{\'e}seaux $\Lambda, \Lambda'$  dans $K^d$ et choisissons deux
%entiers
 $r,k\in\ZM$  tels
tels que 
\ini\begin{equation}\label{condreseaux}
\Lambda\subseteq
\varpi^{-r} \Lambda'\subseteq \varpi^{-r-k}\Lambda.
\end{equation} 
\begin{lemme} \label{blob} Il existe une famille de morphismes
  $\dd$-\'equivariants de foncteurs
$$%\ini\begin{equation} \label{blob}
  \alpha_{\Lambda'|\Lambda}^{n'|n}:\;\;\hbox{Def}_{\Lambda,n} 
\To{} \hbox{Def}_{\Lambda',n'}
%\end{equation}
$$
indic\'ee par les quadruplets $(\Lambda,\Lambda',n,n')$ o\`u
  $\Lambda, \Lambda'$ sont des r\'eseaux de $K^d$ et $n,n'\in\NM$ sont
  tels que $n-n'\geq d(\Lambda,\Lambda')$, et telle que 
\begin{enumerate}
\item si $\Lambda=\Lambda'$, alors $
  \alpha_{\Lambda'|\Lambda}^{n'|n}$ est le morphisme de restriction de
  la structure de niveau.
\item pour tout autre $\Lambda'', n''$ avec $n'-n''\geq
    d(\Lambda',\Lambda'')$, on a $\alpha_{\Lambda''|\Lambda}^{n''|n}=
    \alpha_{\Lambda''|\Lambda'}^{n''|n'} \circ
    \alpha_{\Lambda'|\Lambda}^{n'|n}$.
\item les fibres g\'en\'eriques des morphismes de sch\'emas formels
  $\wh\MC_{LT,\Lambda,n} \To{}\wh\MC_{LT,\Lambda',n'}$ associ\'es sont
  \'etales, finies et surjectives (au sens de Berkovich, par exemple).
\end{enumerate}
\end{lemme}

\`A partir de ce lemme, on 
 d{\'e}finit l'action de $G$ de la mani\`ere suivante :
tout $g\in G$ induit un
isomorphisme  $\hbox{Def}_{\Lambda,n} \simto
\hbox{Def}_{g\Lambda,n}$ qui  envoie le triplet $(X,\rho,\psi)$ sur le
triplet $(X,\rho,\psi\circ g^{-1})$. On pose alors pour tous $n,n' \in \NM$
``ad{\'e}quates'', c'est-\`a-dire tels que $n-n'\geq d(g\OC^d,\OC^d)$
\ini\begin{equation}
  \label{acG}
 g^{n'|n} :\;\; \hbox{Def}_{\OC^d,n} \simto \hbox{Def}_{g\OC^d,n}
 \To{\alpha_{\OC^d|g\OC^d}^{n'|n}} \hbox{Def}_{\OC^d,n'}.
\end{equation}
En passant aux fibres g{\'e}n{\'e}riques, on en d\'eduit les $g^{n'|n}$ de \ref{resumeLT}.
Notons que  leurs degrés ne dépendent pas de $g$ et sont égaux à $[1+\varpi^{n'}M_d(\OC):1+\varpi^{n}M_d(\OC)]=q^{d^2(n-n')}$.

\medskip 

\begin{proof} 
Remarquons que dans l'\'egalit\'e ii), le quadruplet
$(\Lambda,\Lambda'',n,n'')$ v\'erifie bien $n-n''\geq
d(\Lambda,\Lambda'')$ par l'in\'egalit\'e triangulaire.
 
Fixons maintenant $(\Lambda,\Lambda',n,n')$ et choisissons $(k,r)$
v\'erifiant \ref{condreseaux} et tels que $n-n' \geq k+r$.
 %  Alors pour tous entiers $n,n'$ tels que $n\geq n'+k+r$, on
%  d{\'e}finit un morphisme 
%\ini\begin{equation} \label{blob}
%    \alpha_{\Lambda'|\Lambda}^{n'|n}:\;\;\hbox{Def}_{\Lambda,n} \To{}
%    \hbox{Def}_{\Lambda',n'}
%\end{equation}
On d\'efinit un morphisme $(\alpha_{\Lambda'|\Lambda}^{n'|n})_{k,r}$
comme dans l'\'enonc\'e du lemme
en associant {\`a} $(X,\rho,\psi)$ au-dessus de $R\in\CC$ le triplet
$(X'_{k,r},\rho'_{k,r},\psi'_{k,r})$ au-dessus de $R$ d{\'e}fini par
\begin{itemize}
\item $X'_{k,r}:=X/\psi(\varpi^{-r}\Lambda'/\Lambda)$, dont l'existence (et
  la d{\'e}finition pr{\'e}cise) est assur{\'e}e par \cite[Lemma
  4.4]{Drinell}, ou \cite[Lemma II.2.4]{HaTay}.
\item $\rho'_{k,r} :\XM \To{\varpi^{-r}} \XM \To{\rho} X\otimes_R R/\MG_R
  \To{can} X'\otimes_R R/\MG_R $.
\item $\psi'_{k,r} :\;\; \varpi^{-n'} \Lambda'/\Lambda' \To{\times
    \varpi^{-r}} \varpi^{-r-n'} \Lambda'/\varpi^{-r}\Lambda' \injo
  \xymatrix{ \varpi^{-n}\Lambda/\varpi^{-r}\Lambda' \ar@{-->}[r] &
    (\MG_R,X') \\ \varpi^{-n}\Lambda/\Lambda \ar@{->>}[u]
    \ar[r]^{\psi} & (\MG_R,X) \ar[u]^{can}}$
\end{itemize}
Il est clair que ce morphisme est $\dd$-\'equivariant.
Il nous faut voir qu'il ne d\'epend pas du couple $(k,r)$ tel que
$n-n'\geq k+r$. Soit $(k',r')$ un autre tel couple, et supposons que
$r' \geq r$.
%L'ind\'ependance de $k$ est claire sur la
%d\'efinition. Pour voir l'ind\'ependance de $r$,  il faut
 % remarquer que lorsque $r'\geq r$, 
Dans ce cas la multiplication par
  $\varpi^{r'-r}$ dans le $\OC_K$-module formel
  $X/\psi(\varpi^{-r}\Lambda'/\Lambda)$ induit un isomorphisme
  $X/\psi(\varpi^{-r'}\Lambda'/\Lambda) \simto
  X/\psi(\varpi^{-r}\Lambda'/\Lambda)$. On v\'erifie alors ais\'ement
  que cet isomorphisme induit un
  isomorphisme de triplets $(X_{k',r'}',\rho'_{k',r'},\psi'_{k',r'})
  \simto (X_{k,r}',\rho'_{k,r},\psi'_{k,r})$.

La propri\'et\'e i) est imm\'ediate. Pour la propri\'et\'e ii),
choisissons $(k,r)$ comme ci-dessus et $(k',r')$ tel que $n'-n''\geq
k'+r'$ et $\Lambda'\subseteq \varpi^{-r'} \Lambda''\subseteq \varpi^{-r'-k'}\Lambda'$. Alors
soit $(k'',r''):=(k'+k,r'+r)$, on a bien $n-n''\geq k''+r''$ et
$\Lambda\subseteq \varpi^{-r''} \Lambda''\subseteq
\varpi^{-r''-k''}\Lambda$ et on v\'erifie sur la d\'efinition que 
$(\alpha_{\Lambda''|\Lambda}^{n''|n})_{k'',r''}=
    (\alpha_{\Lambda''|\Lambda'}^{n''|n'})_{k',r'} \circ
    (\alpha_{\Lambda'|\Lambda}^{n'|n})_{k,r} $.
Il reste \`a montrer l'\'etale finitude et la surjectivité des morphismes induits sur les
fibres g\'en\'eriques. Pour cela, on utilise le fait que lorsque on a
deux morphismes analytiques  tels que $g\circ f$ est \'etale, alors
($g$  \'etale) $\Rightarrow$ ($f$ étale), et ($f$ étale surjectif) $\Rightarrow$  ($g$ étale). Rappelons maintenant que les
morphismes de  restriction de la structure de niveau sont \'etales surjectifs (et
finis) en fibre g\'en\'erique. Ainsi, \'etant donn\'e un quadruplet
$(\Lambda,\Lambda',n,n')$ tel que $n-n'\geq d(\Lambda,\Lambda')=:\delta$,
l'\'egalit\'e $\alpha_{\Lambda'|\Lambda'}^{0|n'}\alpha_{\Lambda'|\Lambda}^{n'|n}=
\alpha_{\Lambda'|\Lambda}^{0|\delta}\alpha_{\Lambda|\Lambda}^{\delta|n}$
montre que $\alpha_{\Lambda'|\Lambda}^{n'|n}$ est \'etale  en
fibre g\'en\'erique \ssi\ $\alpha_{\Lambda'|\Lambda}^{0|\delta}$
l'est, et donc \ssi\ $\alpha_{\Lambda'|\Lambda}^{2\delta|3\delta}$
l'est !
Ce dernier est le premier morphisme de la suite
$$ \hbox{Def}_{\Lambda,3\delta} \To{}\hbox{Def}_{\Lambda',2\delta}
\To{}\hbox{Def}_{\Lambda,\delta} \To{} 
\hbox{Def}_{\Lambda',0}. $$
En fibre g\'en\'erique, le caract\`ere \'etale de la compos\'ee des deux
derni\`eres fl\`eches
$\alpha_{\Lambda'|\Lambda}^{0|\delta} \circ
\alpha_{\Lambda|\Lambda'}^{\delta|2\delta} =
\alpha_{\Lambda'|\Lambda'}^{0|2\delta}$ implique que celle du milieu,
$\alpha_{\Lambda|\Lambda'}^{\delta|2\delta}$, est non-ramifi\'ee. Mais
alors, toujours en fibre g\'en\'erique, le caract\`ere \'etale et fini de 
 la compos\'ee des deux premi\`eres fl\`eches implique par \cite[3.2.9]{Bic2} que
$\alpha_{\Lambda'|\Lambda}^{2\delta|3\delta}$ est \'etale et finie.
La surjectivité se voit de la m\^eme mani\`ere.
\end{proof}

%{\em Remarque :} Les $\alpha_{\Lambda'|\Lambda}^{n'|n}$ {\em ne commutent
%  pas}, en g\'en\'eral, avec les morphismes d'oubli $\hbox{Def}_{\Lambda,n}\To{}
%\hbox{Def}$.

%\medskip
\medskip

La donn\'ee de descente se d\'efinit comme pour la tour de Drinfeld (voir \cite[3.48]{RZ}), et ici encore se recolle avec l'action de l'inertie sur les 
$$\mltn^{{ca}}:=\mltn\wh\otimes_\knr \ka,\;\; n\in \NM.$$
%On a donc obtenu un syst\`eme projectif $\MC^{{ca}}_{LT}$ de

\alin{Morphisme de p{\'e}riodes} \label{perLT}
Dans le cas $LT$, le morphisme de p{\'e}riodes n'est d{\'e}fini qu'en
fibre g{\'e}n{\'e}rique. Il a {\'e}t{\'e} d'abord d{\'e}fini dans
\cite[par. 23]{HG} {\`a} partir de calculs explicites. La d{\'e}finition
la plus visiblement intrins{\`e}que 
est celle de  \cite[ch. 5]{RZ} mais elle fait appel aux techniques
cristallines des groupes $p$-divisibles et n'est {\'e}crite qu'en
in{\'e}gales caract{\'e}ristiques\footnote{Cependant, A. Genestier sait
  adapter au cas d'{\'e}gales caract{\'e}ristiques, par exemple en
  utilisant le module de coordonn{\'e}es}.

Dans chacune des d{\'e}finitions, la construction repose sur l'existence
d'un $\knr$-espace vectoriel $\MM$ de dimension $d$  muni d'une
action de $\dd$ et  attach{\'e} au $\OC$-module formel $\XM$ (son
module de Dieudonn{\'e} dans le cas d'in{\'e}gales caract{\'e}ristiques ou
son module de coordonn{\'e}es en {\'e}gales caract{\'e}ristiques) et d'un
isomorphisme $\dd$-{\'e}quivariant de $\OC_{\mlto}$-modules
$$ \MM \otimes_\knr \OC_{\mlto} \simto M_{X_u} $$
o{\`u} $M_{X_u}$ est un fibr{\'e} vectoriel au-dessus de $\mlto$ attach{\'e}
{\`a}  l'objet
universel $X_u$ (obtenu comme  fibre g{\'e}n{\'e}rique de l'alg{\`e}bre de Lie
de l'extension vectorielle universelle de l'objet
universel $X_u$ au-dessus de 
$\wh\MC_{LT,0}$, resp. de son module de coordonn{\'e}es en {\'e}gales
caract{\'e}ristiques). Notons alors $\hbox{Lie}_{X_u}$ le 
$\OC_{\mlto}$-module inversible obtenu comme  fibre 
g{\'e}n{\'e}rique de l'alg{\`e}bre de Lie de $X_u$
au-dessus de $\wh\MC_{LT,0}$ ; la compos{\'e}e
$$ 
 \MM \otimes_\knr \OC_{\mlto} \simto M_{X_u} \To{} \hbox{Lie}_{X_u}$$
d{\'e}finit un morphisme de $\knr$-espaces analytiques
$$ \xi^{nr}:\;\;\mlto \To{} \PM(\MM) \simeq \PM^{d-1,nr}$$
qui est le morphisme de p{\'e}riodes voulu.
%\begin{fact} \label{equperLT}
%  Le morphisme $\xi$ est {\'e}tale surjectif et $\dd\times W_K$-{\'e}quivariant. Il est aussi
%  $G$-{\'e}quivariant au sens suivant : pour tout $g\in G$ et tous $n,n'$
%  ad{\'e}quates, le diagramme
%$$\xymatrix{ \mltn \ar[r]^{g^{n'|n}} \ar[d]_{\xi_n} &
%  \MC_{LT,n'} \ar[ld]^{\xi_{n'}} \\  \PM(\MM)
%  & \\ }$$
%est commutatif (bien-s{\^u}r on a not{\'e} $\xi_n$ la compos{\'e}e de $\xi$
%avec le morphisme d'oubli de la structure de niveau).
%\end{fact}
%\begin{proof}
Le fait que $\xi^{nr}$ est étale au sens rigide-analytique est prouvé dans
\cite[5.17]{RZ}. Il se trouve qu'il est aussi étale au sens de
Berkovich car son bord relatif est vide. La surjectivité de $\xi^{nr}$ est prouvée, au moins pour les points classiques (rigides analytiques)  dans
\cite[23.5]{HG}, mais la preuve de {\em loc.cit} montre aussi la surjectivité pour les points de Berkovich. Les propriétés d'équivariance se trouvent dans
\cite[5.37]{RZ} ou \cite[23.28]{HG}.   
%\end{proof}

Par ailleurs, les isomorphismes $\MM(\XM) \To{Frob}
\MM(\phi^*\XM)=\phi^*\MM(\XM)$ munissent $\MM$ et donc $\PM(\MM)$
d'une donn\'ee de descente \`a la Weil compatible \`a $\xi^{nr}$.
%structure $\phi^\ZM$-{\'e}quivariante au-dessus de $\knr$ (ou
%donn{\'e}e de descente {\`a} la Weil). 
%Sur $\MM$ elle n'est pas effective, mais elle l'est sur
%$\PM(\MM)$ : notons en effet $S_K^{1/d}$ le $K$-espace analytique
%obtenu par analytification de la $K$-vari{\'e}t{\'e} de Severi-Brauer
%associ{\'e}e {\`a} l'alg{\`e}bre centrale simple sur $K$ d'invariant
%$1/d$. Alors 
%\begin{fact}
 % Il existe un isomorphisme $\dd\times W_K$-{\'e}quivariant de
 % $\knr$-espaces analytiques $$\PM(\MM) \simto S^{1/d}_K\wh\otimes_K
 % \knr.$$
%\end{fact}
%\begin{proof}
Explicitons-l\`a dans le cas où $K$ est de caractéristique nulle : alors le $\OC_K$-groupe formel $\XM$ sur $\o\FM_p$ est $p$-divisible et $\MM$ est l'isocristal associé. Sa dimension sur $\wh{K^{nr}}$ est la hauteur $d$ de $\XM$, et son unique pente est $1/d$. En particulier $\MM$ est irréductible et donc ``défini sur $\FM_p$" au sens suivant : il existe une base de $\MM$ dans laquelle le Frobenius est donné par la matrice 
$\Phi=\left(\begin{array}{cccc} 0 & 1 & 0 & 0 \\ 0 & \ddots & \ddots  & 0 \\ 0 & 0 & \ddots & 1 \\ \varpi & 0 & 0 & 0 
 \end{array}\right)$. Ainsi, la forme rationnelle de $\PM^{d-1}$ correspondant à la donnée de descente sur $\PM(\MM)$ est celle associée au $1$-cocycle $\gal(\wh{K^{nr}}/K)=\phi^\ZM \To{} \aut{}{\PM^{d-1}}=PGL_{d-1}(\wh{K^{nr}})$ qui envoie $\phi$ sur la matrice $\Phi$. Cette forme rationnelle est la variété de Severi-Brauer d'invariant $1/d$.

Dans le cas où $K$ est de caractéristique positive, le même raisonnement s'applique au module de coordonnées $\MM$.
%\end{proof}

%\alin{Cohomologie} \label{cohoLT}

\subsection{D{\'e}finition  du
  $R\Gamma_c$} \label{defRG} 

Dans cette section, nous  d{\'e}finissons les complexes de
cohomologie des tours $\mlt^{{ca}}$ et $\mdr^{{ca}}$. On les obtient comme
{\'e}valuation en un faisceau constant du foncteur d{\'e}riv{\'e} d'un
certain foncteur ``sections {\`a} 
support compact sur la tour''.  Comme on veut que ce foncteur soit {\`a} valeurs dans la
cat{\'e}gorie des $G\times \dd\times W_K$-modules, le plus facile (surtout
du c{\^o}t{\'e} $LT$)  est de
prendre pour source de ce foncteur la cat{\'e}gorie des faisceaux
{\'e}quivariants sur (la descente de) l'espace des p{\'e}riodes. 
Nous utiliserons sans commentaires, mais avec des r\'ef\'erences pr\'ecises, le formalisme pr\'esent\'e dans l'appendice B de \cite{Dat1}, et qui est essentiellement d\^u \`a Berkovich et Jannsen. Pour le confort du lecteur, rappelons sans les d\'efinir les principales notations de cet appendice :
\begin{itemize}
        \item $\wt{X_{et}}$ d\'esigne le topos \'etale d'un espace analytique $X$ et $\Gamma_!(X,-)$ les sections \`a support compact.
        \item $\wt{X_{et}}(G)$ d\'esigne le topos form\'e par les faisceaux \'etales $G$-\'equivariants d'un espace alg\'ebrique $X$ muni d'une action {\em continue} d'un groupe topologique $G$ \cite[B.1.4]{Dat1}. On sait alors que le foncteur des sections \`a support compact se factorise par  la cat\'egorie (le topos) $\wt{G}$ des ensembles discrets munis d'une action continue de $G$ et nous notons $\Gamma_!^{\infty_G}(X,-) : \wt{X_{et}}(G) \To{} \wt{G}$  le foncteur obtenu.
        \item $\Mo{\Lambda}{\TC}$ d\'esigne la cat\'egorie ab\'elienne des $\Lambda$-modules d'un topos $\TC$ \cite[B.1.5]{Dat1}.
        \item Si $\Lambda$ est une extension finie (sous-entendu locale et plate) de $\ZM_l$ d'uniformisante $\lambda$, on note $\Lambda_\bullet$ le pro-anneau $(\Lambda/\lambda^n)_{n\in\NM}$ et
$\Mo{\Lambda_\bullet}{\TC}$ d\'esigne la cat\'egorie ab\'elienne des $\Lambda_\bullet$-modules du topos $\TC$ ({\em i.e.} des syst\`emes projectifs $(\FC_n)_{n\in\NM}$ de $\Lambda$-faisceaux sur $\TC$  tels que $\lambda^n\FC_n =0$). On a un foncteur $\limproj^{\TC} :\Mo{\Lambda_\bullet}{\TC}\To{} \Mo{\Lambda}{\TC}$ \cite[B.2.2]{Dat1}.
        \item On note $\Gamma_{!,\bullet}^{\infty_G} := \Gamma_!^{\infty_G} \circ \limproj^{\wt{X_{et}}(G)} : \Mo{\Lambda_\bullet}{\wt{X_{et}}(G)} \To{} \Mo{\Lambda}{G}$ \cite[B.2.4]{Dat1}. 
\end{itemize}

%\`A ce propos, lorsqu'un groupe topologique agit contin\^ument sur un espace de Berkovich nous dirons simplement ``faisceau \'equivariant" pour ``faisceau \'etale \'equivariant {\em discret}" (au sens de Berkovich), et nous renvoyons \`a \cite[B.1.4]{Dat1} pour des explications d\'etaill\'ees.

\alin{Sections \`a supports compacts sur la tour} \label{cohoDr}
Comme les constructions sont formellement identiques pour chacune des deux tours, 
%s'appliquant de la m{\^e}me mani{\`e}re {\`a} la tour de 
%Drinfeld $\mdr^{{ca}}$  ou {\`a} celle de Lubin-Tate $\mlt^{{ca}}$, 
nous ommettrons parfois la notation $Dr$ ou $LT$.
Par exemple,
 pour tout $n\in \NM$, $\xi_n$ d{\'e}signera soit la compos{\'e}e
$$ \xi_{Dr,n} :\;\; \mdrn^{{ca}} \To{} \mdro^{{ca}} \To{} \Omega^{d-1,nr}_K \To{}
\PC_{Dr}:=\Omega^{d-1}_K $$
soit la compos{\'e}e
$$ \xi_{LT,n} :\;\; \mltn^{{ca}} \To{} \mlto^{{ca}} \To{} \PM^{d-1,nr}_K \To{}
\PC_{LT}:=S^{1/d}_K. $$
Il s'agit de morphismes d'espaces ``analytiques au-dessus de $K$'' au
sens de \cite{Bic2}, les deux premiers espaces \'etant
$\ka$-analytiques, le troisi\`eme $\knr$-analytique, et le dernier $K$-analytique. 
Dans chacun des cas, l'espace des p{\'e}riodes $\PC$ est muni d'une action
 continue du groupe $J$ en posant  $J_{Dr}:=G$ et $J_{LT}:=\dd$. 
 % et d'une action triviale de $W_K$. 
 %Le morphisme $\xi_{n}$ est
 %formellement {\'e}tale et $J$-{\'e}quivariant.
Comme les morphismes $ J$-\'equivariants $\MC_n^{{ca}} \To{\pi_{n,n'}}
\MC_{n'}^{{ca}}$  sont  finis, on a $\pi_{n,n',!}=\pi_{n,n',*}$. Donc  pour tout objet $\FC \in \wt{\PC_{et}}(J)$,  
on a des applications naturelles $J$-\'equivariantes 
$$ \pi_{n,n'}^*:\;\;\Gamma_!^{\infty_J}(\MC_{n'}^{{ca}},\xi_{n'}^*(\FC)) \To{}
\Gamma_!^{\infty_J}(\MC_{n'}^{{ca}},\pi_{n,n',!}\pi_{n,n'}^*\xi_{n'}^*(\FC)) =
\Gamma_!^{\infty_J}(\MC_n^{{ca}},\xi_n^*(\FC))$$
% Dans le cas $Dr$,
 %il est aussi $\dd$-{\'e}quivariant pour l'action triviale de $\dd$ sur
 %$\Omega^{d-1}_K$ et l'action d{\'e}finie en \ref{actionDr} sur les $\mdrn^{{ca}}$. Dans le
 %cas $LT$, c'est plut{\^o}t le syst{\`e}me inductif des $\xi_{LT,n}$ qui est
 %$G$-{\'e}quivariant au sens de \ref{equperLT}.
%Si $\FC$ est un faisceau ab\'elien de torsion sur $\wt{\PC_{et}}(J)$, alors $\xi_n^*(\FC)$ est un faisceau $W_K$-\'equivariant sur $\wt{{\MC_{n,et}^{{ca}}}}(J)$ et le groupe ab\'elien des sections \`a support compact $\Gamma_!(\MC_n^{{ca}},\xi_n^*(\FC))$ est muni d'une action de $J\times W_K$, lisse sur $J$, que nous notons $\Gamma_!^{\infty_J}(\MC_n^{{ca}},\xi_n^*(\FC))$, suivant les conventions de\cite[B.1.5]{Dat1}.

%et on peut donc poser
%$$ \Gamma_c(\MC^{{ca}},\FC):= \limi{n\in \NM}
%\Gamma_!^{\infty_J}(\MC_n^{{ca}},\xi_n^*(\FC)) $$
%qui est un ensemble  muni d'une action continue (lisse) de $J \times W_K$.
\begin{fact} \label{actiontours}
  L'ensemble $\limi{n\in \NM}
\Gamma_!^{\infty_J}(\MC_n^{{ca}},\xi_n^*(\FC))$ est muni d'une action continue
  de $G\times \dd\times W_K$ %continue (lisse) sur les deux premiers facteurs et
  et triviale  sur le sous-groupe ${K^\times_{diag}}\times 1$. % (rappelons que ${K^\times_{diag}}$ a \'et\'e
%d\'efini en \ref{defDelta}).
\end{fact}
\begin{proof} Dans les deux cas, on a visiblement une action continue de $(J\times W_K)$. 
Dans le cas $Dr$, on a aussi imm\'ediatement l'action de $\dd$.
On s'occupe donc du cas $LT$, o\`u il reste \`a expliquer l'action de $G$.

Fixons $g\in G$ et $n \geq n'$ tels que
le morphisme $g^{n'|n} : \mltn^{{ca}} \To{} \MC_{LT,n'}^{{ca}}$
soit d\'efini. Celui-ci \'etant fini, on a 
$(g^{n'|n})_!=(g^{n'|n})_*$, donc l'\'egalit\'e $\xi_{LT,n'}\circ g^{n'|n} = \xi_{LT,n}$ de \ref{resumeLT} ii) fournit le morphisme
\begin{eqnarray*}
g^{n'|n,*}:\;\;  \Gamma_!^{\infty_J}(\MC_{LT,n'}^{{ca}},\xi_{LT,n'}^*(\FC)) &  \To{} &
%\Gamma_!^{\infty_J}(\MC_{LT,n'}^{{ca}},({g^{n'|n}})_! (g^{n'|n})^*
%\xi_{LT,n'}^*(\FC)) \\&=  & 
\Gamma_!^{\infty_J}(\MC_{LT,n}^{{ca}},\xi_{LT,n}^*(\FC)).
\end{eqnarray*}
 D'apr\`es les propri\'et\'es \ref{resumeLT} i)(b) i et ii, 
%Les propri\'et\'es des morphismes $g^{n'|n}$ mentionn\'ees lors de leur
%construction en \ref{actionLT} assurent la commutativit\'e des
%diagrammes du type 
%$$\xymatrix{\Gamma_!^{\infty_J}(\MC_{LT,n'}^{{ca}},\xi_{LT,n'}^*(\FC)) \ar[d]_{\pi_{m',n'}^*}
%  \ar[r]^{g^{n'|n,*}} &
%\Gamma_!^{\infty_J}(\MC_{LT,n}^{{ca}},\xi_{LT,n}^*(\FC)) \ar[d]^{\pi_{m,n}^*} \\
%\Gamma_!^{\infty_J}(\MC_{LT,m'}^{{ca}},\xi_{LT,m'}^*(\FC)) \ar[r]^{g^{m'|m,*}} &
%\Gamma_!^{\infty_J}(\MC_{LT,m}^{{ca}},\xi_{LT,m}^*(\FC))
%}$$
%o\`u $m\geq n$ et $m'\geq n'$ sont tels que $g^{m'|m}$ est
%d\'efini. Ainsi 
on obtient \`a la limite inductive un morphisme
$g^*:\; \Gamma_c(\MC^{{ca}},\FC) \To{} \Gamma_c(\MC^{{ca}},\FC)$ qui commute \`a l'action de $\dd\times W_K$, Cela
%, et si $g'$
%est un autre \'el\'ement de $G$, alors $(gg')^*={g'}^*g^*$ de sorte
%qu'on obtient une action {\em \`a droite} de $G$ sur
%$\Gamma_!^{\infty_J}(\MC^{{ca}},\FC)$. 
d\'efinit une action {\em \`a droite} de $G$, et 
{\em pour nous ramener \`a une action \`a gauche,
 nous ferons agir $g$ sur $\Gamma_!^{\infty_J}(\MC^{{ca}},\FC)$ par son inverse
${g}^{-1}$.}
La trivialit\'e de cette action sur $K^\times_{diag}$ vient de la fin de \ref{resumeLT} i).
%Soit maintenant $z\in K^\times$ et $z_G$, resp. $z_D$, son image dans
%le centre de $G$, resp. dans celui de $\dd$. 
%Alors les actions g\'eom\'etriques de
%$z_{G}$ et $z_D$ sur la tour stabilisent chaque
%$\MC^{{ca}}_{LT,n}$ et sont inverses l'une de l'autre. 
%Il en est donc de m\^eme sur les sections $\Gamma_!^{\infty_J}(\MC^{{ca}}_{LT,n},\FC)$,
%d'o\`u la
%trivialit\'e de l'action de ${K^\times_{diag}}\times 1$.
\end{proof}

Posons $GD:=(G\times \dd)/K^\times_{diag}$. On a donc obtenu un foncteur
``sections \`a supports compacts sur la tour" 
$$ \application{ \Gamma_c(\MC^{ca},-):\;}{\wt{\PC_{et}}(J)}{\wt{(GD\times W_K)}}
{\FC}{\limi{n\in \NM}
\Gamma_!^{\infty_J}(\MC_n^{{ca}},\xi_n^*(\FC))}
$$
dont le but est la cat\'egorie des $GD\times W_K$-ensembles continus. Ce foncteur commute aux limites projectives finies et envoie donc $\Lambda$-modules sur $\Lambda$-modules pour tout anneau $\Lambda$ (constant sur $\wt{\PC_{et}}(J)$).

\alin{Cohomologie} \label{rgamma}
%Si $\Lambda$ est un anneau de torsion, on d{\'e}finit ainsi un foncteur
%exact {\`a} gauche 
%$$\Gamma_c(\MC^{{ca}},-):\;\; \Mo{\Lambda}{\wt{\PC_{et}}(J)} \To{}
%\Mo{\Lambda}{GD\times W_K},$$
%o\`u le dernier terme d\'esigne la cat{\'e}gorie  des
%$\Lambda$-modules munis d'une action lisse de $GD=(G\times \dd)/{K^\times_{diag}}$ et d'une action
%commutante de $W_K$.
La construction du syst\`eme inductif 
$(\Gamma_!^{\infty_J}(\MC^{{ca}}_n,\xi_n^*(\FC)))_{n\in\NM}$ et de l'action de
$G\times \dd\times W_K$ sur icelui reposaient sur la finitude 
 des morphismes $\pi_{n,n'}$ et $g^{n'|n}$.
Ces m\^emes propri\'et\'es permettent de d\'efinir de mani\`ere analogue pour tout faisceau ab\'elien et tout $q\in\NM$ un syst\`eme inductif
$(H^q_c(\MC^{{ca}}_n,\xi_n^*(\FC)))_{n\in\NM}$ muni d'une action de
$G\times \dd\times W_K$.

\begin{fact}
 Pour tout $q\in\NM$, il y a un isomorphisme de foncteurs $\Mo{\ZM}{\wt{\PC_{et}}(J)} \To{} \Mo{\ZM}{GD\times W_K}$ canonique :
$$ R^q\Gamma_c(\MC^{{ca}},-) \simto \limi{n\in\NM} H^q_c(\MC_n^{{ca}},\xi_n^*-).$$
\end{fact}
\begin{proof}
En effet, par commutation de la cohomologie aux limites inductives
filtrantes, on a $$R^q\Gamma_c(\MC^{{ca}},-) =\limi{}
R^q\left(\Gamma_!^{\infty_J}(\MC_n^{{ca}},-)\circ \xi_n^*\right).$$ 
Comme on sait par Berkovich que $R^q\Gamma_!^{\infty_J}(X,-) = H^q_c(X,-)$ (voir \cite[B.1.7]{Dat1}), il nous suffira de montrer que $R^q(\Gamma_!^{\infty_J}(\MC_n^{{ca}},-)\circ
\xi_n^*)=R^q\Gamma_!^{\infty_J}(\MC_n^{{ca}},-)\circ \xi_n^*$ pour tout $n$. Or on a une factorisation $\xi_n^* : \Mo{\ZM}{\wt{\PC_{et}}(J)}  \To{\beta^*} \Mo{\ZM}{\wt{\PC{\wh\otimes}_K \wh{K^{ca}}}_{et}(J)} \To{{\xi_n^{ca}}^*} \Mo{\ZM}{\wt{\MC_{n,et}^{{ca}}}(J)}$ o\`u $\xi_n^{ca}$ est un morphisme \'etale de $\wh{K^{ca}}$-espaces analytiques et $\beta$ est le changement de base (ce n'est pas un morphisme d'espaces analytiques mais un morphisme des sites \'etales concern\'es). Ainsi ${\xi_n^{ca}}^*$ envoie injectifs sur injectifs puisqu'il est adjoint \`a droite du foncteur exact ${\xi_n^{ca}}_!$, et $\beta^*$ envoie suffisamment d'injectifs sur des injectifs, par exemple tous ceux de la forme $\beta_*(\IC)$, avec $\IC$ injectif.
\end{proof}

En particulier, lorsque $\Lambda$ est un anneau de torsion, le complexe $R\Gamma_c(\MC^{ca},\Lambda)$ a la bonne cohomologie.

%\alin{Coefficients $l$-adiques} \label{defrgadic}

Supposons maintenant que $\Lambda$ est une extension finie de $\ZM_l$ et fixons 
$\FC_\bullet=(\FC_m)_{m\in \NM} \in \Mo{\Lambda_\bullet}{\wt{\PC_{et}}(J)}$.
Comme les morphismes $\pi_{n,n'}$ pour $n\geq n'$ sont \'etales, on a 
$$ \mathop{\limproj}\nolimits^{\wt{\MC_{n,et}}(J)}(\xi_n^{*}\FC_\bullet) \simeq \pi_{n,n'}^*\mathop{\limproj}\nolimits^{\wt{\MC_{n',et}}(J)}(\xi_{n'}^{*}\FC_\bullet). $$
Ainsi, la finitude de ces m\^emes $\pi_{n,n'}$ %permet
% donc comme dans le cas
%de torsion de d\'efinir des morphismes de $\Lambda$-modules
%$$ \pi_{n,n'}^*:\;\; \Gamma_{!,\bullet}^{\infty_J}(\MC^{{ca}}_{n'},\xi_{n'}^{*}(\FC_\bullet)) \To{}
%\Gamma_{!,\bullet}^{\infty_J}(\MC^{{ca}}_n,\xi_{n}^{*}(\FC_\bullet)) $$
%et de poser 
fournit encore les morphismes de transition intervenant dans la limite suivante
$$ \Gamma_{c,\bullet}(\MC^{{ca}},\FC_\bullet):= \limi{n\in \NM}
\Gamma_{!,\bullet}^{\infty_J}\left(\MC_n^{{ca}},\xi_{n}^{*}(\FC_\bullet)\right).$$
Ce $\Lambda$-module est muni par d\'efinition d'une action {\em lisse} de $J$
\`a gauche et d'une action commutante  de $W_K$ \`a droite (qui n'est plus n\'ecessairement lisse \`a cause de la limite projective).
On construit l'action du troisi\`eme groupe ($\dd$ pour le cas $Dr$ ou
$G$ pour le cas $LT$) exactement comme dans \ref{actiontours}. %le cas de torsion.
 En particulier dans le cas $LT$, on utilise le caract\`ere \'etale des $g^{n'|n}$  pour pouvoir \'ecrire
$$ \mathop{\limproj}\nolimits^{\wt{\MC_{n,et}}(J)}(\xi_n^{*}\FC_\bullet) \simeq {g^{n'|n}}^*\mathop{\limproj}\nolimits^{\wt{\MC_{n',et}}(J)}(\xi_{n'}^{*}\FC_\bullet). $$
gr\^ace \`a \ref{resumeLT} ii).
On d\'efinit ainsi un foncteur 
$$ \Gamma_{c,\bullet}(\MC^{{ca}},-):\;\; \Mo{\Lambda_\bullet}{\wt{\PC_{et}}(J)} \To{}
\Mo{\Lambda}{GD\times W_K^{disc}} $$
qui est exact {\`a} gauche. Dans le terme de droite, l'exposant $disc$ signifie qu'on oublie la topologie de $W_K$ ; en effet, on a perdu la lissit\'e de l'action de $W_K$ dans la limite projective.
% d\'esigne la cat\'egorie des $\Lambda$-modules  munis d'une action lisse $GD$ et d'une action de $W_K$ (). 

Supposons maintenant que le $\Lambda_\bullet$-faisceau $\FC_\bullet$ soit un
$\Lambda$-syst\`eme local au sens de \cite[B.2.1]{Dat1}, et notons 
$H^q_c(\MC_n^{{ca}},\xi_n^*((\FC_m)_m))$ la cohomologie \`a supports
compacts de son image inverse sur $\MC_n^{{ca}}$, d\'efinie dans ce cas par Berkovich ({\em cf loc. cit}).
Toujours les m\^emes propri\'et\'es de propret\'e des $\pi_{n,n'}$ et
$g^{n'|n}$ permettent de d\'efinir un syst\`eme inductif
$(H^q_c(\MC_n^{{ca}},\xi_n^*((\FC_m)_m)))_{n\in\NM}$ muni d'une
action du groupe $G\times\dd\times W_K$.

\begin{fact}
  Pour tout $q\in \NM$ et tout $\Lambda$-syst{\`e}me local $J$-{\'e}quivariant
$(\FC_m)_{m\in\NM}$ sur $\PC$,
%  $J$-{\'e}quivariant discret/lisse,
 on a un isomorphisme canonique
  $(G\times \dd\times W_K)$-{\'e}quivariant :
$$ R^q\Gamma_{c,\bullet}(\MC^{{ca}},\FC_\bullet) \simto \limi{n\in \NM}
H^q_c\left(\MC_n^{{ca}},\xi_n^*((\FC_m)_m)\right) $$
\end{fact}
\begin{proof}
  Comme dans le cas de torsion, la commutation de la cohomologie aux limites
  inductives filtrantes et le fait que $\xi_n^*$ envoie suffisamment d'injectifs sur des injectifs montrent que 
\begin{eqnarray*}
  R^q\Gamma_{c,\bullet}(\MC^{{ca}},\FC_\bullet) & \simeq & \limi{n\in\NM}
\left(R^q\left(\Gamma_{!,\bullet}^{\infty_J}(\MC_n^{{ca}},-) \circ
    \xi_n^* \right)(\FC_\bullet) 
 \right) \\
& \simeq &  \limi{n\in\NM} \left(
R^q\left(\Gamma_{!,\bullet}^{\infty_J}(\MC_n^{{ca}},-) \right)
(\xi_n^*(\FC_\bullet)) \right)
\end{eqnarray*}
Mais d'apr{\`e}s \cite[Prop B.2.5 ii)]{Dat1} appliqu{\'e} aux espaces analytiques
quasi-alg{\'e}briques (car lisses) $\MC_n^{{ca}}$, on a $ R^q\Gamma_{!,\bullet}^{\infty_J}(\MC_n^{{ca}},\xi_n^*(\FC_\bullet)) \simeq
H^q_c\left(\MC_n^{{ca}},\xi_n^*(\FC_m)_m\right)$ pour tout $\Lambda$-syst{\`e}me
local $J$-\'equivariant $\FC_\bullet$
 \end{proof}

Nous pouvons maintenant poser :
$$  R\Gamma_c(\MC^{{ca}},\Lambda):=R\Gamma_{c,\bullet}(\MC^{{ca}},\Lambda_\bullet) \in D^b_{\Lambda}(GD\times W_K^{disc}).$$
%C'est un objet de la cat{\'e}gorie
%d{\'e}riv{\'e}e born{\'e}e $D^b_{\Lambda}(GD)$ de la cat\'egorie
%ab\'elienne $\Mo{\Lambda}{GD}$ des $\Lambda(GD)$-modules
%isses, qui est muni d'une action de $W_K$ et dont 
Par ce qui pr\'ec\`ede,  sa cohomologie est $G\times
\dd\times W_K$-isomorphe {\`a} la limite inductive $H^q_c(\MC^{{ca}},\Lambda)$ des
$H^q_c(\MC_n^{{ca}},\Lambda)$, $n\in\NM$. En particulier, le morphisme canonique
$R\Gamma_c(\MC^{ca},\ZM_l)\otimes_{\ZM_l}\Lambda \To{} R\Gamma_c(\MC^{ca},\Lambda)$ est un isomorphisme. Nous posons enfin 
$$ R\Gamma_c(\MC^{{ca}},\o\QM_l):= \o\QM_l\otimes_{\ZM_l}
R\Gamma_c(\MC^{{ca}},\ZM_l) \;\; \in D^b_{\o\QM_l}(GD\times W_K^{disc}). $$
Par commutation du produit tensoriel aux limites inductives, la
cohomologie de ce
complexe est donn{\'e}e par 
$ R^q\Gamma_c(\MC^{{ca}},\o\QM_l) \simeq H^q_c(\MC^{{ca}},\o\QM_l) :=
\limi{n\in\NM} H^q_c(\MC_n^{{ca}},\o\QM_l).$

\subsection{Le th\'eor\`eme de Faltings-Fargues}
Dans \cite{FaltDrin}, Faltings a esquiss\'e un lien remarquable entre les tours de Drinfeld et de Lubin-Tate. Ce lien est expliqu\'e en d\'etail et amplifi\'e dans un travail en cours de Fargues \cite{FarFal} dans le cas o\`u $K$ est $p$-adique. 
% au cas d'\'egales caract\'eristiques par L. Fargues, A. Genestier et V. Lafforgue, lors d'un groupe de travail \`a l'IHES s'\'etalant sur un semestre. 
Si on pouvait donner un sens g\'eom\'etrique raisonnable aux notations $\MC^{ca}:=\limproj \MC^{ca}_n$, la formulation naturelle de ce ``lien" serait de dire que les objets $\MC^{ca}_{Dr}$ et $\MC^{ca}_{LT}$ sont isomorphes, et de mani\`ere $G\times\dd\times W_K$-\'equivariante. 
Malheureusement ces objets n'ont vraiment pas de sens, et le r\'esultat doit \^etre formul\'e en termes de sch\'emas formels $p$-adiques \'enormes et n'ayant plus du tout les propri\'et\'es de finitude normalement requises pour leur associer des espaces rigides. Mais une autre mani\`ere de formuler ce lien sans introduire de sch\'emas formels est la suivante :
\begin{theo} (Faltings-Fargues) \label{faltings}
  Il existe une \'equivalence de topos $\wt{\PC_{Dr,et}}(G)\simeq \wt{\PC_{LT,et}}(\dd)$ qui \'echange les deux foncteurs $\Gamma_c(\MC^{{ca}}_{Dr},-)$ et $\Gamma_c(\MC^{{ca}}_{LT},-)$. 
\end{theo}
C'est le th\'eor\`eme 13.2 du chapitre ``cohomologie" de \cite{FarFal}.
Sous cette forme, on obtient le corollaire imm\'ediat
\begin{coro}
Il existe un isomorphisme $R\Gamma_c(\MC^{ca}_{LT},\ZM_l)\simeq R\Gamma_c(\MC^{ca}_{Dr},\ZM_l)$ dans $D^b_{\ZM_l}(GD\times W_K^{disc})$. %, compatible \`a l'action de $W_K$.
\end{coro}

\subsection{Variantes et décompositions}

Dans cette section, on discute deux variantes de la construction de la section précédente dont on tire quelques propriétés des complexes $R\Gamma_c(\MC^{{ca}},\Lambda)$ où $\Lambda$  désigne une $\ZM_l$-alg\`ebre finie pour un $l\neq p$. % est inversible. %, et qui sera soit de torsion, soit local complet, soit {\'e}gal {\`a} $\o\QM_l$. 
On s'intéresse ensuite à l'action du centre de la catégorie des $\Lambda(GD)$-modules lisses sur ces complexes, puis dans le cas $LT$, on montre une propriété de finitude et on introduit une variante ``sans supports".

\alin{Premi\`ere variante} \label{variante}
Dans \cite{CarAA}, \cite{HaTay} et \cite{Boyer1}, les auteurs consid{\`e}rent
plut{\^o}t la cohomologie des tours
$(\MC_{?,n}^{(0),{ca}})_{n\in\NM}$ pour $?=LT$
ou $Dr$,  o\`u l'exposant $(0)$ d\'esigne la composante connexe o\`u la quasi-isog\'enie du probl\`eme de modules est de hauteur nulle, {\em cf} \ref{modulesDr} et \ref{modulesLT}.
Ces espaces ne sont pas stables sous l'action
 $GD\times W_K$ ; leur stabilisateur est le sous-groupe 
 $(GD\times W_K)^0$ form\'e des triplets $(g,d,w)$ tels que $|\hbox{det}(g)\hbox{Nr}(d)^{-1}|_K=|w|$. En d\'efinissant les $R\Gamma_c$ correspondants comme dans le paragraphe \ref{defRG}, on constate que les inclusions $\MC_n^{(0)}\injo \MC_n$ induisent par adjonction un isomorphisme
\ini\begin{equation} \label{lien}
 \cind{(GD\times
  W_K)^0}{GD\times W_K}{ R\Gamma_c(\MC^{(0),{ca}},\Lambda)} \simto
R\Gamma_c(\MC^{{ca}},\Lambda) 
\end{equation}
dans $D^b_\Lambda(GD\times W_K^{disc})$.

\alin{Deuxi\`eme variante} \label{variante2}
On peut aussi quotienter les tours $(\MC_n)_n$ par l'action (libre) du sous-groupe $\varpi^\ZM$ engendr\'e par l'uniformisante $\varpi \in K^\times \subset G$.  Les quotients $\MC_n/\varpi^\ZM$ sont non-canoniquement isomorphes \`a
$\bigsqcup_{h=0}^{d-1} \MC_n^{(h)}$.
Comme $\varpi^\ZM$ est central dans $G\times\dd\times W_K$, ces
quotients sont encore munis d'une action du produit triple.
%\begin{nota}
Pour all\'eger un peu les notations, nous surlignerons les quotients par $\varpi^\ZM$, en posant $\o{GD}:=GD/\varpi^\ZM$  et 
$ \o\MC:=\MC/\varpi^\ZM $
que nous d\'ecorerons \'eventuellement d'indices $Dr$, $LT$ ou $n$.
%\end{nota} 
En r\'ep\'etant les constructions du paragraphe \ref{defRG}, on obtient  un complexe (deux complexes, en fait)
$R\Gamma_c(\o\MC^{ca},\Lambda) \in D^b_\Lambda(\o{GD}\times W_K^{disc})$ qui se d\'eduit de $R\Gamma_c(\MC^{ca},\Lambda)$ par la formule
\ini
\begin{equation}\label{lien1}
 \Lambda {\mathop\otimes\limits^L}_{\Lambda[\varpi^\ZM]} R\Gamma_c(\MC^{ca},\Lambda)
\simto R\Gamma_c(\o\MC^{ca},\Lambda)
\end{equation} 
(isomorphisme ``canonique'' dans $D^b_\Lambda(\o{GD}\times W_K^{disc})$) 
%o\`u intervient le foncteur (exact \`a droite) des ``$\varpi^\ZM$-coinvariants'' :
%$$\application{}{\Mo{\Lambda}{GD}}{\Mo{\Lambda}{\o{GD}}}{V}
%{\Lambda\otimes_{\Lambda[\varpi^\ZM]}V}$$ 
o\`u le produit tensoriel est pris par rapport au morphisme
d'augmentation $\Lambda[\varpi^\ZM]\To{}\Lambda$. 
%Ce foncteur est adjoint \`a gauche du foncteur d'inclusion $\omega :\;\; \Mo{\Lambda}{\o{GD}} \To{} \Mo{\Lambda}{GD}$. 
On a aussi un lien avec la variante pr\'ec\'edente :
\ini\begin{equation}\label{lien2}
 \cind{(GD\times
  W_K)^0\varpi^\ZM}{GD\times W_K}{ R\Gamma_c(\MC^{(0),{ca}},\Lambda)} \simto
R\Gamma_c(\o\MC^{{ca}},\Lambda). 
\end{equation}

\alin{D{\'e}composition selon le centre de la cat{\'e}gorie des $\Lambda
  \dd$-modules lisses} 
 Notons $\ZG_\Lambda( GD)$ le centre de la
cat{\'e}gorie ab{\'e}lienne $\Mo{\Lambda}{GD}$ des $\Lambda GD$-modules lisses. 
 Comme le centre d'une cat{\'e}gorie
ab{\'e}lienne agit encore sur la cat{\'e}gorie d{\'e}riv{\'e}e born{\'e}e
associ{\'e}e (et s'identifie m{\^e}me au centre de celle-ci), on obtient une
action canonique de $\ZG_\Lambda( GD)$ sur
$R\Gamma_c(\MC^{{ca}},\Lambda)$. En particulier tout idempotent de
$\ZG_\Lambda( GD)$ fournit un facteur direct de $R\Gamma_c(\MC^{{ca}},\Lambda)$.
De plus, en utilisant les m{\^e}mes notations pour les groupes $\dd$ et $G$, on a des morphismes
canoniques $\ZG_\Lambda(\dd) \To{} \ZG_\Lambda( GD)$ et $\ZG_\Lambda( G)
\To{} \ZG_\Lambda( GD)$ induits par
les inclusions $\dd \injo GD$ et $G\injo GD$.

La suite de sous-ensembles  $1+\varpi^n\od$, $n\in\NM$ de $\dd$ est un syst{\`e}me
fondamental de voisinages de l'unit{\'e} form{\'e} de pro-$p$-sous-groupes
ouverts compacts distingu{\'e}s de $\dd$. Comme $p$ est inversible dans
$\Lambda$, ces pro-$p$-sous-groupes ouverts d{\'e}finissent des idempotents de l'alg{\`e}bre
de Hecke $\HC_\Lambda(\dd)$. Comme ils sont distingu{\'e}s, ces
idempotents sont centraux : on note
$\varepsilon_{n,D}$ leurs images dans $\ZG_\Lambda( \dd)$ et
$\varepsilon_{n,D}R\Gamma_c(\MC^{{ca}},\Lambda) \in D^b_\Lambda( GD)$ le facteur direct 
de $R\Gamma_c(\MC^{{ca}},\Lambda)$ associ{\'e}.

{\em Interpr{\'e}tation g{\'e}om{\'e}trique :}
Du côté $Dr$, on prouve facilement (comme dans le lemme \ref{niveauLT} ci-dessous, en plus simple) que l'inclusion de foncteurs
$\Gamma_!^{\infty_{G}}(\mdrn^{{ca}},-) \circ \xi_{Dr,n}^* \injo \Gamma_c(\mdr^{{ca}},-)$
induit un isomorphisme dans $D^b_\Lambda (GD\times W_K^{disc})$
$$ R\Gamma_c(\MC_{Dr,n}^{{ca}},\Lambda) \simto \varepsilon_{n,D}R\Gamma_c(\mdr^{{ca}},\Lambda).$$
Par contre, il n'y a pas d'interpr{\'e}tation g{\'e}om{\'e}trique {\'e}vidente
dans le cas $LT$.

\alin{D{\'e}composition selon le centre de la cat{\'e}gorie des $\Lambda
  G$-modules lisses} 
  %Commen{\c c}ons par le cas d'un anneau 
On suppose toujours que 
$\Lambda$ est un anneau o{\`u} $p$ est inversible. Il existe une
d{\'e}composition ``par le niveau''  similaire {\`a} celle que nous avons
utilis{\'e}e pour $\dd$ ci-dessus, mais beaucoup moins triviale. Notons
$\Mo{\Lambda}{G}_n$ la sous-cat{\'e}gorie pleine de $\Mo{\Lambda}{G}$ formée des objets engendr{\'e}s
par leur vecteurs invariants sous le pro-$p$-sous-groupe ouvert
$H_n:=1+\varpi^nM_d(\OC)$ de $G=GL_d(K)$, et notons $\HC_\Lambda(G,H_n)$
l'alg{\`e}bre de Hecke de la paire $(G,H_n)$ {\`a} coefficients dans $\Lambda$.
\begin{fact} \label{hn}
  La sous-cat{\'e}gorie $\Mo{\Lambda}{G}_n$ est ``facteur direct'' de la
  cat{\'e}gorie $\Mo{\Lambda}{G}$. Les foncteurs
$$%{}{\Mo{\Lambda}{G}_n}{\HC_\Lambda(G,H_n)-\hbox{mod}}
{V}\mapsto {V^{H_n}}
\;\;\hbox{ et } \;\;
%\application{}{\HC_\Lambda(G,H_n)-\hbox{mod}}{\Mo{\Lambda}{G}_n}
{M}\mapsto{\cind{H_n}{G}{1}\otimes_{\HC_\Lambda(G,H_n)} M}
$$
sont des {\'e}quivalences ``inverses''  l'une de l'autre entre
$\Mo{\Lambda}{G}_n$ et la cat{\'e}gorie des $\HC_\Lambda(G,H_n)$-modules.
\end{fact}
\begin{proof} Nous renvoyons \`a l'appendice de \cite{finitude} (lemme A.3) o\`u une d\'ecomposition par le niveau est obtenue pour des groupes r\'eductifs plus g\'en\'eraux, en reprenant un argument de Vign\'eras reposant sur les filtrations de Moy et Prasad. Dans les notations du paragraphe A.2 de \cite{finitude}, on peut choisir pour  
``points optimaux" les isobarycentres des facettes dont l'adh\'erence contient le point sp\'ecial associ\'e \`a $GL_d(\OC_K)$. Pour un tel point $x$ et tout $r> n$, on a $G_{x,r^+}\subseteq H_n$, donc la somme $P_n$  des repr\'esentations $P(r)$, $r<n$ de {\em loc. cit.} est dans $\Mo{\Lambda}{G}_n$. R\'eciproquement, toute repr\'esentation de $\Mo{\Lambda}{G}_n$, et en particulier $\cind{H_n}{G}{\Lambda}$, contient un type non-ramifi\'e de niveau $>n$, donc est quotient d'une somme
de copies de $P_n$. D'apr\`es le lemme A.3 et le paragraphe A.1 de \cite{finitude}, la repr\'esentation $P_n$ d\'ecoupe un facteur direct de $\Mo{\Lambda}{G}$ qui n'est donc autre que $\Mo{\Lambda}{G}_n$.
\end{proof}

On associe ainsi {\`a} tout entier $n\in\NM$ un idempotent central
$\varepsilon_{n,G}$ de $\ZG_\Lambda( G)$ qui ``projette" la catégorie $\Mo{\Lambda}{G}$ sur sa sous-catégorie $\Mo{\Lambda}{G}_n$. 
Faisant agir cet idempotent sur la catégorie $\Mo{\Lambda}{GD}$, on en obtient un ``facteur direct" $\varepsilon_{n,G}\Mo{\Lambda}{GD}$. On peut interpréter ce facteur direct en termes de modules lisses sur l'algèbre de Hecke $\HC_\Lambda(GD,H_n)$ des mesures localement constantes à supports compacts et à valeurs dans $\Lambda$ qui sont invariantes à droite et à gauche sous l'image de $H_n\To{} GD$
(algèbre qui n'a pas d'unité, mais ``suffisamment d'idempotents"). En effet, en notant toujours $\cInd{H_n}{GD}$ l'induction lisse à supports compacts, on déduit du fait précédent que les foncteurs 
$$%{}{\Mo{\Lambda}{G}_n}{\HC_\Lambda(G,H_n)-\hbox{mod}}
{V}\mapsto {V^{H_n}}
\;\;\hbox{ et } \;\;
%\application{}{\HC_\Lambda(G,H_n)-\hbox{mod}}{\Mo{\Lambda}{G}_n}
{M}\mapsto{\cind{H_n}{GD}{1}\otimes_{\HC_\Lambda(GD,H_n)} M}
$$
sont inverses l'un de l'autre.
Cette fois-ci l'interpr{\'e}tation
g{\'e}om{\'e}trique se fait du c{\^o}t{\'e} $LT$ :
\begin{lemme}\label{niveauLT}
 L'inclusion de foncteurs
$\Gamma_!^{\infty_D}(\mltn^{{ca}},-) \circ \xi_{LT,n}^{*}\injo \Gamma_c(\mlt^{{ca}},-)$ induit
un isomorphisme  dans
$D^b( \HC_\Lambda(GD,H_n)\times W_K^{disc})$
$$ R\Gamma_c(\mltn^{{ca}},\Lambda) \simto
R\Gamma_c(\mlt^{{ca}},\Lambda)^{H_n}$$
qui induit {\`a} son tour un isomorphisme dans $D^b_\Lambda( GD\times W_K^{disc})$
$$ \cind{H_n}{GD}{1}\otimes_{\HC_\Lambda(GD,H_n)}R\Gamma_c(\mltn^{{ca}},\Lambda) \simto 
\varepsilon_{n,G} R\Gamma_c(\mlt^{{ca}},\Lambda).$$
\end{lemme}
\begin{proof}
Nous traitons le cas $\Lambda=\ZM_l$, les autres cas s'en d\'eduisant par extension des scalaires. Il suffit de montrer que pour tout $\Lambda_\bullet$-faisceau étale $\dd$-équivariant
 $\FC_\bullet$ sur $\PC_{LT}=S^{1/d}_K$, le morphisme canonique  
$\Gamma_{!,\bullet}^{\infty_D}(\mltn^{{ca}},\xi_{LT,n}^{*}(\FC_\bullet)) \To{} \Gamma_{c,\bullet}(\mlt^{{ca}},\FC_\bullet)$
induit un isomorphisme $\dd\times W_K$-équivariant
\ini
\begin{equation}
\label{inv}
 \Gamma_{!,\bullet}^{\infty_D}(\mltn^{{ca}},\xi_{LT,n}^{*}(\FC_\bullet)) \simto \Gamma_{c,\bullet}(\mlt^{{ca}},\FC_\bullet)^{H_n}. 
 \end{equation}
En effet, cela montrera d'abord  l'existence de l'action de $\HC_\Lambda(GD,H_n)$ sur le terme de gauche et permettra donc de définir $R\Gamma_c(\mltn^{{ca}},\Lambda)$ comme un objet $W_K$-équivariant de $D^b( \HC_\Lambda(GD,H_n))$. Puis cela montrera aussi les autres assertions du lemme, puisque le foncteur des points fixes sous $H_n$ est bien-sûr exact.

Vue la d\'efinition de $\Gamma_{c,\bullet}(\MC^{ca},-)$, pour prouver \ref{inv} il suffit de montrer que pour $m\geq n$, l'application canonique induit un isomorphisme 
$$ \Gamma_{!,\bullet}^{\infty_D}(\mltn^{{ca}},\xi_{LT,n}^{*}(\FC_\bullet)) \simto \Gamma_{!,\bullet}^{\infty_D}(\MC_{LT,m}^{{ca}},\xi_{LT,m}^{*}(\FC_\bullet))^{H_n}.$$
Rappelons que le morphisme 
$\pi_{m,n}:\;\MC_{LT,m}^{{ca}}\To{} \MC_{LT,n}^{{ca}}$ est étale galoisien de groupe $H_n/H_m$ et qu'on a
$\pi_{m,n}^* (\limproj^{\MC_n(D)} \xi_{LT,n}^{*}(\FC_\bullet))\simto \limproj^{\MC_m(D)}\xi_{LT,m}^{*}(\FC_\bullet)$.
Ainsi l'isomorphisme cherch\'e est tautologique si on enl\`eve l'indice $c$ (supports compacts), mais on peut rajouter cet indice $c$, par propret\'e et surjectivit\'e de $\pi_{m,n}$.
 \end{proof}

\begin{coro} \label{admicohoent}
  Les $\Lambda$-modules $H^q_c(\o\MC_{LT}^{{ca}},\Lambda)$, $q\in \NM$, 
  sont {\em $\o{GD}$-admissibles}, et m\^eme $G$-admissibles.
\label{finit} Par cons\'equent, le $\Lambda$-module
$\endo{D^b(\o{GD})}{\varepsilon_{n,G}R\Gamma_c(\o\MC_{LT}^{{ca}},\Lambda)}$
est de type fini. 
\end{coro}
\begin{proof}
Commen{\c c}ons par la premi\`ere assertion. Gr\^ace \`a  \ref{lien2}, il nous suffit de voir que pour tout
  $n\in\NM$, le $\Lambda$-module  $H^q_c(\mltn^{(0),{ca}},\Lambda)$ est de type fini. 
Pour cela, rappelons que $\mltn^{(0)}$ est la fibre g{\'e}n{\'e}rique d'un
  sch{\'e}ma formel $\wh\MC_{LT,n}^{(0)}$. On sait par diverses versions du th\'eor\`eme de Serre-Tate, {\em cf} \cite[II.2.7]{HaTay} et \cite[7.4.4]{Boyer1}, que ce dernier est
  isomorphe au compl{\'e}t{\'e} formel d'un sch{\'e}ma alg{\'e}brique $S_n$ propre sur $\OC^{nr}$
  (une certaine vari{\'e}t{\'e} de Shimura ou une vari{\'e}t{\'e} de
  ``Drinfeld-Stuhler'') en un point ferm{\'e} $x_n$ de sa fibre sp{\'e}ciale. Il s'ensuit que $\MC^{ca}_n$ s'identifie \`a un ouvert 
 %   Notons  $S_n^{an}$  le $\ka$-espace analytique associ{\'e} (par GAGA ou par fibre g{\'e}n{\'e}rique de la compl{\'e}tion formelle $\wh{S_n}$ le long de la fibre sp{\'e}ciale : c'est la m{\^e}me chose car $S_n$ est propre). On a un plongement
%  canonique $\mltn^{(0)} \injo S_n^{an}$ qui identifie $\mltn^{(0)}$ avec l'image
 % r{\'e}ciproque de $x_n$ par le morphisme de sp{\'e}cialisation $S_n^{an}
 % \To{} \o{S_n}$. Cette image r{\'e}ciproque est ouverte dans $S_n^{an}$
 % et son compl{\'e}mentaire s'identifie {\`a} la fibre g{\'e}n{\'e}rique du sous-sch{\'e}ma
 % formel ouvert de $\wh{S_n}$ d{\'e}fini par $\o{S_n}\setminus \{x\}$. Ce
 % compl{\'e}mentaire est donc un domaine analytique de $S^{an}_n$.  Comme
  %$S^{an}_n$ est compact, il s'ensuit que $\mltn^{(0)}$ est un ouvert 
  {\em  distingu{\'e}} (diff\'erence de deux espaces
  compacts) de l'analytifi\'e ${S^{ca,an}_n}$ de $S_n^{ca}$. 
On en d\'eduit gr\^ace \`a la suite exacte associ\'ee \`a la diff\'erence des
deux espaces compacts que sa cohomologie  est
de type fini ({\em cf} \cite[cor. 5.6]{Bic3} dans le cas de torsion et \cite[4.1.15 et  4.1.17]{Fargues} dans le cas $l$-adique).
%{\em Remarque :} bien que cela n'apparaisse pas clairement (c'est
%cach{\'e} dans les r{\'e}f{\'e}rences {\`a} \cite{Bic3} et \cite{Fargues}), tout cela repose sur la
%comparaison des cycles {\'e}vanescents de Berkovich pour les sch{\'e}mas
%formels  et les cycles {\'e}vanescents alg{\'e}briques.
%\begin{nota} \label{notaidemp} On notera
%  $\varepsilon_{n,GD}:=\varepsilon_{n,G}\varepsilon_{n,D}=
%  \varepsilon_{n,D}\varepsilon_{n,G} \in \ZG_\Lambda(GD)$ et $\Mo{\Lambda}{GD}_n$ la
%  sous-cat\'egorie facteur direct  de $\Mo{\Lambda}{GD}$ associ\'ee
%  \`a cet idempotent.  
%\end{nota}

La seconde assertion
%D'apr\`es le lemme \ref{niveauLT} et le fait rappel\'e ci-dessus, 
%il suffit de prouver le
vient du r\'esultat suivant, de nature ``th\'eorie des repr\'esentations" :

%\medskip

\noindent{\em Soit $\Lambda$ un anneau noeth\'erien o{\`u} $p$ est inversible.
  Soient $A^\bullet, B^\bullet \in D^b_\Lambda(\o{GD})$ deux complexes {\`a}
  cohomologie {$G$-admissible et de niveau fini}. 
  %({\em i.e.} dans
  %$\Mo{\Lambda}{\o{GD}}_n$ pour un $n\in\NM$).
  Alors  $\hom{A^\bullet}{B^\bullet}{D^b_\Lambda(\o{GD})}$ est un
  $\Lambda$-module de type fini.}

%\medskip

En effet, gr\^ace aux suites spectrales habituelles, il suffit de voir que chaque
$\ext{i}{\HC^q(A^\bullet)}{\HC^k(B^\bullet)}{\Lambda\o{GD}}$
 est de type fini sur $\Lambda$.
Gr{\^a}ce {\`a} la suite exacte 
$$1\To{} G/\varpi^\ZM \To{} \o{GD} \To{} \dd/K^\times \To{} 1$$
et la compacit{\'e} de $\dd/K^\times$, il suffit encore de prouver que
les $\Lambda$-modules $\ext{i}{\HC^q(A^\bullet)}{\HC^k(B^\bullet)}{\Lambda \o{G}}$ 
sont de type fini.
% Remarquons que $A$ et $B$ sont munis d'une
%action de $G$ prolongeant 
%celle de $G^0$ obtenue en plongeant $G$ dans $\o{GD}$ par
%l'application $g\mapsto (g,\varpi^{-\nu(g)})$ o{\`u} $\varpi$ est une
%uniformisante de $\dd$. De sorte que 
%$$\ext{i}{A}{B}{\Lambda G^0} \simeq \ext{i}{\cind{G^0}{G}{A}}{B}{G}.$$
Comme  le $\Lambda$-module $\hom{P}{\HC^k(B^\bullet)}{\o{G}}$ est de type fini d\`es que $P$ est projectif de type fini dans $\Mo{\Lambda}{\o{G}}$,
il suffit  de prouver que $\HC^q(A^\bullet)$ admet une r\'esolution par des objets projectifs de type fini. On peut pour cela invoquer deux arguments. 
Soit on reprend l'argument de Vign\'eras pour prouver \cite[Thm 2.11]{Vigsheaves} en rempla{\c c}ant l'ingr\'edient \cite[2.2]{Vigsheaves} par \ref{hn}. Soit on utilise la noeth\'eriannit\'e de $\Mo{\Lambda}{G}$ prouv\'ee dans \cite{finitude}, puisque $\HC^q(A^\bullet)$ est par hypoth\`ese de type fini.
\end{proof}

\alin{Une variante sans supports, du c\^ot\'e LT} Dans le paragraphe \ref{defRG}, la raison qui nous force \`a utiliser la cohomologie \`a supports compacts est notre besoin d'un complexe vivant dans la cat\'egorie d\'eriv\'ee $D_\Lambda(GD)$ des $\Lambda(GD)$-modules {\em lisses}. 
Si l'on remplace $\Gamma_!$ par $\Gamma$, on perd la lissit\'e du groupe $J$ et on s'expose \`a des probl\`emes de coefficients $l$-adiques. Du c\^ot\'e $Dr$, cela est r\'edhibitoire, mais du c\^ot\'e $LT$, le complexe obtenu nous permettra de nous raccrocher aux travaux de Boyer.

 Notons $GD^{disc}$ le groupe $GD$ muni de la topologie qui co\"incide avec la topologie naturelle de $G \subset GD$ et induit la topologie discr\`ete du quotient $GD/G=D^\times/K^\times$. Si l'on suit la proc\'edure de \ref{cohoDr} en y rempla{\c c}ant $\Gamma_!$ par $\Gamma$, on obtient un foncteur 
 $$\Gamma(\MC^{(0),ca}_{LT},-) :\; \wt{S_{K,et}^{1/d}}(D^{disc}) \To{} \wt{(GD^{disc}\times W_K^{disc})^0}$$
 puis un complexe $R\Gamma(\MC^{(0),ca}_{LT},\Lambda)\in D^b_\Lambda(GD^{disc}\times W_K^{disc})$. Si $\Lambda$ est une $\ZM_l$-alg\`ebre finie, on a 
% Dans le cas de torsion, sa cohomologie s'identifie aux modules de cohomologie $H^*(\MC^{(0),ca}_{LT},\Lambda)=\limind\,H^*(\MC^{(0),ca}_{LT,n},\Lambda)$. Dans le cas o\`u $\Lambda$ est une extension finie de $\ZM_l$, c'est encore vrai si l'on d\'efinit ceux-ci par la formule 
$$R^q\Gamma(\MC^{(0),ca}_{LT},\Lambda)=H^q(\MC^{(0),ca}_{LT},\Lambda):=\limi{n} \limp{m} H^q(\MC^{(0),ca}_{LT,n},\Lambda/l^m).$$ 
 %En effet, on a
%$$ R\Gamma(\MC^{(0),ca}_{LT},\Lambda)= R\Gamma_\bullet(\MC^{(0),ca}_{LT},\Lambda_\bullet) \simeq R\limproj \circ R\Gamma(\MC^{(0),ca}_{LT},\Lambda_\bullet) $$
%(la premi\`ere \'egalit\'e par d\'efinition et la seconde par commutation de $\Gamma$ aux limites projectives) et on sait 
puisque pour chaque $n$, le syst\`eme des $H^q(\MC^{(0),ca}_{LT,n},\Lambda/l^m)$ est  AR-$l$-adique, {\em cf} \cite[5.9]{Fargues}.

  Notons avec un exposant $^\vee$ le foncteur contragr\'e\-diente sur $\Mo{\Lambda}{GD^{disc}\times W_K^{disc}}$ et son d\'eriv\'e sur $D^b_\Lambda$.

 \begin{lemme} \label{dualite}
 Si $\Lambda=\ZM/l^n\ZM$ ou $\QM_l$, il existe un isomorphisme canonique dans $D^b_\Lambda((GD^{disc}\times W_K^{disc})^0)$ 
$$ R\Gamma(\MC^{(0),ca}_{LT},\Lambda)[2d-2](d-1) \simto R\Gamma_c(\MC^{(0),ca}_{LT},\Lambda)^\vee .$$
\end{lemme}

\begin{proof} Commen{\c c}ons par construire un accouplement 
\ini\begin{equation} \label{acc}
R\Gamma(\MC^{(0),ca}_{LT},\Lambda)\otimes^L_{\Lambda} R\Gamma_c(\MC^{(0),ca}_{LT},\Lambda)\To{} \Lambda[2-2d](1-d)
\end{equation}
dans $D^b_\Lambda((GD^{disc}\times W_K)^0)$. Pour cela, remarquons que si on se donne trois objets ab\'eliens de $\wt{\PC_{LT,et}}(D^{disc})$ munis d'un morphisme $\FC\otimes\GC \To{} \HC$, on r\'ecolte sur les sections le morphisme
$$ \Gamma(\MC_{LT}^{ca},\FC) \otimes \Gamma_c(\MC_{LT}^{ca},\GC) \To{} \Gamma_c(\MC_{LT}^{ca},\HC) $$
qui, par \'etale-finitude des $g^{m|n}$, est $GD^{disc}\times W_K$-\'equivariant. Dans le cas o\`u $\Lambda$ est de torsion, en choisissant une r\'esolution de Godement plate et multiplicative de l'anneau $\Lambda$, on obtient dans $D^b_\Lambda((GD^{disc}\times W_K^{disc})^0)$ un morphisme
$$R\Gamma(\MC^{(0),ca}_{LT},\Lambda)\otimes^L_{\Lambda} R\Gamma_c(\MC^{(0),ca}_{LT},\Lambda)
\To{} R\Gamma_c(\MC^{(0),ca}_{LT},\Lambda).$$
Lorsque $\Lambda=\ZM_l$, une r\'esolution comme-ci-dessus fournit aussi une r\'esolution flasque  de l'anneau $(\ZM/l^n)_{n\in\NM}$ dans $\Mo{\Lambda_\bullet}{\wt{\PC_{LT,et}}(D^{disc})}$ et on obtient encore un tel morphisme.
On d\'efinit ensuite \ref{acc} par composition avec la fl\`eche de $D^b_\Lambda((GD^{disc}\times W_K^{disc})^0)$ 
$$ R\Gamma_c(\MC^{(0),ca}_{LT},\Lambda) \To{\tau_{\leq_{2d-2}}} H^{2d-2}_c(\MC^{(0),{ca}}_{LT},\Lambda)[2-2d] \simto \Lambda[2-2d](1-d)$$
o\`u la premi\`ere fl\`eche est justifi\'ee par le fait que $H^q_c=0$ pour $q>2d-2$, et la seconde est donn\'ee par la famille des morphismes ``traces" $H^{2d-2}_c(\MC_{LT,n}^{ca},\Lambda)\To{\hbox{Tr}_n} \Lambda(1-d)$ normalis\'es par le facteur $|H_1/H_n|^{-1}$ (pour avoir la bonne variance). 
On en d\'eduit ensuite formellement le morphisme de l'\'enonc\'e du lemme.

Supposons maintenant que $\Lambda=\ZM/l^n\ZM$ ou $\QM_l$. Puisque cet anneau est auto-injectif, le morphisme de l'\'enonc\'e induit en cohomologie des morphismes $H^q(\MC_{LT}^{(0),ca},\Lambda)\To{} H^{2d-2-q}_c(\MC_{LT}^{(0),ca},\Lambda)^\vee(1-d)$ qui par construction, co\"incident sur les $H_n$-invariants avec les isomorphismes de dualit\'e de Poincar\'e d\^us \`a Berkovich \cite[Cor 2.3.ii)]{Bic4}
$H^q(\MC_{LT,n}^{(0),ca},\Lambda)\simto H^{2d-2-q}_c(\MC^{(0),ca}_{LT,n},\Lambda)^\vee(1-d)$ (voir aussi \cite[5.9.2]{Fargues} dans le cas $l$-adique).
\end{proof}

\section{R{\'e}alisation cohomologique des correspondances} \label{real}

Dans cette partie, nous prouvons le th\'eor\`eme \ref{main} modulo une estimation sur l'action de l'inertie qui sera obtenue dans la partie suivante.
Nous identifierons toujours les centres respectifs de $G=GL_(K)$ et $D^\times = D^\times_d$ \`a $K^\times$, via les plongements canoniques. Un ``caract\`ere central" d'une repr\'esentation  de $G$ ou $\dd$ est donc un caract\`ere de $K^\times$.

Fixons $\rho\in \Irr{\o\QM_l}{\dd}$ et notons $\omega$ son caract\`ere central.
Nous allons \'etudier dans un premier temps le complexe $$R\Gamma_c[\rho]:=R\Gamma_c(\MC^{ca},\o\QM_l)\otimes^L_{\o\QM_l\dd} \rho \,\in D^b_{\omega}(G\times W_K^{disc}).$$ %pour $\rho\in\Irr{\o\QM_l}{\dd}$ de caract\`ere central $\omega$. 
La cat\'egorie d\'eriv\'ee est celle des $\o\QM_l$-repr\'esentations de $G\times W_K$ dont la restriction \`a $G$ est  lisse et de caract\`ere central $\omega$ (rappelons que $R\Gamma_c(\MC^{ca},\o\QM_l)$ vit dans $D^b_{\o\QM_l}(GD\times W_K)$, {\em cf} \ref{rgamma}).
Par \ref{variante}, les $H^j_c(\MC^{ca},\o\QM_l)$ sont projectifs en tant que repr\'esentations lisses de $\dd$ et on a donc 
$$ \HC^j(R\Gamma_c[\rho])  \mathop{\simeq}\limits_{G\times W^{disc}_K}  H^j_c(\MC^{ca},\o\QM_l) \otimes_{\o\QM_l\dd} \rho.$$
Pour faire court, nous noterons simplement $H^j_c[\rho]$ ces objets de cohomologie.
Par \ref{admicohoent}, ils sont $G$-admissibles (et de caract\`ere central $\omega$), sont non-nuls seulement si $d-1\leq i \leq 2d-2$, et nous allons maintenant les d\'ecrire.

\subsection{Description de la cohomologie, d'apr\`es Boyer}

\alin{Repr{\'e}sentations elliptiques ``cohomologiques''}
%Soit $\rho \in \Irr{\o\QM_l}{\dd}$, 
on a d{\'e}fini et classifi{\'e} en \ref{classrepell} les
repr{\'e}sentations de $G$ elliptiques ``de type $\rho$" not\'ees $\pi_\rho^I$, pour $I\subseteq S_{\rho}$. Parmi celles-ci seules certaines apparaissent dans
la cohomologie des espaces modulaires, et m\'eritent donc d'\^etre appel\'ees  
``cohomologiques''. Pour les d{\'e}crire on utilise la bijection
$S_{\rho}\simto \I{d_\rho-1}$ d{\'e}crite en \ref{conv} et on pose 
pour $0\leq i\leq d_\rho -1$ % et $1\leq j\leq d_\rho$ :
$$ \pi^{\leq i}_\rho:=\pi^{\I{i}}_\rho \;\; \hbox{ et } \pi^{> i}_\rho:=\pi^{\{i+1,\cdots, d_\rho-1\}}_\rho \;\;
 .$$
On convient naturellement que $\pi^{\leq 0}_\rho=\pi^{\geq d_\rho}_\rho=\pi^\emptyset_\rho=JL_d(\rho)$ qui est la s{\'e}rie
discr{\`e}te de $G$ associ{\'e}e {\`a} $\rho$ par la correspondance de Jacquet-Langlands.

\begin{theo} (Boyer) \label{conjHarris}\label{isocoho}
  Pour  $0\leq i\leq d-1$, il existe un isomorphisme de $G\times \dd\times W_K$-modules 
$$ H^{d-1+i}_c[\rho] \simto 
\pi^{\leq i}_{\rho}\otimes \sigma_{d/d_\rho}(\pi_\rho^\vee)|-|^{{{d_\rho-d}\over 2}-i}.$$
Les notations $\pi_\rho$ et $d_\rho$ sont celles de \ref{notationsrho}.
\end{theo}
\begin{proof}
  Cette ``preuve" va essentiellement consister \`a expliciter un dictionnaire entre les notations de \cite{Boyer2} et celles du pr\'esent travail.
Par \ref{variante}, on a %les $H^j_c(\MC^{ca},\o\QM_l)$ sont projectifs en tant que repr\'esentations lisses de $\dd$ et on a donc  on a
\begin{eqnarray*}
 H^j_c[\rho]  % &  \mathop{\simeq}\limits_{G\times W_K} & H^j_c(\MC^{ca},\o\QM_l) \otimes_{\o\QM_l\dd} \rho\\ 
 & \mathop{\simeq}\limits_{G\times W_K} & H^j_c(\MC^{(0),ca},\o\QM_l) \otimes_{\o\QM_l\OC_D^\times} \rho
 \end{eqnarray*}
o\`u  l'action de $G\times W_K$ sur le dernier terme vient de la suite exacte $\OC_D^\times \injo (G\times\dd \times W_K)^0 \twoheadrightarrow G\times W_K$. 
Par dualit\'e et gr\^ace \`a \ref{dualite} on obtient
\begin{eqnarray*}
 H^j_c[\rho] & \mathop{\simeq}\limits_{G\times W_K} &  \hom{H^j_c(\MC^{(0),ca},\o\QM_l) \otimes_{\o\QM_l} \rho}{\o\QM_l}{\OC_D^\times}^\vee \\
 & \mathop{\simeq}\limits_{G\times W_K} &  \hom{\rho}{H^j_c(\MC^{(0),ca},\o\QM_l)^\vee}{\OC_D^\times}^\vee \\
 & \mathop{\simeq}\limits_{G\times W_K} &  \hom{\rho}{H^{2d-2-j}(\MC^{(0),ca}_{LT},\o\QM_l)}{\OC_D^\times}^\vee(1-d).
 \end{eqnarray*}
Boyer, comme Harris et Taylor, travaille sur les cycles \'evanescents du sch\'ema formel 
$\wh{\MC}_{LT}^{(0)}$. Rappelons que d'apr\`es Berkovich \cite[Cor. 2.3.ii)]{Bic4}, on a pour tout $j\in\NM$ :
$$ H^j(\MC_{LT}^{(0),ca},\o\QM_l)  \mathop{\simeq}\limits_{G\times W_K}  R^j\Psi_\eta(\wh{\MC}_{LT}^{(0)},\o\QM_l). $$
Il d\'ecoule donc de la discussion ci-dessus que pour $0\leq i\leq d-1$, on a
$$ H^{d-1+i}_c[\rho] \mathop{\simeq}\limits_{G\times W_K} \Psi^{d-1-i}_{K,l,d}(\rho)^\vee(1-d) = \UC^{d-1-i}_{K,l,d}(\rho)^\vee(1-d)$$ 
dans les notations respectives de Harris-Taylor \cite[p. 87]{HaTay} et Boyer \cite[2.1.12]{Boyer2}.
Il ne  nous reste plus qu'\`a appliquer le th\'eor\`eme 4.1.3 de  \cite{Boyer2}, en y sp\'ecialisant les ``variables" $(s,\pi)$  en $(d_\rho,\pi_\rho^\vee) $, et 
en utilisant le dictionnaire : 
\def\llr{\longleftrightarrow}
  $$ %\pi_\rho \llr \pi^\vee, \,\,\,\,
  \rho \llr JL^{-1}(St_s(\pi))^\vee,\;\; 
    \sigma_{d/d_\rho}(\pi_\rho)\llr \hbox{rec}_K(\pi)^\vee,\;\; 
 %   d_\rho \longleftrightarrow s $$
   % et en ce qui concerne les repr\'esentations elliptiques,
%    \pi^i_\rho =\pi_\rho^{\{1,\cdots,i\}} \llr 
 %   [\overrightarrow{i},\overleftarrow{s-i-1}]_{\pi_0} \;\;\hbox{ et }\;\;
 \pi_{\rho^\vee}^{> i}   %\pi_{\rho^\vee}^{\{i+1,\cdots, d_\rho-1\}} 
    \llr [\overleftarrow{i},\overrightarrow{s-i-1}]_{\pi}
$$ 
Le th\'eor\`eme 4.1.3 de  \cite{Boyer2} nous dit alors que
$$ \UC^{d-d_\rho+i}_{K,l,d}(\rho) \mathop{\simeq}\limits_{G\times W_K} \pi_{\rho^\vee}^{> i} %\pi_{\rho^\vee}^{\{i+1,\cdots, d_\rho-1 \}}
 \otimes \sigma_{d/d_\rho}(\pi_\rho)|-|^{-(d-d_\rho+2i)/2}. $$ 
En changeant $i$ en  $d_\rho-1-i$, en dualisant et en tenant compte de l'\'egalit\'e $(\pi_{\rho^\vee}^{\geq d_\rho-i})^\vee %${\pi_{\rho^\vee}^{\{d_\rho-i,\cdots, d_\rho-1\}}}^\vee 
= \pi_\rho^{\leq i}$ donn\'ee par \ref{dualiteell}, on obtient la formule de l'\'enonc\'e.
\end{proof}

\subsection{Scindages et endomorphismes de $R\Gamma_c[\rho]$}

Nous allons d\'ecrire l'alg\`ebre des endomorphismes de $R\Gamma_c[\rho]$.

\begin{coro} \label{continuite} 
\begin{enumerate}
  \item Le complexe $R\Gamma_c[\rho]$ est scindable dans $D^b_\omega(G)$.
  \item Si $\pi\in \Irr{\omega}{G}$ est telle que  $R\hom{R\Gamma_c[\rho]}{\pi}{G,\omega}\neq 0$  alors $\pi$ est elliptique de type $\rho$.
\end{enumerate}
\end{coro}

\begin{proof}
Le point ii) est une cons\'equence imm\'ediate du th\'eor\`eme de Boyer et de \ref{theoextell} i).
Le point i) peut se d\'emontrer  de deux mani\`eres diff\'erentes, soit par le calcul des Ext de \ref{theoextell} ii), soit par un argument de poids, exactement comme dans la preuve de \cite[Prop. 4.2.2]{Dat1} \`a laquelle nous renvoyons le lecteur.  
\end{proof}

Notons que $R\Gamma_c[\rho]$ {\em n'est pas} scindable dans $D^b_\omega(G\times W_K^{disc})$, et que c'est justement l\`a tout le sel de cette histoire.
Nous poserons dans la suite 
\ini\begin{equation}\label{defsigmaprime}
 \sigma'(\rho^\vee):=\sigma_{d/d_\rho}(\pi_\rho^\vee)|-|^{\frac{d_\rho-d}{2}} \;\hbox{ et } \;
 \CC_\rho := \bigoplus_{i=0}^{d_\rho} \pi^{\leq i}_{\rho}[-(d-1+i)] \in D^b_\omega(G),
  \end{equation}
nous noterons $V_{\sigma'(\rho^\vee)}$ le $\o\QM_l$-espace sur lequel agit $\sigma'(\rho^\vee)$,
et nous appellerons {\em scindage} de $R\Gamma_c[\rho]$ tout isomorphisme
\ini\begin{equation}\label{defscin}
 \alpha:\;\; R\Gamma_c[\rho] \simto \CC_\rho\otimes_{\o\QM_l} V_{\sigma'(\rho^\vee)} \;\;\hbox{ dans }\;\; D^b_\omega(G).
 \end{equation}
 qui induit les isomorphismes $\ref{isocoho}$ sur les groupes de cohomologie.

\begin{rema}\label{triv}
Lorsque $\DC$ est une cat\'egorie $R$-lin\'eaire, $R$ \'etant un anneau commutatif unitaire, et $L$ est un $R$-module libre de type fini, on a un (une classe d'isomorphisme d') endofoncteur $C\mapsto (C \otimes_R L)$ de $\DC$. On v\'erifie imm\'ediatement qu'on a un isomorphisme canonique de $R$-alg\`ebres
$$ \endo{\DC}{C\otimes L} \simto \endo{\DC}{C} \otimes_R \endo{R}{L} . $$
\end{rema}
Ainsi, tout scindage induit un isomorphisme
$$\alpha_*:\;\; \endo{D^b_\omega(G)}{R\Gamma_c[\rho]} \simto \endo{D^b_\omega(G)}{\CC_\rho} \otimes_{\o\QM_l} \endo{\o\QM_l}{V_{\sigma'(\rho^\vee)}}. $$
Par ailleurs, on a la description de l'alg\`ebre des endomorphismes de $\CC_\rho$ 
$$ \endo{D^b_\omega(G)}{\CC_\rho} = \bigoplus_{i\geq
  j} \ext{i-j}{\pi^{\leq i}_{\rho}}{\pi^{\leq j}_{\rho}}{G,\omega}, $$ 
le produit sur l'espace de droite {\'e}tant donn{\'e} par le $\cup$-produit. Notons alors  
 $\TC_{d_\rho}$  l'alg\`ebre des $\o\QM_l$-matrices triangulaires sup\'erieures de taille $d_\rho$, et $(E_{ij})_{i\leq j}$ sa base ``canonique" :
 \begin{fact} \label{choixgen}
Pour tout choix de g\'en\'erateurs $\beta_{i,i+1} \in \ext{1}{\pi^{\leq i+1}_\rho}{\pi^{\leq i}_\rho}{G,\omega}$, l'application $E_{i,i+1}\mapsto \beta_{i,i+1}$  induit un  isomorphisme 
%$$ \beta :\;\; \TC_{d_\rho} \simto \bigoplus \ext{i-j}{\pi^{\leq i}_{\rho^\vee}}{\pi^j_{\rho^\vee}}{\o{G}}.$$
$$ \beta :\;\; \TC_{d_\rho} \simto \endo{D^b_\omega(G)}{\CC_\rho}.$$ 
Changer ce choix de g\'en\'erateurs revient simplement \`a composer $\beta$ avec la conjugaison par une matrice diagonale.
\end{fact}
\begin{proof}
Ceci d\'ecoule de la description du cup-produit \ref{theoextell} ii). On renvoie \`a \cite[4.2.4]{Dat1} pour plus de d\'etails.
\end{proof}

Par la remarque \ref{triv}, un isomorphisme $\beta$ comme ci-dessus induit un isomorphisme d'alg\`ebres encore not\'e $\beta$ 
$$ \beta :\;\; \endo{\o\QM_l}{V_{\sigma'(\rho^\vee)}}\otimes_{\o\QM_l} \TC_{d_\rho} \simto \endo{D^b_\omega(G)}{\CC_\rho\otimes_{\o\QM_l} V_{\sigma'(\rho^\vee)}}.$$
Pour r\'ecapituler, on a obtenu
\begin{prop} \label{propnoneq} 
\`A tout scindage $\alpha$ comme en \ref{defscin} et tout choix de g\'en\'erateurs 
%des droites $\ext{1}{\pi^{i}_{\rho^\vee}}{\pi^{i-1}_{\rho^\vee}}{\o{G}}$ 
\ref{choixgen} est associ\'e un
 isomorphisme de $\o\QM_l$-alg{\`e}bres $\beta^{-1}\circ \alpha_*$ s'inscrivant dans les diagrammes commutatifs 
$$\xymatrix{ \endo{D^b_\omega(G)}{R\Gamma_c[\rho]} \ar[r]^\sim_{\beta^{-1}\circ \alpha_*} \ar[d]^{can}
& \endo{\o\QM_l}{V_{\sigma'(\rho^\vee)}} \otimes_{\o\QM_l} \TC_{d_\rho}
\ar[d]^{\id\otimes E_{ii}^*} \\
 \endo{G}{H^{d-1+i}_c[\rho]} \ar[r]^{(\ref{isocoho})}  &
 \endo{\o\QM_l}{V_{\sigma'(\rho^\vee)}} },$$
o\`u  $i\in\{0,\cdots, d_\rho\}$ et $E^*_{ii}$ d\'esigne la coordonn\'ee diagonale $(i,i)$ de $\TC_{d_\rho}$
 \end{prop}

\subsection{Action de $W_K$ sur $R\Gamma_c[\rho]$}
 On \'etudie ici  le morphisme canonique
$$ \gamma_\rho : \;\; W_K \To{}\aut{D^b_\omega(G)}{R\Gamma_c[\rho]} $$
qui d{\'e}crit l'action de $W_K$ sur $R\Gamma_c[\rho]$.
Le r{\'e}sultat final que nous visons est une version $W_K$-{\'e}quivariante
de la proposition \ref{propnoneq} :

\begin{theo} \label{propeq}
Pour tout rel\`evement de Frobenius g\'eom\'etrique $\phi$, il existe un unique scindage $\alpha_\phi$ comme en \ref{defscin} et
% un unique sous-ensemble $I_\phi\subseteq S_{d_\rho}=\{1,\cdots,d_\rho-1\}$, 
un %isomorphisme $\beta_\phi$ associ\'e \`a un certain
unique choix de g\'en\'erateurs \ref{choixgen} induisant un isomorphime $\beta_\phi$,
%des droites $\ext{1}{\pi^{\leq i}_{\rho^\vee}}{\pi^{i-1}_{\rho^\vee}}{\o{G}}$,
  tels que le diagramme suivant commute :
$$\xymatrix{\endo{D^b_\omega(G)}{R\Gamma_c[\rho]} \ar[r]^\sim_{ \beta_\phi^{-1}\circ \alpha_{\phi*}} 
& \endo{\o\QM_l}{V_{\sigma'(\rho^\vee)}} \otimes_{\o\QM_l} \TC_{d_\rho} \\
W_K \ar[u]^{\gamma_\rho}  \ar[ur]_{ %\sigma_d(JL(\rho^\vee))}
\sigma'(\rho^\vee)\otimes \tau_{d_\rho}^{\emptyset,\phi}}
 & }.$$
De plus, le choix de g\'en\'erateurs, et donc $\beta_\phi$, est ind\'ependant de $\phi$.
\end{theo}

%\begin{rema} 
Rappelons que la notation $\tau_{d_\rho}^{\emptyset,\phi}$ a \'et\'e introduite au-dessus du lemme \ref{corell}. Sa classe d'isomorphisme est, \`a torsion pr\`es par $|-|^{\frac{d_\rho-1}{2}}$, la repr\'esentation sp\'eciale de dimension $d_\rho$. En cons\'equence, la classe d'isomorphisme de la repr\'esentation $\sigma'(\rho^\vee)\otimes \tau_{d_\rho}^{\emptyset,\phi}$ est la correspondante de Langlands $\sigma_d(JL_d(\rho^\vee))$ tordue par $|-|^{\frac{1-d}{2}}$.

Par ailleurs, nous avons fait le choix contestable de fixer un g\'en\'erateur $\mu$ de $\ZM_l(1)$ dans les constructions qui interviennent pour l'\'enonc\'e ci-dessus. On laissera au lecteur int\'eress\'e le soin de suivre l'effet d'un tel choix dans ces constructions et on renvoie \`a la partie 4 de \cite{Dat1}, pour une discussion plus ``canonique" d'une situation analogue.
%\end{rema}

Le reste du  paragraphe est consacr{\'e} {\`a} la preuve du th\'eor\`eme ci-dessus, en admettant la proposition \ref{estimN} ci-dessous, laquelle sera prouv\'ee dans la section suivante.

\begin{lemme} (monodromie quasi-unipotente) \label{monqunip}
  Il existe un unique endomorphisme nilpotent $$N_\rho\in
  \endo{D^b_\omega(G)}{R\Gamma_c[\rho]}$$ tel que, si 
  $I_\rho\subseteq I_K$ d\'esigne le noyau de $\sigma'(\rho^\vee)_{|I_K}$, alors
 $$\forall i\in I_\rho,\;\; \gamma_\rho(i) = \hbox{exp}(N_\rho t_\mu(i)) .$$
 De plus on a 
 \ini
\begin{equation}
  \label{eqNrho}
\forall w\in W_K,\;\;  \gamma_\rho(w)N_\rho\gamma_\rho(w)^{-1} = q^{-\nu(w)} N_\rho=|w|N_\rho.
\end{equation}
\end{lemme}
\begin{proof} 
Commen{\c c}ons par remarquer que l'\'enonc\'e est insensible \`a la torsion par les caract\`eres non-ramifi\'es. En effet, si $\psi : \ZM\To{} \o\QM_l$ est un caract\`ere, alors d'apr\`es \ref{variante}, on a un isomorphisme 
$$ R\Gamma_c[\rho\otimes(\psi\circ |\hbox{Nr}|_{_K})] \simeq R\Gamma_c[\rho] \otimes (\psi\circ |\hbox{det}|_{_K})(\psi\circ |-|^{-1}_{_{W}}) \;\;\hbox{ dans } \; D^b_{\omega(\psi\circ |.|^d_{_K})}(G\times W_K^{disc}).$$
Quitte \`a tordre $\rho$ par un caract\`ere non-ramifi\'e, on peut donc supposer  $\omega(\varpi)=1$. On a alors 
 $$R\Gamma_c[\rho]=R\Gamma_c(\o\MC^{ca},\o\QM_l)\otimes^L_{\o\QM_l\o{D}} \rho, $$
en posant $\o{D}:=\dd/\varpi^\ZM$.
Choisissons de plus  $n\in\NM $ tel que $JL(\rho)^{H_n}\neq 0$ (rappelons que $H_n=1+\varpi^nM_d(\OC_K)$), alors d'apr\`es la description de la cohomologie \ref{conjHarris}, et en utilisant les notations de \ref{hn}, on a m\^eme
$$R\Gamma_c[\rho] = (\varepsilon_{n,G} R\Gamma_c(\o\MC^{ca},\o\QM_l))\otimes^L_{\o\QM_l\o{D}} \rho.$$
Par d\'efinition, le morphisme $\gamma_\rho$ donnant l'action de  $W_K$ se factorise par 
\ini\begin{equation}
\label{fac}
 \endo{D^b_{\ZM_l}(\o{GD})}{\varepsilon_{n,G}R\Gamma_c(\o\MC^{{ca}},\ZM_l)}
 \To{\otimes_{\ZM_l\o{D}} \rho}  \endo{D^b_{\omega}
   (G)}{R\Gamma_c[\rho]}.
\end{equation}
On sait par \ref{admicohoent} que le $\ZM_l$-module de gauche est de type fini.

Notons maintenant $\NC(\rho)$ le noyau du morphisme canonique
$$  \endo{D^b_\omega(G)}{R\Gamma_c[\rho]} \To{can} \prod_{i}
\endo{G}{H^{d-1+i}_c[\rho]} $$
et 
 $\NC(\rho)^0$ son intersection  avec l'image de \ref{fac}, qui est donc un $\ZM_l$-module de type fini. On sait que $\NC(\rho)$ est form{\'e}
d'{\'e}l{\'e}ments nilpotents d'ordre $\leq d_\rho$, {\em cf} \cite[A.1.4]{Dat1} par exemple, et que le groupe $\UC^0(\rho):=1+\NC^0(\rho)$ est naturellement un pro-$l$-groupe, {\em cf} \cite[A.2.1]{Dat1}. Plus pr\'ecis\'ement, on a
$$ \UC(\rho)^0 \simto \limp{n} \UC(\rho)^0/(\id+l^n\NC(\rho)^0). $$

Par d\'efinition, le groupe $I_\rho$ de l'\'enonc\'e est un sous-groupe d'indice fini de
$I_K$ qui est aussi dans le noyau de l'action de $W_K$ sur chaque $H^i_c[\rho]$, par la description de la cohomologie  \ref{conjHarris}. Par cons\'equent
$\gamma_\rho$ envoie $I_\rho$ dans le sous-groupe $\UC(\rho)^0$ de 
$\aut{}{R\Gamma_c[\rho]}$.
% et m\^eme dans le sous-groupe $\UC(\rho)^0:=1+\NC(\rho)^0$ puisque l'action de $W_K$ sur $R\Gamma_c[\rho]$ est induite par celle sur $R\Gamma_c(\o\MC,\ZM_l)$. 
%Mais puisque $\NC(\rho)^0$ est $l$-adiquement complet, l'application canonique
%$$ \UC(\rho)^0 \To{} \limp{n} \UC(\rho)^0/(\id+l^n\NC(\rho)^0) $$
%est un isomorphisme de groupe qui munit $\UC(\rho)^0$ d'une structure de pro-$l$-groupe.
Comme un morphisme de groupes entre groupes profinis est automatiquement continu, le morphisme $I_\rho \To{\gamma_\rho} \UC(\rho)^0$ se factorise par le plus grand pro-$l$-quotient de $I_\rho$.
Mais celui-ci n'est autre que l'image de $I_\rho$ par
le morphisme  $t_\mu : I_K \To{} \ZM_l$. 
Cette image est d'indice fini, donc de la forme $l^m\ZM_l$, et l'on peut par cons\'equent \'ecrire  
$$ \forall i\in I_\rho,\;\; \gamma_\rho(i)= u^{t_\mu(i)/l^m} $$
o\`u $u \in \UC(\rho)^0$ est l'image par $\gamma_\rho$ de n'importe quel \'el\'ement de $I_\rho$ s'envoyant sur $l^m$ par $t_\mu$.

Posant alors $N_\rho:= l^{-m} \hbox{log}(u) \in \NC(\rho)$, (le logarithme \'etant bien d\'efini puisque $u$ est unipotent),
%En composant avec le logarithme, on obtient
%un morphisme de groupes $\hbox{log}\circ \gamma_\rho:\;\; I_\rho \To{}
%\NC(\rho)$ \`a valeurs dans le $\o\QM_l$-espace vectoriel $\NC(\rho)$. On voudrait montrer que %ce morphisme est de la forme
%\ini\begin{equation} \label{forme}
on obtient
$$  \gamma_\rho (i) = \hbox{exp}(t_\mu(i).N_\rho) ,\;\;\forall i\in I_\rho. $$
%\end{equation}
%pour un certain $N_\rho\in \NC(\rho)$. 
%D'apr{\`e}s la proposition \ref{continuite}, l'hypoth{\`e}se
%d'admissibilit{\'e} de la cohomologie enti{\`e}re implique  que 
%{\em l'image du morphisme  $
%  \gamma_\rho$ (ou de mani{\`e}re {\'e}quivalente, celle de $\hbox{log}\circ \gamma_\rho$)
%  est contenue dans un sous $\ZM_l$-module de type fini de
%  $\NC(\rho)$}. 
%Soit $Im$ un tel $\ZM_l$-module que l'on munit de sa topologie $l$-profinie.
%Comme $I_\rho$ est profini, le morphisme de groupes $I_\rho \To{\hbox{log}\circ \gamma_\rho} Im $ est automatiquement continu et se factorise donc par le plus grand $l$-quotient de $I_\rho$, donc par $t_\mu:I_\rho\To{} \ZM_l$. 
%On peut ainsi le mettre sous la forme \ref{forme}.
Puisque pour tous $(w,i)\in W_K\times I_\rho$, on a $t_\mu(wiw^{-1})=q^{-\nu(w)}t_\mu(i)$, l'endomorphisme $N_\rho$ doit satisfaire
l'{\'e}quation \ref{eqNrho}.
Enfin, son unicit\'e est \'evidente.

%Pour l'unicit\'e, soit $N'$ un autre endomorphisme nilpotent, et supposons que le %${\wt\gamma_\rho}'$ associ\'e soit trivial sur le sous-groupe ouvert $I'$ de $I_\rho$. Alors,  pour tout $i\in I'$, on a $\gamma_\rho(i)=\hbox{exp}(N_\rho t_\mu(i))=\hbox{exp}(N't_\mu(i))$ donc $N'=N_\rho$. 

%Quitte {\`a} tordre $\rho$ par un caract{\`e}re non ramifi{\'e} de $D$ on peut la
%supposer  son caract{\`e}re central d'ordre fini. Il existe alors un
%sous-corps $E_\rho$ de $\o\QM_l$ fini sur $\QM_l$ tel que $\rho$ soit
%d{\'e}finie sur $E_\rho$. L'idempotent primitif $[\rho]$ de $\ZG(\o\QM_lD)$ associ{\'e}
%{\`a} $\rho$ appartient alors {\`a} $\ZG(E_\rho D)$ et agit par cons{\'e}quent sur
%$D^b(E_\rhoGD)$. Il est alors clair 
%$$ R\Gamma_c[\rho] = R\Gamma_c(\breve\Sigma_n,\o\QM_l)[\rho] \simeq \o\QM_l \otimes_{E_\rho}
%R\Gamma_c(\breve\Sigma_n, E_\rho)[\rho], $$
%l'isomorphisme (canonique) {\'e}tant compatible {\`a} l'action de $W_K$. 

%Comme le groupe $\aut{D^b(\o\QM_lGD)}{R\Gamma_c[\rho]}$ est le
%groupe des automorphismes d'une $\o\QM_l$-alg{\`e}bre 

\end{proof}

Soit $\phi$ un rel\`evement de Frobenius g\'eom\'etrique. Consid\'erons 
l'application
\ini
\begin{equation}\label{gammaphi} 
\application{\gamma^\phi_\rho(w):\;\;}{W_K}{\endo{D^b_\omega(G)}
  {R\Gamma_c[\rho]}^\times}{w}{\gamma_\rho(w) \hbox{exp}(-N_\rho t_\mu(i_\phi(w))) \hbox{ o{\`u} } w=\phi^{\nu(w)}i_\phi(w)}.
\end{equation}
L'{\'e}quation \ref{eqNrho} assure que l'application $\gamma^\phi_\rho$  est  un morphisme de groupes $W_K \To{} 
\aut{}{R\Gamma_c[\rho]}$. Par construction, celui-ci se factorise par
$W_K\To{} W_K/I_\rho$ qui est discret.

\begin{lemme} \label{lemmeeq}
%(M\^emes hypoth\`eses que la prop. \ref{propeq}).
Il existe un unique scindage $\alpha_\phi$ comme en \ref{defscin}  tel que le diagramme suivant commute :
$$\xymatrix{\endo{D^b_\omega(G)}{R\Gamma_c[\rho]} \ar[r]^\sim_{\beta^{-1}\circ \alpha_{\phi*}} 
& \endo{\o\QM_l}{V_{\sigma'(\rho^\vee)}} \otimes_{\o\QM_l} \TC_{d_\rho} \\
W_K \ar[u]^{\gamma^\phi_\rho}  \ar[ur]_{\sigma'(\rho^\vee)\otimes \tau^{ss}_{{d_\rho}}} & }.$$
Rappelons que la repr\'esentation $\tau^{ss}_{{d_\rho}}$ se factorise par la diagonale de $\TC_{d_\rho}$, donc  la commutation du diagramme est ind\'ependante du choix de g\'en\'erateurs de \ref{choixgen} pour d\'efinir $\beta$.
\end{lemme}

\begin{proof}
Comme dans le lemme \ref{monqunip}, notons $I_\rho\subset I_K$ le noyau de $\sigma'(\rho^\vee)_{|I_K}$, qui est aussi celui de $(\gamma^\phi_\rho)_{|I_K}$, par construction. 
%Soit $\phi$ un rel{\`e}vement de Frobenius g\'eom\'etrique dans $W_K$. 
Puisque $I_K/I_\rho$ est fini, on peut trouver un entier
$n\in\NM$ tel que (l'image de) $\phi^n$ soit {\em central} dans le groupe
$W_K/I_\rho$.
Par le lemme de Schur, l'automorphisme $\sigma'(\rho^\vee)(\phi^n)$ de $V_{\sigma'(\rho^\vee)}$ est un scalaire que l'on notera $\xi\in\o\QM_l^\times$.
Par la description \ref{isocoho} de la cohomologie, on voit que l'endomorphisme $\HC^{d-1+j}(\gamma^\phi_\rho(\phi^n))$ de $H^{d-1+j}_c[\rho]$ ({\em i.e} l'action de $\phi^n$ sur $H^{d-1+j}_c[\rho]$) est annul\'e par le polyn\^ome $X-q^{nj}\xi$. Par  \cite[lemme A.1.4]{Dat1}, il existe donc un unique scindage $\alpha_\phi$ tel que (pour tout choix de $\beta$)  
$$ \beta^{-1}\circ \alpha_{\phi*}(\gamma^\phi_\rho(\phi^n)) = \xi \otimes \hbox{Diag}(1,q^n,\cdots,q^{nd_\rho}) \in \endo{\o\QM_l}{V_{\sigma'(\rho^\vee)}}\otimes \TC_{d_\rho}.$$
Or, le commutant de $\xi \otimes \hbox{Diag}(1,q^n,\cdots,q^{nd_\rho})$ dans  $\endo{\o\QM_l}{V_{\sigma'(\rho^\vee)}}\otimes \TC_{d_\rho}$ est $\endo{\o\QM_l}{V_{\sigma'(\rho^\vee)}}\otimes \o\QM_l^{d_\rho}$ o\`u $\o\QM_l^{d_\rho}$ d\'esigne ici la sous-alg\`ebre des matrices diagonales de $\TC_{d_\rho}$. Ainsi, puisque $\phi^n$ est central dans $W_K/\ker(\gamma^\phi_\rho)$, on a 
$$ \hbox{Im} (\beta^{-1}\circ \alpha_{\phi*}\circ \gamma^\phi_\rho) \subset \endo{\o\QM_l}{V_{\sigma'(\rho^\vee)}}\otimes \o\QM_l^{d_\rho}. $$
Mais par le diagramme commutatif de la proposition \ref{propnoneq} et la $W_K$-\'equivariance des isomorphismes \ref{isocoho}, on en d\'eduit que 
\begin{eqnarray*}
 \forall w\in W_K,\;\; (\beta^{-1}\circ \alpha_{\phi*})(\gamma^\phi_\rho(w)) & = & \sum_{i=0}^{d_\rho} \sigma'(\rho^\vee)(w)|w|^{-i} \otimes E_{ii} \\
 & = & \sum_{i=0}^{d_\rho} \sigma'(\rho^\vee)(w) \otimes |w|^{-i}E_{ii}
\end{eqnarray*}
o\`u $E_{ii}$ d\'esigne toujours la matrice \'el\'ementaire de coordonn\'ees $(i,i)$ de $\TC_{d_\rho}$. Le scindage $\alpha_\phi$ fait donc bien commuter le diagramme du lemme. Pour l'unicit\'e, remarquons que toute autre solution $\alpha'$ v\'erifie la conclusion du point ii) de \cite[lemme A.1.4]{Dat1} pour $\phi^n$, et co\"incide donc avec $\alpha$ par l'assertion d'unicit\'e de ce m\^eme point.

\end{proof}

\alin{Preuve du th\'eor\`eme \ref{propeq}}
On garde les notations pr\'ec\'edentes ; en particulier, $\alpha_\phi$ est le scindage de $R\Gamma_c[\rho]$ du lemme \ref{lemmeeq}, et
$N_\rho$ 
l'endomorphisme nilpotent de $R\Gamma_c[\rho]$ donn{\'e} par le lemme
\ref{monqunip}. 
%On voudrait montrer que pour un bon $\beta=\beta_\phi$ (correspondant \`a un choix de g\'en\'erateurs des $\ext{i}{\pi^i_{\rho^\vee}}{\pi^{i-1}_{\rho^\vee}}{\o{G}}$), il existe un certain sous-ensemble $I\subseteq S_{d_\rho}$ tel que
%$$(\beta_\phi^{-1}\circ\alpha_{\phi*})(N_\rho) = (\id_{V_{\sigma'(\rho^\vee)}}\otimes N_{I}) \in
%\endo{\o\QM_l}{V_{\sigma'(\rho^\vee)}} \otimes_{\o\QM_l} 
%\TC_{d_\rho}$$
%o\`u $N_I$ a \'et\'e d\'efini en \ref{deftauI}.
Choisissons un $\beta$ arbitraire pour commencer.
Comme $N_\rho$ induit l'endomorphisme nul en cohomologie, on a
$$(\beta^{-1}\circ\alpha_{\phi*})(N_\rho) \in \sum_{i>j} \endo{\o\QM_l}{V_{\sigma'(\rho^\vee)}}
\otimes E_{ij}, $$
en notant $(E_{ij})_{i\geq j}$ la base ``canonique" de $\TC_{d_\rho}$. 
D'autre part, on d\'eduit de \ref{eqNrho} la relation 
\ini\begin{equation}\label{rel}
 \gamma^\phi_\rho(w) N_\rho \gamma^\phi_\rho(w)^{-1} = |w| N_\rho 
\end{equation}
pour tout $w\in W_K$. Appliqu{\'e}e {\`a} l'{\'e}l{\'e}ment $\phi^n$ de la preuve du
lemme \ref{monqunip}, pour lequel on a $$(\beta^{-1}\circ\alpha_{\phi*})(\phi^n)= \sum_{i=0}^{d_\rho-1} \xi q^{ni}(\id \otimes E_{ii}),$$
 cette relation implique que 
$$(\beta^{-1}\circ\alpha_{\phi*})(N_\rho) \in \sum_{i=1}^{d_\rho-1} \endo{\o\QM_l}{V_{\sigma'(\rho^\vee)}}
\otimes E_{i-1,i}.$$ 
  {\'E}crivons donc $(\beta^{-1}\circ\alpha_{\phi*})(N_\rho)= \sum_i M_i\otimes E_{i-1,i}$ avec $M_i\in
\endo{\o\QM_l}{V_{\sigma'(\rho^\vee)}}$. 
Par d{\'e}finition de $\alpha_\phi$, on a pour tout $w\in W_K$ 
$$ (\beta^{-1}\circ\alpha_{\phi*})(\gamma^\phi_\rho(w)) = \sum_{i=0}^{d_\rho-1} 
\sigma'(\rho^\vee)(w) \otimes |w|^{-i}E_{ii}. $$
La relation \ref{rel} implique donc que chaque $M_i$ commute avec
$\sigma'(\rho^\vee)(w)$ pour tout $w$. Par le lemme de Schur, on a donc
$M_i\in\o\QM_l$. Quitte \`a changer $\beta$, on peut alors supposer que $M_i$ est soit nul, soit \'egal \`a $1$. Pour achever la preuve du th\'eor\`eme \ref{propeq}, il nous reste \`a prouver que les $M_i$ ainsi obtenus sont tous inversibles, auquel cas,
le changement de $\beta$ est d'ailleurs unique. 
De mani\`ere \'equivalente, il  nous reste \`a prouver :
\begin{prop} \label{estimN}
$N_\rho^{d_\rho-1} \neq 0$.
\end{prop}
Nous reportons la preuve \`a la section \ref{monodromie}.

\subsection{Preuve du th{\'e}or{\`e}me \ref{main}}
Nous prouvons ici le th\'eor\`eme \ref{main}, \`a partir du th\'eor\`eme \ref{propeq}.
Nous noterons simplement $R\Gamma_c:=R\Gamma_c(\MC^{ca},\o\QM_l)$ l'objet de $D^b_{\o\QM_l}(GD\times W_K^{disc})$ construit en \ref{rgamma}. 

\begin{lemme}  \label{reduc1}
Les deux \'enonc\'es ci-dessous sont \'equivalents \`a l'\'enonc\'e  du th\'eor\`eme \ref{main} :
\begin{enumerate}
        \item Pour toutes repr\'esentations $\pi\in\Irr{\o\QM_l}{G}$, $\rho\in\Irr{\o\QM_l}{\dd}$ de m\^eme caract\`ere central, on a :
        $$\HC^*\left(R\hom{R\Gamma_c}{\pi\otimes
    \rho^\vee}{D^b(G D)}\right)
\mathop{\simeq}\limits_{W_K} \left\{\begin{array}{rl} 
  \sigma_d(\pi)|-|^{\frac{d-1}{2}} & \hbox{ si }  \rho=LJ_d(\pi) \\
                             0 & \hbox{ sinon } \end{array}\right.. $$
        \item Pour toutes repr\'esentations $\pi\in\Irr{\o\QM_l}{G}$, $\rho\in\Irr{\o\QM_l}{\dd}$ de m\^eme caract\`ere central $\omega$, on a :
        $$\HC^*\left(R\hom{R\Gamma_c[\rho]}{\pi}{D^b_\omega(G)}\right)
\mathop{\simeq}\limits_{W_K} \left\{\begin{array}{rl} 
  \sigma_d(\pi)|-|^{\frac{d-1}{2}} & \hbox{ si }  \rho=LJ_d(\pi) \\
                             0 & \hbox{ sinon } \end{array}\right.. $$
                              
\end{enumerate}
La d\'efinition de $\HC^*$ a \'et\'e donn\'ee dans l'introduction. %La notation $D^b_\omega(G)$ d\'esigne la cat\'egorie d\'eriv\'ee des $\o\QM_lG$-repr\'esentations lisses de caract\`ere central $\omega$.
\end{lemme}

\begin{proof}
Soit $\omega$ le caract\`ere central commun. On a une factorisation
$$ \hom{-}{\pi\otimes\rho^\vee}{GD}:\;\;\Mo{\o\QM_l}{GD}
\To{\hom{-}{\pi}{G}} \Mo{\omega}{\dd} 
\To{\hom{\rho}{-}{\dd}} \o\QM_l-\hbox{e.v.} $$
o\`u les deux derni\`eres cat\'egories sont semi-simples.
En particulier, on a un isomorphisme
$$  \hom{\rho}{\HC^*(R\hom{R\Gamma_c}{\pi}{G})}{\dd} \simto 
\HC^*(R\hom{R\Gamma_c}{\pi\otimes\rho^\vee}{GD}) $$
qui bien-s\^ur est $W_K$-\'equivariant.
Il est alors clair que l'\'enonc\'e du th\'eor\`eme A implique le point i).
 Compte tenu du fait que la famille de
foncteurs $\hom{\rho}{-}{\dd}$, $\rho \in \Irr{\omega}{\dd}$
est fid\`ele sur la cat\'egorie $\Mo{\omega}{\dd}$, la
r\'eciproque est encore vraie. Le  point ii) est \'equivalent au point i), gr\^ace \`a la factorisation
$$ \hom{-}{\pi\otimes\rho^\vee}{GD}:\;\;\Mo{\o\QM_l}{GD}
\To{-\otimes_{\QM_l\dd} \rho} \Mo{\omega}{G} 
\To{\hom{-}{\pi}{G}} \o\QM_l-\hbox{e.v.}. $$
\end{proof}

Fixons \`a nouveau $\rho\in \Irr{\o\QM_l}{\dd}$ de caract\`ere $\omega$, ainsi qu'un sous-ensemble $I\subseteq S_{\rho}$, et posons
$$ \HC_\rho^I:= \HC^*\left( R\hom{R\Gamma_c[\rho]}{\pi_\rho^I}{D^b_\omega(G)}   \right). $$
C'est un $\o\QM_l$-espace vectoriel gradu\'e de dimension finie muni
de l'action par automorphismes de degr\'e $0$ de $W_K$ induite par
l'action $\gamma_{\rho} :\; W_K\To{}
\aut{D^b_\omega(G)}{R\Gamma_c[{\rho}]}$. En vertu de \ref{reduc1} ii) et \ref{continuite} ii), le th\'eor\`eme \ref{main} sera prouv\'e une fois qu'on aura identifi\'e $\HC_\rho^I$ \`a $\sigma_d(\pi^I_\rho)|-|^{\frac{d-1}{2}}$.

Notons $N^I_\rho$ l'endomorphisme de degr\'e $0$ (d'espace vectoriel
gradu\'e) de $\HC_\rho^I$ induit fonctoriellement par l'endomorphisme
$N_{\rho}$ de $R\Gamma_c[{\rho}]$ du lemme \ref{monqunip}.
Choisissons un rel\`evement de Frobenius g\'eom\'etrique $\phi$ 
et notons $\HC_\rho^{I,\phi}$ la repr\'esentation gradu\'ee lisse de $W_K$
induite par l'action $\gamma^\phi_{\rho}:\; W_K\To{}
\aut{D^b_\omega(G)}{R\Gamma_c[{\rho}]}$ d\'efinie en \ref{gammaphi}. 
% L'endomorphisme $N_\rho^I$ est
%l'op\'erateur de monodromie de la $W_K$-repr\'esentation gradu\'ee $\HC^I_\rho$
%et la connaissance de cette derni\`ere \'equivaut \`a celle du 
Le couple $(\HC_\rho^{I,\phi},N_\rho^I)$ est la repr\'esentation de Weil-Deligne associ\'ee \`a $\HC_\rho^I$ (et aux choix de $\phi$ et $\mu$).

En utilisant le scindage $\alpha_\phi$ de \ref{lemmeeq}, on obtient
 la premi\`ere ligne de
\begin{eqnarray*}
 \HC_\rho^{I,\phi} & \simeq_{W_K} & \bigoplus_{i=0}^{d_\rho-1} \HC^*\left(R\hom{\pi_\rho^{\leq i}\otimes \sigma'({\rho^\vee})|-|^{-i}}{\pi_\rho^I }{D^b_\omega(G)} \right)[d-1+i] \\
% & \simeq_{W_K} & \bigoplus_{i=0}^{d_\rho-1}
% \ext{\delta(i,I)}{\pi_\rho^i}{\pi_\rho^I}{G,\omega} \otimes_{\o\QM_l}
% \left(\sigma'(\rho^\vee)|-|^{-i}\right)^\vee [d-1+i-\delta(i,I)] \\
 & \simeq_{W_K} & \sigma'(\rho^\vee)^\vee \otimes\left(\bigoplus_{i=0}^{d_\rho-1}
 \ext{\delta(i,I)}{\pi_\rho^{\leq i}}{\pi_\rho^I}{{G,\omega}} \otimes_{\o\QM_l}
 |-|^{i} [d-1+i-\delta(i,I)]\right)
 \end{eqnarray*}
 et la deuxi\`eme ligne vient  du calcul d'extensions de \ref{theoextell}.

On veut maintenant expliciter l'endomorphisme $N_\rho^I$. On proc\`ede
exactement comme dans \cite[4.4]{Dat1}.
Notons $\beta_{i-1,i} \in \ext{1}{\pi_\rho^{\leq i}
  }{\pi_\rho^{\leq i-1}}{{G,\omega}}$, o\`u $i\in
\{1,\cdots,d_\rho-1\}$ les g\'en\'erateurs qui d\'efinissent l'isomorphisme $\beta_\phi$ du th\'eor\`eme  \ref{propeq}  et choisissons 
pour chaque $j\in \{0,\cdots,d_\rho-1\}$ un
g\'en\'erateur de la droite
$\ext{\delta(i,I)}{\pi^{\leq i}_\rho}{\pi^I_\rho}{{G,\omega}}$.
Alors l'action (\`a droite) de $\alpha_{\phi*}(N_{\rho})$ sur $\HC_\rho^I$ est
donn\'ee par 
$$ e_i \mapsto e_i \cup \beta_{i,i+1} \in \o\QM_l e_{i+1}. $$
La formule pour le $\cup$-produit de \ref{theoextell} ii) montre que
$e_i\cup \beta_{i,i+1}$ est non-nul \ssi\
$\delta(i+1,I)=\delta(i,I)+1$ et des consid\'erations
\'el\'ementaires montrent que cela se produit \ssi\ $i+1\in I^c$
(compl\'ementaire de $I$ dans $S_{\rho}=\{1,\cdots,d_\rho-1\}$).

Oublions maintenant la graduation de $\HC^I_\rho$. La
discussion de \cite[4.4]{Dat1} montre alors que 
$$ \HC^I_\rho \simeq_{W_K} \sigma'(\rho^\vee)^\vee \otimes \tau_{d_\rho}^{I}, $$
et compte tenu de la d\'efinition \ref{defsigmaprime} de $\sigma'(\rho^\vee)$, de l'\'egalit\'e $\pi_{\rho^\vee}=\pi_\rho^\vee$ et de la compatibilit\'e de la correspondance de Langlands aux contragr\'edientes, on obtient
\begin{eqnarray*}
 \HC^I_\rho & \simeq_{W_K}  & \sigma_{d/d_\rho}(\pi_\rho)|-|^{\frac{d-d_\rho}{2}}\otimes \tau_{d_\rho}^{I} \\
& \simeq_{W_K} & \sigma_{d/d_\rho}(\pi_\rho)\otimes \tau_{d_\rho}^{I}|-|^{\frac{1-d_\rho}{2}}|-|^{\frac{d-1}{2}} \simeq \sigma_d(\pi_\rho^I)|-|^{\frac{d-1}{2}} 
\end{eqnarray*}
par la description \ref{corell} de la correspondance de Langlands pour les repr\'esentations elliptiques.

\begin{rema}
\'Ecrivons $I=\{i_1,\cdots, i_{|I|}\}$ dans l'ordre croissant et posons $i_0:=0$, $i_{|I|+1}=d_\rho$ et  $d_k:=i_{k+1}-i_k$ pour $k=0,\cdots, |I|$. Alors la discussion de \cite[4.4]{Dat1} montre qu'en tant que $W_K$-repr\'esentation {\em gradu\'ee}, on a
$$ \HC_\rho^I [1-d_\rho] \simeq \bigoplus_{k=0}^{|I|} \left(\sigma_{d/d_\rho}(\pi_\rho)\otimes \tau_{d_k}^\emptyset |.|^{\frac{d-d_\rho}{2}+i_k}\right) [-|I|+2k]. $$
Ainsi les composantes sont ind\'ecomposables et rang\'ees par ordre croissant de poids de 2 en 2.
\end{rema}

\section{Monodromie et vari\'et\'es de Shimura} \label{monodromie}

Un des buts de cette section est de prouver la proposition \ref{estimN} et donc de terminer la preuve du th\'eor\`eme \ref{main}. Ceci fera intervenir le syst\`eme des vari\'et\'es de Shimura \'etudi\'ees par Harris-Taylor (ou son analogue d\^u \`a Drinfeld-Stuhler en \'egales caract\'eristiques), ainsi que la description de la filtration de monodromie de leurs cycles \'evanescents par Boyer.

Dans la deuxi\`eme partie, on prouve la conjecture monodromie-poids pour les vari\'et\'es uniformis\'ees par les rev\^etements de l'espace sym\'etrique de Drinfeld.
Ceci s'applique \`a certaines vari\'et\'es de Shimura unitaires \'etudi\'ees notamment par Harris \cite{HaCusp}, Carayol \cite{CarAA} et Rapoport \cite{RapAA}.

On note toujours $G=GL_d(K)$ et $D=D_d$.
\def\dd{{D^\times}}

\subsection{Vari\'et\'es de Harris-Taylor}

On s'int\'eresse ici aux vari\'et\'es de Shimura \'etudi\'ees par Harris et Taylor dans \cite{HaTay}. Plus pr\'ecis\'ement, on suppose que $K$ est le compl\'et\'e en une place $w$ d'un corps CM $F$ v\'erifiant les hypoth\`eses de \cite[I.7]{HaTay}, on fixe un ``niveau hors $p$" $U^p\subset \GC(\AM^{\infty,p})$ d'un certain groupe unitaire $\GC$ d\'efini aussi en \cite[I.7]{HaTay} et qui \`a la place $w$ est isomorphe \`a $GL_d$, puis pour un entier $n$, on note $S_{HT,n}$ le $\OC_K$-sch\'ema propre et r\'egulier not\'e $X_{U^p,m}$ dans \cite[III.4]{HaTay}, avec $m=(n,0,\cdots,0)$. On notera avec des exposants $\eta$ ou $s$ les fibres g\'en\'eriques ou sp\'eciales de ces objets et $j$ et $i$ les immersions correspondantes, d\'ecor\'ees des indices pertinents.

L'entier $n$ correspond \`a une structure de niveau (\`a la Drinfeld-Katz-Mazur) $H_n=1+\varpi^nM_d(\OC_K)$ en $w$ sur le probl\`eme de modules d\'efinissant la vari\'et\'e de Shimura sur $\OC_K$. Par des  constructions analogues \`a celles utilis\'ees pour la tour de Lubin-Tate, la famille des $(S_{HT,n})_{n\in\NM}$  est munie d'un syst\`eme de morphismes finis $g^{m|n}$. Une mani\`ere agr\'eable de formaliser ceci est d'introduire la cat\'egorie suivante :

\begin{DEf}
On note $\NM(G)$  la cat\'egorie dont les objets sont les entiers naturels et les fl\`eches sont donn\'ees par 
$$\hom{n}{m}{\NM(G)} :=H_m\ba \{g\in G, gH_ng^{-1} \subset H_m\}/H_n ,$$
la composition \'etant induite par la multiplication dans $G$. On note aussi $\NM(G^0)$ la sous-cat\'egorie dont les fl\`eches viennent de $G^0=\ker |\hbox{det}|_K$.
\end{DEf}
Avec cette d\'efinition, la famille des $S_{HT,n}$ 
est l'image d'un $\NM(G)$-diagramme, {\em i.e.} un foncteur de $\NM(G)$ dans les $\OC_K$-sch\'emas. Cette formulation invite \`a utiliser le langage des cat\'egories (co)-fibr\'ees. Par exemple, on dispose de la cat\'egorie $\hbox{Perv}(S_{HT}^{s,ca},\o\QM_l)$ au-dessus de $\NM(G)$, dont les fibres sont les cat\'egories de $\o\QM_l$-faisceaux pervers  sur les $S_{HT,n}^{s,ca}$ et qui est fibr\'ee par les ${^{p}}g^{m|n*}$ et cofibr\'ee par les $g^{m|n}_*$.
Les sections de cette cat\'egorie bi-fibr\'ee seront appel\'ees  {\em Faisceaux Pervers de Hecke}. Ce sont donc des syst\`emes $(\FC_n)_n$ de faisceaux pervers munis de morphismes de transition $\FC_m\To{} g^{m|n}_*\FC_n$.  En vertu de l'exactitude des $g^{m|n}_*$, ils forment une cat\'egorie  ab\'elienne, que nous noterons $\hbox{FPH}$.

Les r\'esultats globaux principaux de Boyer dans \cite{Boyer2} concernent le faisceau pervers de Hecke $R\Psi$ form\'e par le syst\`eme des cycles \'evanescents d\'ecal\'es $R\Psi(S_{HT,n}^\eta,\o\QM_l)[d-1]$.
Celui-ci est muni d'une action de $W_K$ compatible avec celle sur $k^{ca}$. En particulier, on a pour chaque \'etage un op\'erateur de monodromie nilpotent $N_n$, deux filtrations $K_\bullet$ (croissante) et $I^\bullet$ (d\'ecroissante)
d\'efinies respectivement par les noyaux %$K_p:=\ker N_n^{p+1}$ 
et les images %$I^q:=\im N_n^q$ 
des it\'er\'es de $N_n$, ainsi que leur convolution $M_\bullet$, appel\'ee ``filtration de monodromie". La famille des $N_n$ d\'efinit un endomorphisme $W_K$-\'equivariant $N$ de $R\Psi$  et les familles de filtrations induisent des filtrations $W_K$-\'equivariantes $K_\bullet R\Psi$, $I^\bullet R\Psi$  et $M_\bullet R\Psi$ de l'objet $R\Psi$ dans $FPH$.
Remarquons que $N$ n'est {\em a priori} pas nilpotent, et que ces filtrations ne sont {\em a priori} pas finies. Cependant on a :
\begin{fact} (Boyer)
Pour tout $n\in\NM$, l'action de $I_K$ sur le $\o\QM_l$-faisceau constructible $R^j\Psi(S^\eta_{HT,n},\o\QM_l)$ se factorise par un quotient fini.
\end{fact}
Il s'ensuit que l'op\'erateur $N_n$ est nul sur les faisceaux de cohomologie, et puisque $R\Psi(S^\eta_{HT,n},\o\QM_l)$ est d'amplitude cohomologique $[1-d,0]$, que $N_n$ est nilpotent d'ordre $\leq d$. Par cons\'equent $N$ est aussi nilpotent d'ordre $\leq d$ et les filtrations de $R\Psi$ sont finies. 
On peut pr\'eciser cela en d\'ecoupant, suivant Boyer,  les cycles \'evanescents selon l'action de l'inertie (voir aussi la preuve de \ref{subtilite}) :
$$ R\Psi = \bigoplus_{\sigma\in\Irr{\o\QM_l}{W_K}/\sim} R\Psi_\sigma, $$
d\'ecomposition dans la cat\'egorie $\hbox{FPH}$ et o\`u les $\sigma$ sont prises \`a \'equivalence inertielle pr\`es. Alors il r\'esulte de \cite[5.4.7]{Boyer2} que $R\Psi_\sigma$ est d'amplitude cohomologique $[1-d_\sigma,0]$ o\`u $d_\sigma:= \lfloor d/\dim(\sigma)\rfloor $
%Il d\'ecrit ensuite \cite[Thm 5.4.5]{Boyer2} les gradu\'es de la filtration de monodromie $M_\bullet R\Psi_\sigma$, mais il ne sera pas difficile de d\'eduire de son travail les bigradu\'es pour la bifiltration $I^\bullet K_\bullet (R\Psi_\sigma)$, {\em cf} la preuve de \ref{bigrad} ci-dessous. 
%De sa description r\'esulte en particulier que  
et donc que l'op\'erateur de monodromie $N_\sigma$ est nilpotent d'ordre $\leq d_\sigma$, %:= \lfloor d/\dim(\sigma)\rfloor $, 
 que les gradu\'es $\gr^q_I := I^q/I^{q+1}( R\Psi_\sigma)$ et $\gr_p^K:= K_{p+1}/K_p( R\Psi_\sigma)$ sont nuls en dehors de $0\leq p,q < d_\sigma$ et que les gradu\'es $\gr_k^M:=M_{k+1}/M_k(R\Psi_\sigma)$ sont nuls en dehors de $-d_\sigma<k<d_\sigma$. %De plus, lorsque $\sigma=\sigma'(\rho)$ pour $\rho\in \Irr{\o\QM_l}{\dd}$, alors on a $d_\sigma=d_\rho$.
Le th\'eor\`eme 5.4.5 de \cite{Boyer2} d\'ecrit les gradu\'es de la filtration de monodromie $M_\bullet R\Psi_\sigma$, et il n'est pas difficile d'en d\'eduire  les bigradu\'es de la bifiltration $I^\bullet K_\bullet (R\Psi_\sigma)$, {\em cf} la preuve de \ref{bigrad} ci-dessous.
\ali
Bien-s\^ur, le lien entre les r\'esultats globaux de Boyer et les r\'esultats locaux dont nous avons besoin se fait en prenant la ``fibre en un point supersingulier" de $R\Psi_\sigma$.
Plus pr\'ecis\'ement, donnons-nous un point supersingulier $x$ dans $S^{s,ca}_{HT,0}$. On sait que les morphismes $g^{0|n}$ sont  totalement ramifi\'es au-dessus de $x$ et que si $g\in G^0$, la pr\'eimage r\'eduite de $x$ dans $S^{s,ca}_{HT,n}$ est un singleton $\{x_n\}$ ind\'ependant de $g$. Le syst\`eme des $x_n^*(R\Psi_{n,\sigma})$ est donc un foncteur contravariant  de $\NM(G^0)$ dans $D^b(\o\QM_l)$ (on identifiera cette derni\`ere \`a la cat\'egorie des espaces vectoriels \`a graduation finie). En passant \`a la limite inductive, on obtient un espace vectoriel gradu\'e muni d'une action de $G^0$. D'o\`u un foncteur additif
$$ x^*:\;\; FPH \To{} \gr\Mo{\o\QM_l}{G^0}:=\{\o\QM_l\hbox{-repr\'esentations lisses gradu\'ees de } G^0\}. $$

On sait que le compl\'et\'e formel (\'etage par \'etage) de la tour $S^{nr}_{HT}$ le long de la famille $x_n$ s'identifie \`a $\MC_{LT}^{(0)}$, et par le th\'eor\`eme de comparaison de Berkovich, on a donc un isomorphisme $G^0\times I_K$-\'equivariant
$$
 x^*(R\Psi) \simeq \bigoplus_i H^i(\MC_{LT}^{(0),ca},\o\QM_l)[d-1-i]. 
$$
De plus, l'action de $\OC_D^\times$ par automorphismes du compl\'et\'e formel de chaque $S^{nr}_{HT,n}$ en $x_n$ munit $x^*(R\Psi)$ d'une action de $\OC_D^\times$ {\em a priori}, et l'isomorphisme de Berkovich ci-dessus est \'equivariant, puisque canonique.
En fait, l'action de $G^0\times \OC_D^\times \times I_K$ sur l'espace vectoriel gradu\'e $x^*(R\Psi)$ est la restriction d'une action de $(G\times \dd\times W_K)^0$ que l'on obtient en  consid\'erant le compl\'et\'e formel de 
$S^{s,ca}_{HT}$ le long de l'orbite Galoisienne $X=(X_n)_{n\in\NM}$ de $x$ sous $W_K$. En effet, celle-ci est un sous $\NM(G)$-diagramme de $S^{s,ca}_{HT}$ de dimension $0$, dont le compl\'et\'e formel associ\'e est muni d'une action de $\dd$. Le stabilisateur du compl\'et\'e formel en $x$ est justement $(G\times \dd\times W_K)^0$.

Pour d\'ecrire $x^*(R\Psi_\sigma)$, il est plus commode de modifier le foncteur $x^*$ en posant :
$$ ix^*:= \cInd{G^0\varpi^\ZM}{G}\circ x^*:\;\;FPH \To{} \gr\Mo{\o\QM_l}{G/\varpi^\ZM}$$
o\`u $\varpi$ est une uniformisante de $K^\times \subset G$. Alors l'objet $ix^*(R\Psi_\sigma)$ est un facteur direct stable par $G\times \dd\times W_K$ de $H^*(\MC^{ca}_{LT}/\varpi^\ZM,\o\QM_l)$
dont la restriction \`a $I_K$ est $\sigma$-isotypique.
D'apr\`es la description de Boyer \ref{conjHarris} et la dualit\'e \ref{dualite}, on a donc : %on a donc pour toute repr\'esentation $(\rho,V_\rho)\in\Irr{\o\QM_l}{\dd}$ un isomorphisme $(G\times \dd\times W_K)^0$-\'equivariant d'espaces vectoriels gradu\'es
\def\cor{\rightsquigarrow}
\ini
\begin{eqnarray} \label{rpsisigma}
 ix^*(R\Psi_{\sigma}) & \mathop{\simeq}\limits_{G\times \dd\times W_K} & \left(\bigoplus_{\rho\cor\sigma} \left(\bigoplus_i  H^i_c[\rho][-i]\right)^\vee \otimes \rho\right)\,[1-d](1-d) \\
 & \mathop{\simeq}\limits_{G\times \dd\times W_K} & \bigoplus_{\rho\cor\sigma} \left( \bigoplus_{i=0}^{d_\sigma-1}  \pi_{\rho^\vee}^{> i} \otimes \rho\otimes \sigma'(\rho) (-i)[-i]\right)[d_\sigma-1] .
   \end{eqnarray}
o\`u $\rho\cor\sigma$ d\'esigne la conjonction des deux conditions $\omega_\rho(\varpi)=1$ et $\sigma'(\rho)\sim \sigma$ (\'equivalence inertielle) et $\sigma'(\rho)$ est d\'efinie en \ref{defsigmaprime}. En particulier $x^*(R\Psi_\sigma)$ est nul si le degr\'e de $\sigma$ ne divise pas $d$.

On aimerait transf\'erer la bifiltration de $R\Psi_\sigma$. Pour cela, on remarque que le foncteur $x^*$ est ``exact" dans le sens suivant : 
\begin{lemme} \label{exact}
Le foncteur $x^*$ se prolonge en un foncteur de la cat\'egorie des suites exactes 
de $FPH$ dans la cat\'egorie des triangles distingu\'es %$G$-\'equivariants
de $\gr\Mo{\o\QM_l}{G^0}$. %$\o\QM_l G^0$-repr\'esentations lisses gradu\'ees. %'espaces vectoriels gradu\'es.
\end{lemme}
En d'autre termes, une suite exacte $A\injo B\twoheadrightarrow C$ induit 
%$x^*A\To{} x^*B \To{} x^*C \To{+1}$ form\'e de fl\`eches $G$-\'equivariantes telles que la  
une longue suite exacte
%soit distingu\'e %{\em au sens des 
%}, ce qui signifie que la longue suite d'espaces vectoriels
$$ \cdots \To{} x^*A^i\To{} x^*B^i \To{} x^*C^i \To{} x^*A^{i+1} \To{} \cdots .$$

\begin{proof}
Notons temporairement $D^b_c$ la cat\'egorie des sections de la cat\'egorie bi-fibr\'ee au-dessus de $\NM(G)$ dont les fibres sont les $D^b_c(S^{s,ca}_{HT,n},\o\QM_l)$. La cat\'egorie $FPH$ est une sous-cat\'egorie pleine de $D^b_c$, et le foncteur $x^*$ est la restriction d'un foncteur de source $D^b_c$, encore not\'e $x^*$ et d\'efini de la m\^eme mani\`ere.
La cat\'egorie $D^b_c$ est additive, $\ZM$-gradu\'ee, et munie d'une famille \'evidente de triangles distingu\'es : ceux qui \`a chaque \'etage le sont. Le foncteur $x^*$ est exact, {\em i.e.} envoie triangles distingu\'es sur triangles distingu\'es au sens du lemme. Ainsi, il nous suffira de montrer que toute suite exacte $A\injo B \twoheadrightarrow C$ se compl\`ete de mani\`ere unique en un triangle distingu\'e de $D^b_c$. Or, \`a chaque \'etage, on connait l'unicit\'e de $\delta_n : C_n\To{} A_n[1]$ compl\'etant la suite $A_n\To{} B_n \To{} C_n$ en un triangle distingu\'e \cite[Cor 1.1.10 ii)]{BBD}. Cette unicit\'e et la $p$-exactitude des $g^{m|n}_*$ assurent que le syst\`eme $(\delta_n)_n$ est bien un morphisme dans $D^b_c$.
\end{proof}

%Soit $\wt\ZM$ la cat\'egorie associ\'ee \`a l'ensemble ordonn\'e $\ZM\cup\{-\infty,\infty\}$.
%D'apr\`es le lemme ci-dessus, la filtration de monodromie $M_\bullet R\Psi_\sigma$ d\'efinit un {\em objet spectral} au sens de Verdier \cite[II.4]{Verdier} de type  $\wt\ZM$ et \`a valeurs dans $\Mo{\o\QM_l}{G^0}$,
%dont la suite spectrale associ\'ee calcule la cohomologie de $x^*(R\Psi_\sigma)$ munie de son action de $G^0$

%Soit  $\wt\ZM^2$ la cat\'egorie dont les objets sont les couples $(p,q)$ dans $\ZM\cup \{-\infty,+\infty\}$ et les morphismes sont ceux associ\'es \`a l'ordre produit de l'ordre croissant en $p$ et d\'ecroissant en $q$. 
D'apr\`es le lemme ci-dessus, la bifiltration de $R\Psi_\sigma$ induit  une ``bifiltration" sur son image $x^*R\Psi_\sigma$ dans $\gr\Mo{\o\QM_l}{G^0}$ (en termes rigoureux, un  {\em objet spectral} au sens de Verdier \cite[II.4]{Verdier}). % de type  $\wt\ZM^2$  %\`a valeurs dans $D^b(\o\QM_l)$, et muni d'une action lin\'eaire de $G^0$. 
%et \`a valeurs dans $\Mo{\o\QM_l}{G^0}$.
Par d\'efinition, elle est $I_K$-\'equivariante, et par le th\'eor\`eme 2 de \cite{FarMono}, elle est aussi $\OC_D^\times$-\'equivariante. En fait, elle est $(G\times\dd\times W_K)^0$-\'equivariante, comme on le voit en r\'ep\'etant ces arguments pour l'orbite Galoisienne $X$ de $x$, puis en se restreignant \`a $x$.

La suite spectrale associ\'ee \`a la filtration de monodromie de $x^*R\Psi_\sigma$ est enti\`erement d\'ecrite par Boyer. Voici une reformulation de cette description, en termes de bifiltration, et adapt\'ee \`a nos besoins ; par convention, nos gradu\'es sont donn\'es par $\gr_I^q:=I^q/I^{q+1}$ et $\gr^K_p:=K_{p+1}/K_p$.

\begin{theo} \label{bigrad} (Boyer) Soit $\sigma\in \Irr{\o\QM_l}{W_K}$, et rappelons que 
%Pour tous $p,q \in \ZM$, posons $i_\sigma^{p,q}:= d_\sigma -1 -p-q$ o\`u comme ci-dessus,  
$d_\sigma=\lfloor d/\dim(\sigma)\rfloor$.  
\begin{enumerate}
\item Le $\o\QM_l$-espace gradu\'e $x^*(\gr_I^q\gr^K_p R\Psi_{\sigma})$ est non-nul seulement si $p,q\geq 0$ et $p+q<d_\sigma$, auquel cas il est concentr\'e en degr\'e $p+q-d_\sigma+1$ et donn\'e par :
$$ ix^*(\gr_I^q\gr^K_p R\Psi_{\sigma})[p+q-d_\sigma+1] \mathop{\simeq}\limits_{G\times\dd\times W_K} %\left\{\begin{array}{ll} %\pi_{\rho^\vee}^{\{p+q,\cdots,d_\rho-1\}}
\bigoplus\limits_{\rho\cor \sigma} \tau_{\rho^\vee}^{>{p+q}}\otimes\rho\otimes \sigma'(\rho)(q-p) % & \hbox{si }  p,q, i_\sigma^{p,q}\geq 0   \\
            %0 & \hbox{sinon} \end{array} \right.
$$
o\`u  $\tau_\rho^{>i}$ est l'"unique" extension non triviale de $\pi_\rho^{>i+1}$ par $\pi_\rho^{>i}$ dans $\Mo{\omega_\rho}{G}$ pour $0<i <d_\sigma$ et  $\tau_\rho^{>d_\sigma}:=\pi_\rho^{>d_\sigma}=\pi_\rho^\emptyset$.

\item Pour tous $p,q\geq 0$ tels que $p+q<d_\sigma-1$, le morphisme $G\times\dd\times W_K$-\'equivariant 
$$ix^*(\gr_I^q\gr^K_p R\Psi_{\sigma})[p+q-d_\sigma+1] \To{} ix^*(\gr_I^{q+1}\gr^K_p R\Psi_{\sigma})[p+q+1-d_\sigma+1] $$
d\'eduit par d\'ecalage et rotation du triangle distingu\'e $\gr_I^{q+1}\To{} I^q/I^{q+2} \To{} \gr_I^q$ est donn\'e (modulo isomorphismes) sur chaque facteur direct par l'"unique" morphisme $G$-\'equivariant non-nul $\tau_{\rho^\vee}^{>p+q} \To{} \tau_{\rho^\vee}^{>p+q+1}$.
\end{enumerate}
\end{theo}

\begin{proof}
L'assertion sur la nullit\'e en dehors du triangle $p,q\geq 0$ et $p+q<d_\sigma$ est une simple cons\'equence de la d\'efinition des filtrations $I^\bullet $ et$K_\bullet $ et du fait que $N^{d_\sigma}=0$. Nous supposons dor\'enavant que ces in\'egalit\'es sont v\'erifi\'ees et nous allons d'abord montrer que $ix_{eq}^*(\gr_I^q\gr^K_p)$ est concentr\'e en degr\'e $p+q-d_\sigma+1$. Pour cela, le plus commode est de se raccrocher 
\`a la description des $\gr^M_k (R\Psi_\sigma)$ par Boyer. Par d\'efinition, on a $\gr^M_k=\bigoplus_{p-q=k} \gr_I^q\gr^K_p$. Inversement, pour r\'ecup\`erer les bigradu\'es \`a partir du gradu\'e de monodromie, il faut se rappeler que l'op\'erateur $N_\sigma$  envoie $I^qK_p(R\Psi_\sigma)$ dans $I^{q+1}K_{p-1}(R\Psi_\sigma)$ et induit un isomorphisme 
$$ N_\sigma : \gr_I^{q}\gr^K_p(R\Psi_\sigma) \simto \gr_I^{q+1}\gr^K_{p-1}(R\Psi_\sigma) $$
tant que $p>1$. On en d\'eduit les formules :
\ini\begin{equation}\label{forgr}
 \gr_I^q\gr_0^K = \ker\left(\gr_{-q}^M \To{N_\sigma} \gr_{-q-2}^M\right) \;\; \hbox{ et } \;\; N_\sigma^{p}: \gr_I^q\gr_p^K \simto \gr_I^{q+p}\gr_0^K. 
 \end{equation}
D'apr\`es Boyer \cite[5.4.5]{Boyer2}, les gradu\'es de monodromie sont de la forme
$$ \gr^M_k(R\Psi_\sigma) = \mathop{\bigoplus_{|k|\leq t\leq d_\sigma}}\limits_{t \equiv k-1 [2]} \PC(t,\sigma)(?) $$
o\`u $\PC(t,\sigma)$ est le faisceau pervers de Hecke not\'e $\PC(g,t,\pi_v)$ dans {\em loc. cit} avec le dictionnaire suivant : $g \longleftrightarrow \dim(\sigma)$ et $\pi_v \longleftrightarrow \sigma_g^{-1}(\sigma)$, et $(?)$ d\'esigne une torsion \`a la Tate qui ne nous importe pas ici. De par leur d\'efinition, chaque $\PC(t,\sigma)$ est semi-simple de longueur finie et sans multiplicit\'e. De plus, les $\PC(t,\sigma)$ sont deux-\`a-deux ``disjoints". On d\'eduit alors des formules \ref{forgr} et par une r\'ecurrence facile que 
$$ \gr_I^q\gr^K_p (R\Psi_\sigma) = \PC(p+q+1 , \sigma)(?) .$$
Le th\'eor\`eme 5.4.7 de \cite{Boyer2} d\'ecrit (entre autres) la cohomologie de $x^*(\PC(t,\sigma))$ et montre que celle-ci est concentr\'ee en degr\'e $t-d_\sigma$, d'o\`u la premi\`ere assertion du point i).

Consid\'erons maintenant la suite spectrale $(G\times\dd\times W_K)^0$-\'equivariante associ\'ee \`a la filtration de monodromie  
$$ E_1^{i,j}:= \HC^{i+j}(x^*(\gr^M_{-i}(R\Psi_\sigma))) \Rightarrow \HC^{i+j}(x^*R\Psi_\sigma). $$
Du simple fait que la cohomologie de $x^*\gr^q_I\gr^K_p$ est concentr\'ee en degr\'e $p+q+1-d_\sigma$, on tire pour tous $p,q$ :
\begin{enumerate}
        \item $x^*(\gr_I^p\gr^K_q(R\Psi_\sigma))[p+q+1-d_\sigma] = E^{q-p,2p+1-d_\sigma}_1$,
        \item La fl\`eche du point ii) est la diff\'erentielle $d_1^{q-p,2p+1-d_\sigma}$.
\end{enumerate}
Les assertions restantes sont donc cons\'equences du th\'eor\`eme 4.2.3 de \cite{Boyer2}, puisque dans le dictionnaire entre nos notations et celles de {\em loc. cit}, on a
$\tau^{>i}_\rho \longleftrightarrow [\overleftarrow{i}]_{\pi_\rho}\times[\overrightarrow{s-i}]_{\pi_\rho}$.
%D'apr\`es Boyer, \cite[Thm 4.2.4, 4.2.4]{Boyer2}, on a $\hom{\rho}{}{\dd}$

\end{proof}

%La suite spectrale associ\'ee \`a cet objet spectral calcule la cohomologie de $x^*(R\Psi_\sigma)$ munie de son action de $G^0$. En consid\'erant plut\^ot la filtration de monodromie $M_\bullet$ sur $R\Psi_\sigma$, on obtient un objet spectral de type $\wt\ZM$ dont la suite spectrale est enti\`erement d\'ecrite par Boyer...

%Ce r\'esultat peut paraitre inutilisable pour prouver la proposition \ref{estimN} puiqu'il ne fait pas intervenir la cat\'egorie d\'eriv\'ee $D^b_{\o\QM_l}(G^0)$. En fait l'op\'erateur $x^*(N_{\sigma})$ induit par la monodromie de $R\Psi_{\sigma}$ est nul, puisque nul en  cohomologie.
Pour pouvoir utiliser ce th\'eor\`eme, il nous faut maintenant relever le foncteur $x^*$ vers la cat\'egorie $D^b_{\o\QM_l}(G^0)$. Malheureusement, cela pose plusieurs probl\`emes techniques. 

Le premier vient du formalisme $l$-adique. Pour l'expliquer, rappelons que pour $X$ sch\'ema de type fini sur un corps de dimension cohomologique finie, la cat\'egorie $D^b_c(X,\o\QM_l)$ est d\'efinie comme limite inductive des cat\'egories $D^b_c(X,\Lambda)[\frac{1}{l}]$ o\`u $\Lambda$ d\'ecrit les anneaux d'entiers d'extensions finies de $\QM_l$ et les morphismes de transition sont donn\'es par extension des scalaires de $\Lambda[\frac{1}{l}]$ \`a $\Lambda'[\frac{1}{l}]$. De plus, chaque $D^b_c(X,\Lambda)[\frac{1}{l}]$ s'identifie 
naturellement \`a une sous-cat\'egorie triangul\'ee pleine de la cat\'egorie 
$D^+_{\Lambda_\bullet}(\wt{X_{et}})[\frac{1}{l}]$ (via les complexes ``normalis\'es" d'Ekedahl, {\em cf} \cite[Prop 2]{FarMono}). Enfin, chaque
 $D^+_{\Lambda_\bullet}(\wt{X_{et}})$ est munie d'une $t$-structure perverse \cite[Thm 7]{FarMono} qui induit sur $D^b_c(X,\Lambda)[\frac{1}{l}]$ la $t$-structure perverse interm\'ediaire \cite[Cor 3]{FarMono} et pour laquelle les morphismes finis sont $t$-exacts.
Notons alors
$FPH_{\Lambda_\bullet}$ la cat\'egorie des sections de la cat\'egorie bi-fibr\'ee sur $\NM(G)$ dont la fibre en $n$  est le coeur de la $t$-structure perverse sur $D^+_{\Lambda_\bullet}(\wt{S_{HT,n,et}})$. Par $t$-exactitude des foncteurs de transition $g^{m|n}_*$, cette cat\'egorie est ab\'elienne. 
Soit $FPH_{\Lambda_\bullet}^c$ sa sous-cat\'egorie pleine (ab\'elienne et \'epaisse) form\'ee des sections %dont l'image par le foncteur de localisation $[\frac{1}{l}]$
%$D^{+,naif}_{\Lambda_\bullet}(S_{HT})$ la cat\'egorie des sections $(K_n)_{n\in\NM(G)}$ de la cat\'egorie bi-fibr\'ee  dont la fibre en $n\in\NM(G)$  est $D^b_{\Lambda_\bullet}(\wt{S_{HT,n,et}})$, et
%$FPH_{\Lambda_\bullet}^c$ la sous-cat\'egorie pleine form\'ee des sections pour lesquelles l'image de chaque 
$(K_n)_{n\in\NM(G)}$ avec $K_n[\frac{1}{l}] \in D^b_c(S_{HT,n},\Lambda)[\frac{1}{l}]$.
 Plut\^ot qu'avec la cat\'egorie $FPH$ nous travaillerons avec la cat\'egorie ab\'elienne $FPH'$ suivante :
$$ FPH':=\limi{\Lambda} (FPH_{\Lambda_\bullet}^c [\frac{1}{l}]). $$
On a un foncteur \'evident  $FPH' \To{\iota} FPH$, exact et fid\`ele, dont le d\'efaut de pleinitude vient de ce que cette construction ``borne" les d\'enominateurs et l'alg\'ebricit\'e lorsqu'on bouge dans le diagramme.

Par d\'efinition, le complexe des cycles \'evanescents $R\Psi$ est naturellement un objet de $FPH'$, dont la valeur en $n\in\NM(G)$ est $(R\Psi(S_{HT,n},\Lambda_\bullet))_\Lambda$. 
%Ce qui n'est p
Par contre,
% pas du tout \'evident {\em a priori}, c'est de savoir si
 son op\'erateur de monodromie $N$ n'est {\em a priori} pas un endomorphisme de $R\Psi$ dans $FPH'$, car il est d\'efini \`a chaque \'etage par un logarithme  dont les d\'enominateurs croissent avec le niveau $n$. %pourraient \^etre non born\'es. 
Cependant on a le lemme suivant :

\begin{lemme} \label{subtilite}
Soit $\sigma \in \Irr{\o\QM_l}{W_K}$.
Le plongement  $R\Psi_\sigma \injo R\Psi$ dans $FPH$ et l'op\'erateur de monodromie $N_\sigma$ sont dans l'image essentielle de $\iota$.
\end{lemme}

\begin{proof}
Commen{\c c}ons par expliciter la construction de $R\Psi_\sigma$ dans la cat\'egorie $FPH$. Notons pour cela $\gamma: I_K \To{} \endo{FPH}{R\Psi}$ l'action de l'inertie. Par d\'efinition de $N$, l'application $\wt\gamma : i\in I_K \mapsto \gamma(i) \hbox{exp}(-t_l(i)N) $ d\'efinit une action localement constante de $I_K$ sur $R\Psi$ et se prolonge donc en un morphisme de $\o\QM_l$-alg\`ebres $\HC_{\o\QM_l}(I_K) \To{} \endo{FPH}{R\Psi}$ o\`u $\HC_{\o\QM_l}(I_K)$ d\'esigne l'alg\`ebre des distributions localement constantes sur $I_K$. Le facteur $R\Psi_\sigma$ est d\'ecoup\'e via $\wt\gamma$ par l'idempotent de $\HC_{\o\QM_l}(I_K)$ associ\'e \`a $\sigma$.

Choisissons une base  $(I_m)_{m\in\NM}$ de voisinages de l'unit\'e de $I_K$ form\'ee de sous-groupes ouverts de $I_K$ normaux dans $W_K$, et notons $R\Psi_m$ le facteur direct de $R\Psi$ dans $FPH$ d\'ecoup\'e via $\wt\gamma$ par la repr\'esentation triviale de $I_m$. On a une suite exacte $I_m^l\injo I_m \stackrel{t_l}{\twoheadrightarrow} l^{m'}\ZM_l$ avec $I_m^l$ un pro-$l'$-groupe. Comme la restriction de l'action $\gamma$ \`a $I_m^l$ est localement constante, on a un idempotent $\gamma[I_m^l]$ o\`u $[I_m^l]\in \HC_{\ZM_l}(I_m^l)$ est la mesure de Haar normalis\'ee. 
Choisissons un \'el\'ement $T_m\in I_m$  tel que $t_l(T_m)=l^{m'}$. On a alors
$$R\Psi_m= \sum_{t\in\NM} \ker \left((\gamma(T_m)-1)^t_{|\im\gamma[I_m^l]}\right)$$
et son op\'erateur de monodromie est donn\'e par
$$ N_{|R\Psi_m}= \frac{1}{l^{m'}} \log(\gamma(T_m)). $$
On a d\'eja remarqu\'e que $N$ est nilpotent d'ordre $\leq d$.
%d'apr\`es la description de Boyer, $N_\sigma$ est nilpotent d'ordre $d_\sigma$. Il s'ensuit que $N_{|R\Psi_m}$ est nilpotent d'ordre $ d_m=\sup_{I_m\subset \ker \sigma} d_\sigma$, et par cons\'equent on a aussi 
On a donc aussi
$$ R\Psi_m= \ker\left((\gamma(T_m)-1)^{d}_{|\im\gamma[I_m^l]}\right).$$
Or, les endomorphismes $\gamma[I_m^l]$ et $\gamma(T_m)$ de $R\Psi$ sont dans $\endo{FPH'}{R\Psi}$. Ceci montre l'\'enonc\'e du lemme pour $R\Psi_m$ \`a la place de $R\Psi_\sigma$ et on en d\'eduit aussit\^ot l'\'enonc\'e du lemme puisque $R\Psi_\sigma=(R\Psi_m)_\sigma$ d\`es que $I_m\subset \ker\sigma$.

\end{proof}

\begin{prop} \label{equiv}
Il existe un foncteur $x^*_{eq}$ s'inscrivant dans un diagramme essentiellement commutatif
$$\xymatrix{\hbox{FPH'} \ar[r]^-{x^*_{eq}} \ar[d]_{\iota} & D^b_{\o\QM_l}(G^0) \ar[d]^{\oplus \HC^i[-i]} \\ \hbox{FPH} \ar[r]_-{x^*} & \gr\Mo{\o\QM_l}{G^0}},$$
exact au sens du lemme \ref{exact} et 
 tel qu'on ait un isomorphisme canonique et compatible \`a l'action de $W_K$
$$x^*_{eq} R\Psi \simeq R\Gamma(\MC_{LT}^{(0),ca},\o\QM_l)[d-1]\;\;\hbox{ dans } \;D^b_{\o\QM_l}(G^0).$$
%${x^*:\;\;}{F\hbox{FPH}}{D^bF_{\o\QM_l}(G^0)}$ 
\end{prop}
Admettons ce r\'esultat momentan\'ement et posons $ix^*_{eq}:=\cInd{G^0\varpi^\ZM}{G}\circ x^*_{eq}$, qui est donc un foncteur ``exact" de $FPH'$ dans $D^b_{\o\QM_l}(G/\varpi^\ZM)$ et envoie $R\Psi$ sur $R\Gamma(\MC^{ca}_{LT}/\varpi^\ZM,\o\QM_l)[d-1]$ en respectant les actions de $W_K$. D'apr\`es \ref{rpsisigma},
% la description de Boyer \ref{conjHarris} et la dualit\'e \ref{dualite}, on d\'eduit de l'isomorphisme de la proposition   pour toute repr\'esentation $\sigma\in\Irr{\o\QM_l}{W_K}$ un isomorphisme $ I_K$-\'equivariant dans $D^b_{\o\QM_l}(G/\varpi^\ZM)$
on a un isomorphisme $I_K$-\'equivariant
$$ ix^*_{eq}(R\Psi_{\sigma}) \simeq \bigoplus_{\rho\cor\sigma} \left(R\Gamma_c[\rho]^\vee \otimes V_\rho\right) [1-d](1-d), $$
l'action de $I_K$ sur le $\o\QM_l$-espace $V_\rho$ \'etant triviale.

\alin{Preuve de la proposition \ref{estimN}}
Fixons $\rho\in\Irr{\o\QM_l}{\dd}$ comme dans la partie \ref{real} et posons $\sigma:=\sigma'(\rho)$, de sorte que $d_\sigma=d_\rho$.
Par l'isomorphisme pr\'ec\'edent, il suffit de prouver que $$ix_{eq}^*(N_{\sigma})^{d_\sigma-1}\neq 0.$$
En effet, cela impliquera $N_{\rho'}^{d_\rho-1}\neq 0$ pour au moins un $\rho'\sim \rho$, ce qui, par torsion, impliquera alors $N_\rho^{d_\rho-1}\neq 0$.

Par d\'efinition des filtrations, l'op\'erateur $N_\sigma$ envoie $I^qK_p(R\Psi_\sigma)$ dans $I^{q+1}K_{p-1}(R\Psi_\sigma)$ et induit un isomorphisme 
$$ N_\sigma : \gr_I^{q}\gr^K_p(R\Psi_\sigma) \simto \gr_I^{q+1}\gr^K_{p-1}(R\Psi_\sigma) $$
tant que $p>1$. En particulier, puisque les $\gr_I^qR\Psi_{\sigma}$ et $\gr^K_pR\Psi_{\sigma}$ sont nuls pour $p,q$ hors de $\{0,\cdots, d_\sigma-1\}$, on a un diagramme
$$ \xymatrix{ R\Psi_{\sigma} \ar[r]^{N_{\sigma}^{d_\sigma-1}} \ar[d]_{can_+} & R\Psi_{\sigma} \\ \gr_I^0\gr^K_{d_\sigma-1}R\Psi_{\sigma}  \ar[r]^{\sim} & \gr_I^{d_\sigma-1}\gr^K_0 R\Psi_{\sigma} \ar[u]_{can_-} }
$$
o\`u les applications $can$ sont les projections et injections canoniques (notons qu'on a aussi $\gr_I^0\gr^K_{d_\sigma-1} = \gr^M_{d_\sigma-1}$ et $\gr_I^{d_\sigma-1}\gr^K_0=\gr^M_{1-d_\sigma}$). Nous devons donc \'etudier les morphismes $ix_{eq}^*(can_+)$ et $ix^*_{eq}(can_-)$ dans la cat\'egorie $D^b_{\o\QM_l}(G/\varpi^\ZM)$. %Remarquons que l'amplitude cohomologique  de $x^*_{eq}(R\Psi_{\sigma'(\rho)})$ est $[1-d_\rho,0]$.
Pour cela, remarquons tout d'abord que d'apr\`es \ref{bigrad} i), la cohomologie des complexes $ix^*_{eq}(\gr_I^q\gr^K_p R\Psi_{\sigma})$ est concentr\'ee en un seul degr\'e $p+q-1-d_\sigma$. Ce miracle d\'etermine compl\`etement ces complexes dans la cat\'egorie d\'eriv\'ee $D^b_{\o\QM_l}(G/\varpi^\ZM)$, et l'\'enonc\'e de \ref{bigrad} est encore vrai si l'on y remplace $ix^*$ par $ix^*_{eq}$. En voici une cons\'equence :

\begin{lemme} \label{bifilt}  %Soit $\rho\in\Irr{\o\QM_l}{\dd}$ et $i_\sigma^{p,q}$ l'entier d\'efini en \ref{bigrad}. Nous noter
Notons $V'_\sigma:=V_{\sigma}\otimes V_\rho$. % et reprenons la notation $i_\sigma^{p,q}$ de \ref{bigrad}.
\begin{enumerate}
\item Le complexe $x^*_{eq}(I^q\gr^K_p R\Psi_{\sigma})$ est non-nul seulement si $p,q\geq 0$ et $p+q<d_\sigma$, auquel cas il est cohomologiquement concentr\'e en degr\'e $p+q-d_\sigma+1$ et donn\'e par :
$$ ix^*_{eq}(I^q\gr^K_p R\Psi_{\sigma}) \mathop{\simeq}\limits_{} %\left\{\begin{array}{ll} %\pi_{\rho^\vee}^{\{p+q,\cdots,d_\rho-1\}}
\left(\bigoplus\limits_{\rho'\cor\sigma} \pi_{{\rho'}^\vee}^{>p+q} \otimes V'_\sigma\right)[d_\sigma-1-p-q] %& \hbox{si }  p,q, i_\sigma^{p,q}\geq 0   \\
          %  0 & \hbox{sinon} \end{array} \right.
$$
et le morphisme canonique $ix^*_{eq}(I^q\gr^K_p R\Psi_{\sigma})\To{} ix^*_{eq}(\gr_I^q\gr^K_p R\Psi_{\sigma})$ est induit (modulo isomorphismes) sur chaque facteur par l'injection $\tau_{{\rho'}^\vee}^{>p+q}\injo \pi_{{\rho'}^\vee}^{>p+q}$.

\item La ``filtration" par les noyaux $(ix^*_{eq}(K_\bullet R\Psi_{\sigma}))$ est {\em canonique} au d\'ecalage de $1-d_\sigma$ pr\`es, ce qui signifie que pour tout $p$, on a un triangle commutatif
$$ \xymatrix{ix^*_{eq}(K_pR\Psi_{\sigma}) \ar[r] \ar[d]^\sim & ix^*_{eq}(R\Psi_{\sigma}) \\ \tau_{_{< p+1-d_\sigma}} (ix^*_{eq}(R\Psi_{\sigma})) \ar[ru]^{tronq} & }
$$
\end{enumerate}
\end{lemme}
\begin{proof}
On d\'emontre le point i) en fixant $p$ et en faisant une r\'ecurrence descendante de $q=d_\sigma-1-p$ \`a $0$. En effet, le cas $q=d_\sigma-1-p$ est donn\'e par l'\'egalit\'e $I^{d_\sigma-1-p}\gr^K_p = \gr_I^{d_\sigma-1-p}\gr^K_p$ et le point i) de \ref{bigrad}. Supposons la propri\'et\'e d\'emontr\'ee pour $q+1$, et consid\'erons le triangle distingu\'e 
\ini\begin{equation}\label{tri}
 \To{} I^q\gr^K_p [-i_\sigma^{p,q}] \To{} \gr_I^q\gr_p^K[-i_\sigma^{p,q}] \To{} I^{q+1}\gr^K_p[1-i_\sigma^{p,q}]. 
\end{equation}
Par l'hypoth\`ese de r\'ecurrence et par \ref{bigrad} i), les deux derniers objets sont dans le coeur $\Mo{\o\QM_l}{G/\varpi^\ZM}$. Par l'hypoth\`ese de r\'ecurrence et \ref{bigrad} ii), la fl\`eche entre eux  est induite (modulo isomorphismes) sur chaque facteur par l'injection $\pi_{\rho'}^{i_\sigma^{p,q+1}} \injo \tau_{\rho'}^{i_\sigma^{p,q}}$. En particulier cette fl\`eche est  surjective, donc le premier membre $I^q\gr^K_p[-i_\sigma^{p,q}]$ du triangle \ref{tri} est lui aussi dans le coeur $\Mo{\o\QM_l}{G^0}$ et s'identifie au noyau de la deuxi\`eme fl\`eche, ce qui ach\`eve le pas de r\'ecurrence.

Pour le deuxi\`eme point on a %Remarquons tout d'abord que pour tout $p$, la fl\`eche $I^0\gr_p^K \To{} \gr_p^K$ est un isomorphisme.
$\gr^K_p=I^0\gr^K_p$, donc le complexe 
% Or d'apr\`es le th\'eor\`eme \ref{bifilt}, la cohomologie % nous d\'ecrit 
%$x^*(I^0\gr^K_{p}) \in \gr\Mo{\o\QM_l}{G^0}$ de 
%, c'est-\`a-dire qu'il ne nous donne {\em a priori} que la cohomologie de 
%l'objet 
$x^*_{eq}(\gr^K_{p})$ est concentr\'ee en degr\'e $1-d_\sigma+p$, par le point i). %Mais le miracle est que cette cohomologie est concentr\'ee en un seul degr\'e, ce qui nous permet de remonter l'information \`a $D^b_{\o\QM_l}(G^0)$. Comme ce degr\'e est $d_\sigma-1$, 
Le point ii) en d\'ecoule par une r\'ecurrence imm\'ediate.
\end{proof}

Comme les fl\`eches canoniques $I^0\gr_{d_\sigma-1}^K \To{} \gr^0_I\gr^K_{d_\sigma-1}$  et $\HC^0 \To{} \tau_{\geq 0}$ sont des isomorphismes (pour $x_{eq}^*(R\Psi_{\sigma})$), le point ii) du lemme fournit  un triangle commutatif
$$ \xymatrix{ ix_{eq}^*(R\Psi_{\sigma}) \ar[rd]^{tronq_+} \ar[d]_{ix_{eq}^*(can_+)} &  \\ ix_{eq}^*(\gr_I^0\gr^K_{d_\sigma-1})   & \HC^{0}(ix^*_{eq}(R\Psi_{\sigma})) \ar[l]_{\sim} }
$$
De m\^eme, la fl\`eche $\gr_I^{d_\sigma-1} \gr^K_0 \To{} \gr^K_0$ permet de factoriser $ix_{eq}^*(can_-)$
$$ \xymatrix{ ix_{eq}^*(R\Psi_{\sigma})  & \HC^{1-d_\sigma}(ix^*_{eq}(R\Psi_{\sigma}))[d_\sigma-1] \ar[l]_-{tronq_-}  \\ ix_{eq}^*(\gr_I^{d_\sigma-1}\gr^K_{0}) \ar[u]^{ix^*_{eq}(can_-)}  \ar[r]^\alpha &  ix_{eq}^*(\gr^K_{0}) \ar[u]_\sim }
$$
o\`u la fl\`eche $\alpha$ s'identifie \`a la compos\'ee
$$\xymatrix{ ix_{eq}^*(I^{d_\sigma-1}\gr^K_0) \ar[r] \ar@{=}[d]&   ix_{eq}^*(I^{d_\sigma-2}\gr^K_0) \ar[r] \ar@{=}[d] & \cdots \ar[r] & ix_{eq}^*(I^{0}\gr^K_0) \ar@{=}[d] \\ \left(\bigoplus\limits_{\rho'\cor\sigma} \pi_{{\rho'}^\vee}^{>d_\sigma-1}\right) \otimes V'_{\sigma}[0] \ar[r]^{\alpha_{d_\sigma-1}} &  
\left( \bigoplus\limits_{\rho'\cor\sigma} \pi_{{\rho'}^\vee}^{>d_\sigma-2}\right) \otimes V'_{\sigma}[1]  \ar[r]^-{\alpha_{d_\sigma-2}} &
\cdots \ar[r]^-{\alpha_{1}} & \left(\bigoplus\limits_{\rho'\cor\sigma} \pi_{{\rho'}^\vee}^{>0}\right)\otimes V'_{\sigma}[d_\sigma-1]
}, $$
gr\^ace \`a \ref{bifilt} i). Comme le complexe $ix^*_{eq}(R\Psi_\sigma)$ est scind\'e dans $D^b_{\o\QM_l}(G/\varpi^\ZM)$, le morphisme de troncation $tronq_+$ est un \'epimorphisme tandis que $tronq_-$ est un monomorphisme. Il ne reste donc  plus qu'\`a prouver que la compos\'ee ci-dessus est non-nulle. %Que l'on puisse d\'ecrire cette compos\'ee dans $D^b_{\o\QM_l}(G^0)$ rel\`eve du miracle.
Or, pour tout $0<q<d_\sigma$, on a des triangles distingu\'es 
$$\xymatrix{ ix_{eq}^*(I^{q+1}\gr^K_0)[1-d_\sigma] \ar[r] \ar@{=}[d]&   ix_{eq}^*(I^{q}\gr^K_0)[1-d_\sigma] \ar[r] \ar@{=}[d] &  ix_{eq}^*(\gr_I^q\gr^K_0)[1-d_\sigma] \ar@{=}[d] \ar[r] & \\ 
\left(\bigoplus\limits_{\rho'\cor\sigma} \pi_{{\rho'}^\vee}^{>q+1}\right)\otimes V'_{\sigma}[-q-1] \ar[r]^-{\alpha_{q+1}} &  
\left(\bigoplus\limits_{\rho'\cor\sigma} \pi_{{\rho'}^\vee}^{>q}\right) \otimes V'_{\sigma}[-q]  \ar[r] &
 \left(\bigoplus\limits_{\rho'\cor\sigma} \tau_{{\rho'}^\vee}^{>q}\right)\otimes V'_{\sigma}[-q]\ar[r] &
} $$
qui montrent d'apr\`es \ref{bifilt} i) que $\alpha_{q+1}$ est donn\'e sur chaque facteur, et modulo isomorphismes, par un \'el\'ement non-nul de
%$I^{q+1}\To{} I^q \To{} \gr_I^q$ fournissent apr\`es rotation et d\'ecalage, des suites exactes courtes dans $\Mo{\o\QM_l}{G^0\times I_K\times \OC_D^\times}$
%$$
%0\To{} x^*(I^q \gr^K_p )[-i_\sigma^{0,q}] \To{} x^*(\gr_I^q \gr^K_p)[-i_\sigma^{0,q}] \To{} x^*(I^{q+1} \gr^K_p) [-i_\sigma^{0,q}+1] \To{} 0.
%$$ 
%dont la classe provient du g\'en\'erateur de  %$\ext{1}{\pi^{\{p+q+1,\cdots,d_\sigma-1\}}_{\rho^\vee}}{\pi^{\{p+q,\cdots,-1\}}_{\rho^\vee}}{G^0}, $
$\ext{1}{\pi_{{\rho'}^\vee}^{>q+1}}{\pi_{{\rho'}^\vee}^{>q}}{G/\varpi^\ZM}. $
Par \ref{theoextell} ii) et \ref{remaextell}, il s'ensuit que la compos\'ee $\alpha$ est  non-nulle dans $D^b_{\o\QM_l}(G/\varpi^\ZM)$.

\alin{Preuve de la proposition \ref{equiv}} \def\top{\hbox{\sl Top}}
Notons $\wt{\NM(G)}$ le topos des pr\'efaisceaux sur $\NM(G)$. On a un morphisme de topos $(\delta^*,\delta_*)$ du topos classifiant $ \wt{G}$ de $G$ (dont les objets sont les ensembles munis d'une action lisse de $G$) vers $ \wt{\NM(G)}$, donn\'e par 
$$ \delta_*(E)=\left(n\mapsto E^{H_n}\right) \;\; \hbox{ et } \;\; \delta^*(F):= \limi{n\in\NM} F(n). $$
%Posons $\Lambda:=\ZM_l$, on a donc un foncteur $$R\delta^*_\bullet := \delta^*\circ R\limproj^{\NM(G)} :\; D^+_{\Lambda_\bullet}(\NM(G)) \To{} D^+_\Lambda(G)$$
%qui apr\`es tensorisation induit un foncteur encore not\'e $R\delta_\bullet^* : D^+_{\Lambda_\bullet}(\NM(G))\otimes \o\QM_l \To{} D^+_{\o\QM_l}(G)$.

%Notons que le foncteur $\delta^*$ est aussi exact \`a droite sur les objets ab\'eliens.
Par ailleurs notons $\wt{S^{s,ca}_{HT,et}} \To{} \NM(G)$ le topos \'etale fibr\'e associ\'e au $\NM(G)$-diagramme de $k^{ca}$-sch\'emas form\'e par les $(S^{s,ca}_{HT,n})_{n\in\NM}$, et %, selon la notation consacr\'ee, 
$\top(\wt{S^{s,ca}_{HT,et}})$ son topos total. Si $\Lambda$ est une extension de $\ZM_l$,
% on se gardera de confondre la cat\'egorie d\'eriv\'ee 
on a un foncteur ``oubli" $\omega$ de la cat\'egorie d\'eriv\'ee $D^{b,top}_{\Lambda_\bullet}(S_{HT}):=D^b_{\Lambda_\bullet}(\top(\wt{S^{s,ca}_{HT,et}}))$ des 
$\Lambda_\bullet$-modules de $\top(\wt{S^{s,ca}_{HT,et}})$  
vers la cat\'egorie $D^{b,naif}_{\Lambda_\bullet}(S_{HT})$ des sections de la cat\'egorie bi-fibr\'ee sur $\NM(G)$ dont les fibres sont les cat\'egories d\'eriv\'ees $D^b_{\Lambda_\bullet}(S^{s,ca}_{HT,n})$. %Cette derni\`ere est gradu\'ee et munie d'une famille naturelle de triangles distingu\'es, mais {\em n'est pas triangul\'ee} (aucun des axiomes $TR_i$ n'est a priori v\'erifi\'e). 
%on a un foncteur \'evident
%$\omega:\; D^+_{\Lambda_\bullet}(\SC_{HT}) \To{} D^{+, naif}_{\Lambda_\bullet}(S_{HT})$ qui envoie triangles distingu\'es sur triangles distingu\'es. 
Remarquons que, par $t$-exactitude des $g^{m|n}_*$, la cat\'egorie $FPH_{\Lambda_\bullet}^c$ introduite au-dessus de \ref{subtilite} est une sous-cat\'egorie pleine de $D^{b,naif}_{\Lambda_\bullet}(S_{HT})$. On notera  $\FC\PC\HC_{\Lambda_\bullet}^c$ la sous-cat\'egorie pleine de $D^{b,top}_{\Lambda_\bullet}(S_{HT})$ des objets dont l'image par $\omega$ est dans $FPH_{\Lambda_\bullet}^c$.
%Il s'ensuit que la %cat\'egorie $D^{b,naif}_c(S_{HT},\o\QM_l)$ de la preuve du lemme \ref{} et par l\`a la 
%cat\'egorie $FPH$ est une sous-cat\'egorie pleine de $D^{+, naif}_{\Lambda_\bullet}(S_{HT})\otimes \o\QM_l$, dont les suites exactes courtes se prolongent uniquement en triangles distingu\'es.
%Notons alors $\FC\PC\HC$ la sous-cat\'egorie pleine de $D^{b}_{\Lambda_\bullet}(\SC_{HT})\otimes \o\QM_l$ dont les objets sont ceux dont l'image par $\omega$ est (essentiellement) dans $FPH$. 
D'apr\`es l'appendice, proposition \ref{fptopfib}, le foncteur $\omega$ induit une \'equivalence de cat\'egories $\FC\PC\HC_{\Lambda_\bullet}^c \simto FPH_{\Lambda_\bullet}^c$, ce qui fait de $\FC\PC\HC_{\Lambda_\bullet}^c$ une sous-cat\'egorie ab\'elienne admissible de $D^{b,top}_{\Lambda_\bullet}(S_{HT})$.

Le point supersingulier $x$ induit un morphisme de topos $(\top(x)^*,\top(x)_*) :\;\; \wt{\NM(G^0)} \To{} \top(\wt{S^{s,ca}_{HT,et}})$. On a en particulier un foncteur 
$$ \top(x)^*_\bullet  :\; D^{b,top}_{\Lambda_\bullet}(S_{HT}) \To{\top(x)^*} D^b_{\Lambda_\bullet}(\wt{\NM(G^0)}) \To{R\limproj^{\NM(G^0)}} D^b_{\Lambda}(\wt{\NM(G^0)}) $$ 
et une famille de foncteurs
$$ x_{n,_\bullet}^*  :\; D^b_{\Lambda_\bullet}(\wt{S_{HT,n}^{s,ca}}) \To{x_n^*} D^b_{\Lambda_\bullet}(pt) \To{R\limproj} D^b_{\Lambda}(pt) \To{\HC^*} \gr\Mo{}{\Lambda} $$
qui s'inscrivent
dans le diagramme essentiellement commutatif suivant  :
$$\xymatrix{ \FC\PC\HC_{\Lambda_\bullet}^c \ar[d]_{\simeq} \ar@{^(->}[r] & D^{b,top}_{\Lambda_\bullet}(S_{HT}) \ar[d]_{\omega}\ar[r]^{\top(x)^*_\bullet} & D^b_{\Lambda}(\wt{\NM(G^0)})\ar[d] \ar[r]^{\delta^*}  &  D^b_{\Lambda}(G^0) \ar[d]^{\HC^*} \\
FPH_{\Lambda_\bullet}^c \ar@{^(->}[r] & D^{b,naif}_{\Lambda_\bullet}(S_{HT}) \ar[r]^{(x_{n,\bullet}^*)_n} &  \gr\Mo{\Lambda}{\wt{\NM(G^0)}} \ar[r]^{\delta^*} & \gr\Mo{\Lambda}{G^0}
}$$
Dans ce diagramme, la compos\'ee de la ligne du bas induit, apr\`es inversion de $l$ et passage \`a la limite sur $\Lambda$, le foncteur not\'e plus haut $x^*\circ \iota : FPH' \To{} \gr\Mo{\o\QM_l}{G^0}$. Nous poserons donc $x_{eq}^*$ le foncteur obtenu en composant la ligne du haut avec un quasi-inverse de $\omega_{|\FC\PC\HC_{\Lambda^\bullet}^c}$, en inversant $l$ et en passant \`a la limite sur $\Lambda$.

Il ne nous reste plus qu'\`a exhiber un isomorphisme canonique $x^*_{eq}(R\Psi)\simto R\Gamma(\MC_{LT}^{(0),ca},\o\QM_l)$ dans $D^b_{\o\QM_l}(G^0)$.
On va pour cela se ramener \`a l'isomorphisme de Berkovich entre cycles \'evanescents formels et alg\'ebriques. En effet, pour tout $n$, on a un diagramme de cat\'egories
$$\xymatrix{ \wt{S^{\eta,nr}_{HT,n,et}} \ar[d]_{\psi_{*n}} \ar[r]^{x_n^{an,*}} & \wt{\MC^{(0)}_{LT,n,et}}  \ar[d]^{\gamma_{*n}} \\ \wt{S^{s,ca}_{HT,n,et}} \ar@{=>}[ru] \ar[r]^{x_n^*} & \wt\bullet }$$ 
o\`u $\psi_n= i_n^{s,ca,*}\circ j_{n,*}^{\eta,ca}$ est le foncteur cycles \'evanescents ``usuel'' pour les vari\'et\'es alg\'ebriques d\'efinies sur $\wh\OC_K^{nr}$, $x_n^{an,*}$ envoie un faisceau \'etale sur la restriction de son analytifi\'e \`a $\MC_{LT,n}^{(0)}$ (foncteur not\'e $\FC\mapsto \wh\FC$ chez Berkovich \cite[5.1]{Bic3} \cite[3.1]{Bic4}), et $\gamma_n$ est le foncteur des sections globales apr\`es extension des scalaires de $\knr$ \`a $\ka$, qui s'identifie avec le foncteur not\'e $\Psi_\eta$ par Berkovich \cite[p. 373]{Bic4}  puisque la fibre sp\'eciale r\'eduite de $\MC_{LT,n}^{(0)}$ est un point. La double fl\`eche diagonale signifie qu'on a une transformation naturelle canonique $x_n^*\circ \psi_n \To{} \gamma_n\circ x_n^{an,*}$. D'apr\`es \cite[Thm 3.1]{Bic4}, cette transformation induit un isomorphisme $x_n^* R\psi_{n*} (\Lambda_\bullet) \simto R\gamma_{n*}(\Lambda_\bullet)$ dans $D^b(\Lambda_\bullet)$.

Ces cat\'egories, foncteurs, et transformations naturelles, s'organisent au-dessus de $\NM(G^0)$ et fournissent un diagramme de cat\'egories
$$\xymatrix{ \hbox{Top}(\wt{S^{\eta,nr}_{HT,et}}) \ar[d]_{\psi_*} \ar[r]^{x^{an,*}} & \hbox{Top}(\wt{\MC_{LT,et}^{(0)}})  \ar[d]^{\gamma_*} \\ \hbox{Top}(\wt{S^{s,ca}_{HT,et}}) \ar@{=>}[ru] \ar[r]^{\hbox{Top}(x)^*} & \wt{\NM(G^0)} }$$ 
et en particulier un morphisme 
\ini\begin{equation}\label{is1}
\hbox{Top}(x)^* R\psi_*(\Lambda_\bullet) \To{} R\gamma_*(\Lambda_\bullet)  \hbox{ dans } D^+_{\Lambda_\bullet}(\wt{\NM(G^0)}).
\end{equation}
 Ce dernier est un isomorphisme, puisqu'il l'est en restriction \`a chaque \'etage. 

Par d\'efinition,  l'objet  $R\Psi$ de $FPH'$ est donn\'e par le syst\`eme inductif 
$((R\psi_{n*}(\Lambda_\bullet))_{n\in\NM(G)}[\frac{1}{l}])_\Lambda$. Mais comme les morphismes de transitions $(g^{m|n})^\eta$ du $\NM(G^0)$-diagramme $(S^{\eta}_{HT,n})_n$ sont {\em \'etales} (fibre g\'en\'erique), la restriction de $R\psi_*(\Lambda_\bullet)$ \`a l'\'etage $n$ est canoniquement isomorphe \`a $R\psi_{n*}(\Lambda_\bullet)$. En d'autres termes $\omega^{-1}(R\Psi)$ est donn\'e par le syst\`eme inductif
 $((R\psi_*(\Lambda_\bullet)[\frac{1}{l}])_{\Lambda})$ de la cat\'egorie $\limi{}(\FC\PC\HC_{\Lambda_\bullet}^c[\frac{1}{l}])$. 
Vu l'isomorphisme \ref{is1}, il nous suffira donc de prouver que pour tout $\Lambda$, on a un isomorphisme
\ini\begin{equation}\label{is2}
\delta^* R\limproj^{\NM(G^0)} R\gamma_*(\Lambda_\bullet) \simeq R\Gamma(\MC_{LT}^{(0),ca},\Lambda)
\end{equation}
 dans $D^b_{\Lambda}(G^0)$.
On remarque pour cela que le morphisme de p\'eriodes $\xi_{LT} : \MC^{(0)}_{LT,0} \To{} \PC_{LT}$ induit un foncteur exact et pleinement fid\`ele
$$ \xi^* :\; \wt{\PC_{LT,et}} \To{} \hbox{Top}(\wt{\MC_{LT,et}^{(0)}}) $$
dont l'image est constitu\'ee des faisceaux {\em cart\'esiens} du topos total de droite, et que, revenant aux d\'efinitions, on a  les factorisations suivantes :
$$ \Gamma(\MC_{LT}^{(0),ca},-)= \delta^* \circ \gamma_* \circ \xi^* :\;\;  \wt{\PC_{LT,et}} \To{} \wt{G^0} $$
et
$$ \Gamma_\bullet(\MC_{LT}^{(0),ca},-) = \delta^*\circ \limproj^{\NM(G^0)} \circ \gamma_*\circ \xi^* :\; \Mo{\Lambda_\bullet}{\wt{\PC_{LT,et}}} \To{} \Mo{\Lambda}{G^0}.$$
(Remarquer que $\gamma_*$ commute aux limites projectives quelconques, ce qui nous permet de placer la limite projective l\`a o\`u on l'a plac\'ee).
Mais $\xi_{LT}$ est un morphisme \'etale, donc $\xi^*$ envoie injectifs sur injectifs, de sorte que le dernier isomorphisme se d\'erive en $R\Gamma_\bullet(\MC_{LT}^{(0),ca},-) = \delta^*\circ R\limproj^{\NM(G^0)} \circ R\gamma_*\circ \xi^*$, ce qui conclut la preuve de \ref{is2} puisque $\xi^*(\Lambda_\bullet)=\Lambda_\bullet$.

\subsection{Uniformisation $p$-adique et conjecture de puret{\'e}} \label{unifpurete}

Dans cette partie nous prouvons que les vari\'et\'es uniformis\'ees par les rev\^etements  de l'espace sym\'etrique de Drinfeld satisfont la conjecture dite ``monodromie-poids" de Deligne. Cela s'applique \`a certaines vari\'et\'es de Shimura associ\'ees \`a des groupes unitaires d\'efinis globalement par une alg\`ebre \`a division et involution satisfaisant les m\^emes propri\'et\'es que celles intervenant dans les vari\'et\'es de Harris-Taylor, except\'e qu'\`a la place o\`u le compl\'et\'e est $K$, elle doit \^etre totalement ramifi\'ee, au lieu d'\^etre d\'eploy\'ee. Pour de telles vari\'et\'es de Shimura, on en d\'eduit donc le facteur local de la fonction $L$ en la place concern\'ee.

Cependant, nous  \'eviterons d'introduire les nombreuses notations n\'ecessaires \`a la d\'efinitions de ces vari\'et\'es de Shimura, et nous contenterons d'exposer les raisonnements pour des vari\'et\'es ``abstraites" uniformis\'ees par les $\MC_{Dr,n}$. On commence par rappeler l'\'enonc\'e de la conjecture de Deligne.

\alin{Puret\'e de la filtration de monodromie} \label{monpds}
 On sait depuis Grothendieck que sur toute repr\'esentation   $l$-adique continue de dimension finie $(\sigma,V)$ de $W_K$, l'inertie $I_K$ agit de mani\`ere quasi-unipotente, 
c'est-\`a-dire qu'il existe 
un unique endomorphisme nilpotent $N_\sigma$ de $V$ tel qu'il existe un sous-groupe $I_\sigma \subset I_K$ d'indice fini agissant via la formule $ \forall i\in I_\sigma,\;\; \sigma(i)=\hbox{exp}(N_\sigma t_\mu(i))$, o\`u $t_\mu$ a \'et\'e introduit au paragraphe \ref{formcontl}.
L'op\'erateur $N_\sigma$ est donc le logarithme de la partie
unipotente de la monodromie de $\sigma$, mais par abus de langage,
nous l'appelerons simplement ``op\'erateur de monodromie de
$\sigma$"\footnote{Bien-s\^ur $N_\sigma$ d\'epend du g\'en\'erateur
  $\mu$ de $\ZM_l(1)$, mais $N_\sigma\otimes \mu^* : \;
  V\otimes_{\ZM_l} \ZM_l(1) \To{} V$ n'en d\'epend pas.}. 
Il v\'erifie n\'ecessairement l'\'equation habituelle
$wN_\sigma w^{-1}=|w|N_\sigma$ de sorte que pour tout rel\`evement de
Frobenius g\'eom\'etrique $\phi$, 
l'application
\ini\begin{equation}
\label{sigmalisse}
w\mapsto  \sigma^\phi(w):=\sigma(w)\hbox{exp}(-N_\sigma t_\mu(i_\phi(w)))\;\;\hbox{ o\`u } \; w=\phi^{\nu(w)}i_\phi(w) 
 \end{equation}
d\'efinit une repr\'esentation {\em lisse} de $W_K$ sur $V$. Rappelons aussi que par
\cite[8.4.2]{DelAntwerp}, la classe d'isomorphisme de $\sigma^\phi$ ne d\'epend pas du choix de $\phi$.

Ceci s'applique en particulier aux espaces de cohomologie $l$-adique $H^i(X_{K^{ca}},\o\QM_l)$ d'une vari\'et\'e $X$ propre sur $K$.
%dans un syst\`eme local $l$-adique $\LC$.
\`A l'op\'erateur de monodromie $N$ est associ\'ee une filtration croissante stable sous $W_K$ de l'espace $V$ dite ``filtration de monodromie"\footnote{qui, elle, est bien ind\'ependante du choix de $\mu$} $\cdots \subseteq M_iV\subseteq M_{i+1}V \subseteq\cdots $ de longueur finie et caract\'eris\'ee par les propri\'et\'es que $N(M_iV)\subseteq M_{i-2}V$ pour tout $i$, et $N$ induit des isomorphismes $N^i:\;\; \gr^M_{i}V \otimes |-|^i\simto \gr^{M}_{-i}V$ pour tout $i\geq 0$.
Par ailleurs, Deligne a prouv\'e l'existence d'une autre filtration croissante stable sous $W_K$ de $V$, dite ``filtration par les poids", $\cdots \subseteq W_iV \subseteq W_{i+1}V\subseteq \cdots$ caract\'eris\'ee par la propri\'et\'e que les valeurs propres de {\em tout} rel\`evement  de Frobenius g\'eom\'etrique sur $\gr^W_{i}V$ sont des entiers alg\'ebriques dont tous les conjugu\'es complexes sont de norme complexe $q^{i/2}$.

\begin{conj} \label{conmp}(Monodromie-Poids) Si $X$ est propre et lisse sur $K$, 
%et $\LC$ est un syst\`eme local $l$-adique, alors les gradu\'es de la filtration de monodromie sont purs de poids ...
alors pour tout $i\in \NM$, on a
$$ M_i\left(H^j(X^{ca},\o\QM_l)\right)= W_{i+j}\left(H^j(X^{ca},\o\QM_l)\right) .$$
\end{conj}

Lorsque $K$ est d'\'egales caract\'eristiques, l'\'enonc\'e est essentiellement contenu et d\'emontr\'e dans les travaux de Deligne sur les conjectures de Weil. Le cas d'in\'egales caract\'eristiques est tr\`es peu avanc\'e, m\^eme dans les cas de r\'eduction semi-stable.

Remarquons que pour tout $i\in\ZM$ on a $N(W_iV)\subseteq W_{i-2}V$, de sorte que par la caract\'erisation de la filtration de monodromie, la conjecture ci-dessus est \'equivalente \`a l'assertion : {\em Pour tout $i\geq 0$, $N$ induit un isomorphisme  
\ini\begin{equation}\label{equiMP}
N^i:\;\; \gr^W_{i+j}\left(H^j(X^{ca},\o\QM_l)\right)\otimes |-|^i \simto \gr^{W}_{-i+j}\left(H^j(X^{ca},\o\QM_l)\right).
\end{equation}
}

\alin{Uniformisation $p$-adique}
Soit $\Gamma$ un sous-groupe {\em discret, cocompact et sans torsion} de $G$. 
On sait que l'action d'un tel sous-groupe  sur le $\knr$-espace analytique $\MC_{Dr,0}$ est {\em libre}, et il en est donc de m\^eme de l'action sur les rev\^etements $\MC_{Dr,n}$.
Par \cite[lemma 4]{Bic5}, l'espace annel\'e quotient $\MC_{Dr,n}/\Gamma$ est  muni d'une structure de $\knr$-espace analytique.
Notons $\Gamma_Z:=\Gamma\cap K^\times$ ; c'est un sous-groupe discret cocompact de $K^\times$. D'apr\`es \cite[Thm 3.49]{RZ}, la donn\'ee de descente \`a la Weil de $\MC_{Dr,n}$ devient effective sur $\MC_{Dr,n}/\Gamma_Z$, donc {\em a fortiori} sur $\MC_{Dr,n}/\Gamma$. %Nous noterons $S_{\Gamma,n}$ l'espace analytique descendu \`a $K$. 

\begin{fact} (Mustafin, Cherednik, Drinfeld, Rapoport-Zink, Varshawski) 
Le $K$-espace analytique obtenu par descente du quotient $\MC_{Dr,n}/\Gamma$ est alg\'ebrisable. Plus pr\'ecis\'ement, il existe une vari\'et\'e alg\'ebrique $S_{\Gamma,n}$ propre et lisse sur $K$ dont l'analytification lui est canoniquement isomorphe. En particulier, $S_{\Gamma,n}$ est munie d'une action de $\dd$ et d'apr\`es  le th\'eor\`eme de comparaison  GAGA de Berkovich \cite[7.1]{Bic2}, il y a des isomorphismes $\dd\times W_K$-\'equivariants
$$ H^i(S_{\Gamma,n}\otimes_K K^{ca},\o\QM_l) \simto H^i_c(\MC_{Dr,n}^{ca}/\Gamma,\o\QM_l) $$
o\`u le terme de gauche d\'esigne la cohomologie \'etale $l$-adique au sens des vari\'et\'es alg\'ebriques.

\end{fact}

%\begin{proof}
%Dans le cas o\`u $K$ est $p$-adique, la r\'ef\'erence la plus commode est le th\'eor\`eme 6.36 de \cite{RZ}. L'uniformisation y est obtenue par voie modulaire et contient beaucoup plus d'informations que ce dont nous avons besoin ici. 
%Dans le cas d'\'egales caract\'eristiques, le r\'esultat dans la g\'en\'eralit\'e demand\'ee est d\^u \`a Drinfeld, {\em cf} les commentaire de la partie III.5 de \cite{BouCar}.
%\end{proof}

Ce r\'esultat  permet de s'assurer que la cohomologie $l$-adique de $\MC^{ca}_{Dr,n}/\Gamma$ est munie d'une filtration par les poids et d'une filtration de monodromie, mais on aurait aussi pu le voir directement, \`a partir de l'isomorphisme $W_K$-\'equivariant de Hochschild-Serre \cite[Prop B.3.1]{Dat1}
$$\o\QM_l \otimes^L_{\o\QM_l[\Gamma]} R\Gamma_c(\MC_{Dr,n}^{ca},\o\QM_l)
\simto R\Gamma_c(\MC_{Dr,n}^{ca}/\Gamma,\o\QM_l) $$
o\`u, dans l'expression de gauche, le complexe est vu \`a travers le foncteur d'oubli $D^b_{\o\QM_l}(GD\times W_K^{disc})\To{} D^b_{\o\QM_l}((\Gamma\times \dd)/\Gamma_Z^{diag}\times W_K^{disc})$.

Dor\'enavant, nous surlignerons tous les quotients par $\Gamma_Z$ (par exemple ceux de $\Gamma$, $\dd$, $G$, ou $\MC_{Dr,n}$), et pour une repr\'esentation $\rho\in\Irr{\o\QM_l}{\o\dd}$, nous noterons $(-)[\o\rho]$ le foncteur exact $(-)\otimes_{\o\QM_l\o\dd} \rho$ \`a ne pas confondre avec  le foncteur $(-)[\rho]=(-)\otimes^L_{\o\QM_l\dd}\rho$ introduit au d\'ebut de la section \ref{real}. 
% En particulier, si $\rho\in \Irr{\o\QM_l}{\dd/\Gamma_Z(1+\varpi^n\OC_D)}$, on a des isomorphismes $W_K$-\'equivariants
%$$ H^i(S_{\Gamma,n}^{ca},\o\QM_l)[\rho]\simto H^i_c(\MC^{ca}_{Dr,n}/\Gamma,\o\QM_l)[\rho],$$
%la notation $[\rho]$ \'etant celle de \ref{}. 
Avec ces notations on a une d\'ecomposition
$$ R\Gamma_c(\MC_{Dr,n}^{ca}/\Gamma,\o\QM_l) \simeq \bigoplus_{\rho\in\Irr{}{\o\dd/(1+\varpi^n\OC_D)}} R\Gamma_c(\MC_{Dr,n}^{ca}/\Gamma,\o\QM_l)[\o\rho]\otimes \rho^\vee $$
qui induit en cohomologie des d\'ecompositions
$$ H^i(S_{\Gamma,n}^{ca},\o\QM_l) \simeq \bigoplus_{\rho\in\Irr{}{\o\dd/(1+\varpi^n\OC_D)}}  H^i(S_{\Gamma,n}^{ca},\o\QM_l)[\o\rho]\otimes \rho^\vee, $$
o\`u l'action de $W_K$, resp. $\dd$, se fait sur les premiers, resp. seconds, facteurs des produits tensoriels.
L'isomorphisme de Hochschild-Serre se d\'ecompose aussi et donne :
\ini\begin{equation}\label{isoHS}
 \o\QM_l \otimes^L_{\o\QM_l[\o\Gamma]} (R\Gamma_c[\rho])
\simto R\Gamma_c(\MC_{Dr,n}^{ca}/\Gamma,\o\QM_l)[\o\rho] 
\end{equation}
pour toute $\rho\in\Irr{\o\QM_l}{\dd/\Gamma_Z(1+\varpi^n\OC_D)}$.
Soit alors $N_{\rho,\Gamma,j}$ l'endomorphisme du $\o\QM_l$-espace vectoriel $H^{j}(S_{\Gamma,n}^{ca},\o\QM_l)[\o\rho]$ fonctoriellement induit par $N_\rho$ via  l'isomorphisme pr\'ec\'edent, le passage \`a la cohomologie en degr\'e $j$ et le th\'eor\`eme de comparaison GAGA. Comme tous les isomorphismes utilis\'es sont $W_K$-\'equivariants, la d\'efinition de $N_\rho$ montre que le sous-groupe $I_\rho\subset I_K$ du lemme \ref{monqunip} agit sur $H^{j}(S_{\Gamma,n}^{ca},\o\QM_l)[\rho]$ par $i\mapsto \hbox{exp}(N_{\rho,\Gamma,j}t_\mu(i))$, de sorte que $N_{\rho,\Gamma,j}$ est l'op\'erateur de monodromie de la repr\'esentation de $W_K$ sur $H^{j}(S_{\Gamma,n}^{ca},\o\QM_l)[\rho]$.

Fixons dor\'enavant un rel\`evement de Frobenius g\'eom\'etrique $\phi$ et 
notons $H^{j,\phi}(S_{\Gamma,n}^{ca},\o\QM_l)[\o\rho]$ la repr\'esen\-tation lisse de $W_K$ associ\'ee, comme en \ref{sigmalisse}.
En combinant l'isomophisme \ref{isoHS} avec le scindage $\alpha_\phi$ de \ref{lemmeeq}, on constate que la suite spectrale de Hochschild-Serre d\'eg\'en\`ere en des isomorphismes $\dd\times W_K$-\'equivariants %({\em cf} \cite[]{Dat1} pour plus de d\'etails) 
\begin{eqnarray*} 
H^{j,\phi}(S_{\Gamma,n}^{ca},\o\QM_l)[\o\rho] & \simto  & \bigoplus_{i=0}^{d_\rho-1}
\tor{j-d+1-i}{\o\QM_l}{\pi^{\leq i}_\rho}{\o\Gamma}\otimes \sigma'(\rho^\vee)|-|^{-i} \\
& \simeq & \bigoplus_{i=0}^{d_\rho-1} \ext{d-1-j+i}{\pi^{\leq i}_\rho}{\o\QM_l}{\o\Gamma}^* \otimes\sigma'(\rho^\vee)|-|^{-i}\\
& \simeq & \bigoplus_{i=0}^{d_\rho-1}
\ext{d-1-j+i}{\pi^{\leq i}_\rho}{\CC^\infty(G/\Gamma,\o\QM_l)}{\o{G}}^* \otimes \sigma'(\rho^\vee)|-|^{-i}
\end{eqnarray*}
La repr\'esentation $\CC^\infty({G}/\Gamma,\o\QM_l)$ de $G$ est admissible et semi-simple, avec constituents ``unitarisables" : en particulier, les seules repr\'esentations elliptiques qui peuvent y apparaitre sont les s\'eries discr\`etes et les repr\'esentations de Speh locales.
Soit $m_{\rho,\Gamma}^\emptyset$, resp. $m'_{\rho,\Gamma}$, la multiplicit\'e de la s\'erie discr\`ete $\pi^\emptyset_{\rho}$, resp. de la repr\'esentation de Speh $\pi^{S_{\rho}}_{\rho}$  dans cette repr\'esentation $\CC^\infty({G}/\Gamma,\o\QM_l)$. Compte tenu du calcul de $\hbox{Ext}$ entre repr\'esentations elliptiques, on trouve la description suivante :

\noindent{\em Pour $j=d-1$, on a
$$ {H}^{d-1,\phi}(S_{\Gamma,n}^{ca},\o\QM_l)[\o\rho] \mathop{\simeq}\limits_{\dd\times W_K}
 \left\{ \begin{array}{ll}
 \left(\mathop{\bigoplus}\limits_{i=0}^{d_\rho-1}  \sigma'(\rho^\vee)|.|^{-i}\right)^{m_{\rho,\Gamma}^\emptyset} & \hbox{si }
d_\rho \hbox{ est pair} \\
\left(\mathop{\bigoplus}\limits_{i=0}^{d_\rho-1}  \sigma'(\rho^\vee)|.|^{-i}\right)^{m_{\rho,\Gamma}^\emptyset} \oplus \left(\sigma'(\rho^\vee)|.|^{-k}\right)^{m'_{\rho,\Gamma}} & \hbox{si }
d_\rho=1+2k
\end{array}\right.
$$
et pour $j\neq d-1$, on a 
$$ {H}^{j,\phi}(S_{\Gamma,n}^{ca},\o\QM_l)[\o\rho] \mathop{\simeq}\limits_{\dd\times W_K}
 \left\{ \begin{array}{ll}
 0 & \hbox{si }
j+d_\rho-d \hbox{ est impair} \\
\left(\sigma'(\rho^\vee)|.|^{-k}\right)^{m'_{\rho,\Gamma}} & \hbox{si }
j+d_\rho-d=2k\geq 0
\end{array}\right.
$$
}
Observons en particulier que l'action de $W_K$ sur ${H}^{j,\phi}(S_{\Gamma,n}^{ca},\o\QM_l)[\o\rho]$ est semi-simple. On en d\'eduit la remarque suivante :

\begin{rema} 
L'action d'un rel\`evement de Frobenius $\phi$ sur les espaces de cohomologie $H^i(S_{\Gamma,n}^{ca},\o\QM_l)$ est semi-simple.
\end{rema}

Revenons \`a notre probl\`eme initial ; 
par le th\'eor\`eme de ``multiplicit\'es limites" de \cite[1.3]{Rog}, on sait que, quitte \`a remplacer $\Gamma$ par un sous-groupe d'indice fini, on peut supposer $m_{\rho,\Gamma}^\emptyset >0$. 
On peut maintenant \'enoncer
\begin{prop} \label{propequMP}
Avec les notations ci-dessus, supposons $\Gamma$ ``assez petit" pour que $m_{\rho,\Gamma}^\emptyset\neq 0$. Alors les propri\'et\'es suivantes sont \'equivalentes (toujours sous la description \ref{conjHarris}) :
\begin{enumerate}
\item L'endomorphisme nilpotent $N_\rho$ de $R\Gamma_c[\rho]$ d\'efini dans le lemme \ref{monqunip} est d'ordre $d_\rho$ ({\em i.e.} v\'erifie $N_\rho^{d_\rho-1}\neq 0$).
\item L'op\'erateur de monodromie $N_{\rho,\Gamma,d-1}$ de $H^{d-1}(S_{\Gamma,n}^{ca},\o\QM_l)[\o\rho]$ est d'ordre $d_\rho$.
\item La conjecture monodromie-poids est v\'erifi\'ee pour les $H^{j}(S_{\Gamma,n}^{ca},\o\QM_l)[\o\rho]$, $j\in\NM$. 
\end{enumerate}
\end{prop}
\begin{proof}
\'Evidemment, $ii)\Rightarrow i)$, puisque $N_\rho$ est d'ordre au plus $d_\rho$ et induit $N_{\rho,\Gamma,d-1}$. 
Par ailleurs, rappelons  que la repr\'esentation $\sigma'(\rho^\vee)=\sigma_{d/d_\rho}(\tau^0_\rho)|-|^{\frac{d_\rho-d}{2}}$ est pure de poids $d-d_\rho$. 
Ainsi la description de la cohomologie ci-dessus montre que $\gr_{d-d_\rho}^W(H^{d-1}(S_{\Gamma,n}^{ca},\o\QM_l)[\o\rho])\neq 0$, et donc la conjecture monodromie-poids dans sa version \ref{equiMP} implique que $N_{\rho,\Gamma,d-1}^{d_\rho-1}\neq 0$. On a donc $iii)\Rightarrow ii)$. 

Il nous reste \`a prouver que $i)\Rightarrow iii)$. En fait le seul espace de cohomologie qui peut poser probl\`eme pour la conjecture monodromie-poids est celui de degr\'e m\'edian $j=d-1$, les autres \'etant purs. Il nous faut alors expliciter l'action de $N_{\rho,\Gamma,d-1}$ sur 
$H^{d-1}(S_{\Gamma,n}^{ca},\o\QM_l)[\o\rho]$. Mais si on suppose la propri\'et\'e $i)$, alors dans l'isomorphisme
$$H^{d-1}(S_{\Gamma,n}^{ca},\o\QM_l)[\o\rho] \simeq 
\bigoplus_{i=0}^{d_\rho-1}
\ext{i}{\pi^{\leq i}_\rho}{(\pi_{\rho}^\emptyset)^{m_{\rho,\Gamma}^\emptyset}\oplus (\pi_{\rho}^{S_{\rho}})^{m'_{\rho,\Gamma}}}{\o{G}}^* \otimes \sigma'(\rho^\vee)|-|^{-i},
$$
induit par le scindage $\alpha_\phi$, l'op\'erateur $N_\rho$ agit par $\cup$-produit et la description \ref{theoextell} ii) de ce $\cup$-produit montre que $N_\rho$ induit des isomorphismes
$$ \ext{i}{\pi^{\leq i}_\rho}{(\pi_{\rho}^\emptyset)^{m_{\rho,\Gamma}^\emptyset}}{\o{G}}^* \otimes \sigma'(\rho^\vee)|-|^{-i} \simto
\ext{i-1}{\pi_{\rho}^{<i}}{(\pi_{\rho}^\emptyset)^{m_{\rho,\Gamma}^\emptyset}}{\o{G}}^* \otimes\sigma'(\rho^\vee)|-|^{-i+1}
$$ pour tout $i\in \{1,\cdots, d_\rho-1\}$.
Or, toujours par la description de la cohomologie donn\'ee plus haut, on a pour 
%tout $k\in \ZM\setminus \{0\}$
$i\neq \frac{d_\rho-1}{2}$  
$$ \ext{i}{\pi^{\leq i}_\rho}{(\pi_{\rho}^\emptyset)^{m_{\rho,\Gamma}^\emptyset}}{\o{G}}^* \otimes \sigma'(\rho^\vee)|-|^{-i} \simto \gr_{(d-1)+(1-d_\rho+2i)}^W(H^{d-1}(S_{\Gamma,n}^{ca},\o\QM_l)[\o\rho] ) $$
et pour $i=\frac{d_\rho-1}{2}$ (lorsque $d_\rho$ est impair !) on a
$$ \ext{i}{\pi^{\leq i}_\rho}{(\pi_{\rho}^\emptyset)^{m_{\rho,\Gamma}^\emptyset} \oplus (\pi_{\rho}^{S_{\rho}})^{m'_{\rho,\Gamma}}}{\o{G}}^* \otimes \sigma'(\rho^\vee)|-|^{-i} \simto \gr_{d-1}^W(H^{d-1}(S_{\Gamma,n}^{ca},\o\QM_l)[\o\rho] ). $$
Dans tous les cas, la description de l'action de $N_\rho$ par $\cup$-produit montre que pour tout $k\in\ZM$, $N^k$ induit un isomorphisme
$$ N^k :\;\; \gr^W_{d-1+k}\left(H^{d-1}(S_{\Gamma,n}^{ca},\o\QM_l)[\o\rho]\right) \simto \gr^W_{d-1-k}\left(H^{d-1}(S_{\Gamma,n}^{ca},\o\QM_l)[\o\rho]\right) $$
et par \ref{equiMP}, la conjecture monodromie-poids pour la partie $\rho$-covariante $H^{d-1}(S_{\Gamma,n}^{ca},\o\QM_l)[\o\rho]$ en d\'ecoule.
\end{proof}
 
On en d\'eduit le deuxi\`eme r\'esultat principal de l'introduction, le th\'eor\`eme \ref{main2}.

\appendix

\def\ob{\hbox{ Ob}}

\section{Cat\'egories ab\'eliennes admissibles et topos fibr\'es}
\setcounter{subsubsection}{0}
%	\subsection{Topos fibr\'es}
Soit $I$ une petite cat\'egorie et $X\To{} I$ un topos fibr\'e sur $I$. Nous suivrons autant que possible les notations de \cite[ch. VI]{IllCC2}. En particulier, le topos total de $X$ sera not\'e $\top(X)$, le symbole $X^{dis}$ d\'esignera le topos $\bigsqcup_{i\in I} X_i$, et la lettre $e : X^{dis} \To{} \top(X)$ le morphisme de topos \'evident, d\'ecrit en \cite[6.1.1]{IllCC2}. On dispose donc d'une suite de trois foncteurs adjoints $(e_!,e^*,e_*)$ reliant les cat\'egories  $X^{dis}$ et $\top(X)$.

Soit $A$ un (pro)-anneau de $\top(X)$ et $A^{dis}:= e^*(A)$. Les foncteurs $e^*$ et $e_*$ respectent les cat\'egories de modules respectives et y restent adjoints, et $e^*$ admet encore un adjoint \`a gauche pour les modules, que nous noterons $e_!^A$ pour le distinguer de $e_!$. Nous supposerons par la suite que $e_!^A$ est {\em exact}. Il revient au m\^eme de demander que ``les" morphismes de topos annel\'es $(X_{s(\alpha)},A_{s(\alpha)}) \To{} (X_{b(\alpha)},A_{b(\alpha)})$ associ\'es aux fl\`eches $\alpha$ dans $I$ (de source $s(\alpha)$ et but $b(\alpha)$) soient {\em plats}, ce qui est par exemple v\'erifi\'e si $A$ est {\em constant}.

On notera simplement $D^{top}:=D^+_A(\top(X))$ et $D^{dis}:= D^+_{A^{dis}}(X^{dis})$. 
%, tel que le morphisme canonique $e_!e^*A\To{} A$ soit plat. C'est le cas lorsque $A$ est constant. 
On a donc une paire de foncteurs adjoints 
$ (e^*,R{e_*})$ reliant $D^{top}$ et $D^{dis}$.
%les cat\'egories d\'eriv\'ees  $D^+_A(X^{dis})$ et $D^+_A(\top(X))$.
On d\'efinit une cat\'egorie  $D^{naif}$  dont 

\begin{itemize} 
\item les objets sont les paires $(K,\kappa)$ o\`u $K\in D^{dis}$ et $\kappa: K\To{} e^*Re_* K$ est un morphisme tel que
\begin{enumerate}
 \item la compos\'ee $K\To{\kappa} e^*Re_* K \To{Adj} K$ est l'identit\'e, et
 \item les  deux compos\'ees 
$\xymatrix{ K \ar[r]^-\kappa  %\ar[r]\ar@<.8ex>[r]\ar@<-.8ex>[r] &
 &  e^*Re_* K 
  \ar@<.5ex>[r]^-{^{ e^*Re_*\kappa}} \ar@<-.5ex>[r]_-{_{e^*(Adj)Re_*}} & e^*Re_*e^*Re_* K 
  }$ sont \'egales.
 
\end{enumerate}
\item les fl\`eches $(K,\kappa)\To{} (K',\kappa')$ sont les morphismes $K\To{\alpha} K'$ tels que $\kappa'\circ \alpha= e^*Re_*(\alpha)\circ \kappa$. On a donc une suite exacte :
\ini\begin{equation} \label{sex}
%\xymatrix{ 0\ar[r] & \hom{(K,\kappa)}{(K',\kappa')}{D^{naif}}
% \ar[r]  %\ar[r]\ar@<.8ex>[r]\ar@<-.8ex>[r] &
% &   \hom{K}{K'}{D^{dis}}  
%  \ar@<.5ex>[r]^{\delta_1} \ar@<-.5ex>[r]_{\delta_2} &  \hom{K}{e^*Re_* K'}{D^{dis}}\} 
%  }
 0\To{}  \hom{(K,\kappa)}{(K',\kappa')}{D^{naif}}
 \To{}   \hom{K}{K'}{D^{dis}}  
  \To{\delta}   \hom{K}{e^*Re_* K'}{D^{dis}}
  \end{equation}
  o\`u $\delta$ d\'esigne la diff\'erence $\kappa'\circ \alpha-e^*Re_*(\alpha)\circ \kappa$.
\end{itemize}
La cat\'egorie $D^{naif}$ est $\ZM$-gradu\'ee. On a un foncteur $e^*_{naif} : D^{naif} \To{} D^{dis}$ d'oubli du morphisme $\kappa$ et on dira qu'un triangle de $D^{naif}$ est distingu\'e si son image par $e^*_{naif}$ l'est. La cat\'egorie $D^{naif}$ {\em n'est pas triangul\'ee}. Elle s'identifie \`a une sous-cat\'egorie pleine des sections de la cat\'egorie bi-fibr\'ee sur $I$ dont les fibres sont les $D^+_{A_i}(X_i)$.

On a aussi un foncteur   $\omega: \; D^{top} \To{} D^{naif}$ qui envoie un objet $K$ sur la paire $(e^*K,  e^*K \To{e^*(Adj) } e^*Re_*e^* K)$, et qui permet de factoriser  $e^*=e^*_{naif}\circ \omega$.

%\subsection{Cat\'egories ab\'eliennes admissibles}
Soit $C^{dis}$ une sous-cat\'egorie ab\'elienne admissible de $D^{dis}$. Rappelons \cite[1.2.5]{BBD} que cela signifie que (i) $\hom{K}{K'[n]}{D^{dis}}=0$ pour tout $n<0$ et tous $K,K'\in C^{dis}$, et (ii) les suites exactes courtes de $C^{dis}$ se d\'eduisent des triangles distingu\'es par oubli de la fl\`eche de bord.
Nous noterons $C^{naif}$, resp. $C^{top}$, la sous-cat\'egorie pleine de $D^{naif}$, resp. de $D^{top}$ form\'ee des objets $X$ tels que $e_{naif}^*(X)$, resp. $e^*(X)$,  soit isomorphe \`a un objet de  $C^{dis}$. Le foncteur $\omega$ se restreint donc en un foncteur $C^{top} \To{} C^{naif}$.

\begin{lemme}
Supposons $C^{dis}$ stable par $e^*Re_*$. Alors $C^{naif}$ est une cat\'egorie ab\'elienne.
\end{lemme}
\begin{proof}
On d\'efinit les noyaux et conoyaux de la mani\`ere la plus na\"ive qui soit, sachant 
que sous l'hypoth\`ese, $e^*Re_*$ induit un endo-foncteur {\em exact} de $C^{dis}$. Nous laissons la v\'erification des axiomes $\hbox{Ab}_i$ au lecteur.
\end{proof}

Remarquons que l'hypoth\`ese ``stable par $e^*Re_*$" revient \`a demander que pour toute fl\`eche $\alpha$ de $I$, ``le" foncteur  $D^+_{A_{s(\alpha)}}(X_{s(\alpha)})\To{R\alpha_*}  D^+_{A_{b(\alpha)}}(X_{b(\alpha)})$ envoie $C^{dis}_{s(\alpha)}$ and $C^{dis}_{b(\alpha)}$.
Cette hypoth\`ese n'est donc g\'en\'eralement pas v\'erifi\'ee par l'exemple le plus simple de cat\'egorie $C^{dis}$, \`a savoir $\Mo{A^{dis}}{X^{dis}}$. Pourtant dans ce cas encore la cat\'egorie $C^{naif}$ est bien-s\^ur ab\'elienne, et plus g\'en\'eralement, la conclusion du lemme reste vraie si on suppose que $C^{dis}$ stable par $e^* e_!^A$.

\begin{prop} \label{fptopfib}
Soit $C^{dis}$ une sous-cat\'egorie ab\'elienne admissible et stable par $e^*Re_*$ de $D^{dis}$, alors le foncteur $\omega: C^{top} \To{} C^{naif}$ est une \'equivalence de cat\'egories. En particulier, $C^{top}$ est une sous-cat\'egorie ab\'elienne admissible de $D^{top}$.
\end{prop}

\def\sb{\hbox{\bf s}}

\noindent{ \it Preuve de la pleine fid\'elit\'e} : fixons pour cela  deux objets $K,L$ dans $C^{top}$; ce sont donc des complexes de $A$-modules dans $\top(X)$. Choisissons un complexe $I_L^\bullet$ \`a composantes injectives et quasi-isomorphe \`a $Y$. On a donc $\hom{K}{L}{C^{top}}=H^0(\sb(\hom{X}{I_L^\bullet}{C_A(X)}))$ o\`u $C_A(X)$ est la cat\'egorie des complexes de $A$-modules dans $\top(X)$ et $\sb$ d\'esigne le complexe simple associ\'e \`a un complexe double. Soit 
$$ I_L^\bullet \To{} e_*e^*(I_L^\bullet) \To{} (e_*e^*)^2(I_L^\bullet) \To{} \cdots $$
la r\'esolution standard de $I_L^\bullet$ dans $C_A(X)$ associ\'ee \`a la paire adjointe $(e^*,e_*)$. Nous avons suppos\'e que $e_!^A$ est exact, ce qui implique que  $e^*$ envoie   injectifs sur injectifs, tout comme $e_*$. On en d\'eduit une suite spectrale :
$$E_1^{p,q} = H^q(\sb(\hom{K}{(e_*e^*)^{p+1} I_L^\bullet}{C_A(X)})) \Rightarrow H^{p+q}(\sb(\hom{K}{I_L^\bullet}{C_A(X)})), $$
autrement dit,  une suite spectrale 
$$ E_1^{p,q} = \ext{q}{K}{(e_*Re^*)^{p+1} L}{D^{top}} \Rightarrow \ext{p+q}{K}{L}{D^{top}} .$$
Par d\'efinition de cette suite spectrale, on a $E_1^{pq}=0$ si $p<0$. Par l'hypoth\`ese $C^{dis}$ admissible et stable par $e^*Re_*$, on a pour tout $q<0$ et tout $p\geq 0$
$$ E_1^{pq} = \ext{q}{e^*K}{e^*(e_*Re^*)^p L}{D^{dis}} = 0.$$
On en d\'eduit sur la ligne $p=0$ de la suite spectrale une suite exacte :
$$ 0 \mapsto \hom{K}{L}{C^{top}} \To{} \hom{e^*K}{e^*L}{D^{dis}} \To{\delta} \hom{e^*K}{e^*Re_*e^*L}{D^{dis}}. $$
Revenant \`a la d\'efinition de $\omega$, on constate que $\delta$ s'identifie (au signe pr\`es) \`a la fl\`eche not\'ee aussi $\delta$ dans la  suite exacte \ref{sex} appliqu\'ee \`a $\omega K$ et $\omega L$.
Il s'en suit que l'application $\hom{K}{L}{C^{top}} \To{} \hom{\omega K}{\omega L}{C^{naif}}$ est bijective.

\alin{Essentielle surjectivit\'e} On se donne  un objet $(K,\kappa)$ de $C^{naif}$ et on choisit un complexe $I_K^\bullet$ de $A^{dis}$-modules quasi-isomorphe \`a $K$ et \`a composantes injectives,  ainsi qu'un rel\`evement  $I_K^\bullet\To{} e^*e_* I_K^\bullet$ de $\kappa$ en un morphisme de complexes, encore not\'e $\kappa$. %Pour tout entier $p\geq 0$, posons $J^{p\bullet}:= e_*(e^*e_*)^p I_K^\bullet$.

On d\'efinit un syst\`eme de morphismes de la forme suivante :
\ini\begin{equation}\label{comp}
\xymatrix{ e_* I_K^\bullet 
  \ar@<.5ex>[r]^{e_*(\kappa)} \ar@<-.5ex>[r]_{Adj_{e_*I_K}} & e_*(e^*e_*) I_K^\bullet \ar[r] \ar@<.8ex>[r]^{e_*e^*e_*(\kappa)} \ar@<-.8ex>[r]  & e_*(e^*e_*)^2 I_K^\bullet \cdots
  }
\end{equation}
o\`u chaque fl\`eche sup\'erieure se d\'eduit de $\kappa$ et les autres fl\`eches sont de la forme
\begin{eqnarray*}
  (e_*e^*)^{i-1} \hbox{Adj}_{(e_*e^*)^{p-i}e_*I_K} & : & 
  (e_*e^*)^{i-1} (e_*e^*)^{p-i} e_* I_K^\bullet \To{}  (e_*e^*)^{i-1}(e_*e^*)(e_*e^*)^{p-i} e_* I_K^\bullet
 \end{eqnarray*}
pour tout $0<i\leq p$ entiers.
%et pour tout $p>0$ et $i=0$, on pose
%$$ d_p^0 := e_*(e^*e_*)^{p-1}(\kappa) :\; J^{p-1,\bullet} \To{} J^{p,\bullet} .$$
%On a donc un syst\`eme de 
Les axiomes impos\'es \`a $\kappa$ nous disent que, {\em dans la cat\'egorie homotopique} $K^{top}$ des complexes de $A^{top}$-modules, ce syst\`eme se prolonge en un objet cosimplicial, {\em i.e.} un foncteur de la cat\'egorie des ensembles finis ordonn\'es non-vides vers $K^{top}$. Si on pouvait remonter cet objet cosimplicial \`a la cat\'egorie (ordinaire) des complexes de $A^{top}$-modules, on montrerait facilement que le complexe de cochaines associ\'e est un rel\`evement cherch\'e de $(K,\kappa)$ dans $D^{top}$.
Mais ceci n'est g\'en\'eralement pas faisable, et il nous faut utiliser un substitut remarquable introduit par les auteurs de \cite[3.2]{BBD}.

\def\fl{\hbox{\sl Fl}} 

\def\ve{\varepsilon}

\alin{Complexes homotopiquement simpliciaux de \cite[3.2]{BBD}\footnote{Cette terminologie n'est pas dans {\em loc. cit.} mais l'auteur ignore s'il en existe une standard.}} 
Notons $\o\Delta$ la cat\'egorie dont les objets sont les entiers $\geq -1$ et les morphismes sont donn\'es par $\hom{p}{p'}{\o\Delta}:= \{\hbox{applications injectives croissantes }\alpha:\;[0,p]\To{}[0,p']\},$ 
avec la convention que $[0,-1]=\emptyset$ et que $\hom{-1}{p}{\o\Delta}$ est un singleton dont nous noterons $\ve_p$ l'unique \'el\'ement. La source et le but d'une fl\`eche dans $\hbox{Fl}(\o\Delta)$ seront not\'es $s(\alpha)$ et $b(\alpha)$ et la diff\'erence $b(\alpha)-s(\alpha)$ sera not\'ee $|\alpha|$. On note aussi  $\Delta$ la sous-cat\'egorie pleine des entiers $\geq 0$.

Soit $\AC$ une cat\'egorie ab\'elienne. Nous appellerons {\em complexe homotopiquement simplicial} de $\AC$ la donn\'ee 
 d'une suite $(J^{p,\bullet})_{p\in\NM}$  d'objets $\ZM$-gradu\'es de $\AC$ munie  d'une famille $(d(\alpha))_{\alpha\in\fl(\Delta)}$ de morphismes $J^{s(\alpha),\bullet}\To{} J^{b(\alpha),\bullet}[1-|\alpha|]$ satisfaisant la propri\'et\'e 
$$\forall \alpha\in \fl(\Delta),\;\; \sum_{\alpha=\beta\gamma}d(\beta)d(\gamma) =0. $$
Rappelons que cette propri\'et\'e implique que les $d(\id_p)$ sont des diff\'erentielles, et que la famille des  $d(\alpha)$ pour $|\alpha|=1$ est une famille de morphismes de complexes qui, dans la cat\'egorie homotopique, se prolonge en un complexe cosimplicial strict. R\'eciproquement, les auteurs de \cite{BBD} montrent comment {\em sous certaines conditions} un complexe cosimplicial strict de la cat\'egorie homotopique $K(\AC)$ peut se relever en un complexe homotopiquement simplicial. Ces conditions sont v\'erifi\'ees par notre syst\`eme \ref{comp} mais nous aurons besoin d'un rel\`evement assez explicite, {\em cf} \ref{relevexpl} ci-dessous.

\`A tout complexe homotopiquement simplicial  $(J^{*\bullet},(d(\alpha))_{\alpha})$  uniform\'ement born\'e inf\'erieurement en $\bullet$ ({\em i.e} $J^{*q}=0$ pour $q<<0$), les auteurs de \cite{BBD}  associent le complexe simple 
$ (\sb(J)^\bullet, d^\bullet)$ d\'efini par
$$ \sb(J)^n:= \bigoplus_{p+q=n} J^{pq},\;\;\hbox{ et }\;\; d^n := \bigoplus_{p+q=n} \sum_{s(\alpha)=p} d(\alpha)_{|J^{p,q}}, $$
ces sommes \'etant finies. 

Nous aurons besoin de versions ``augment\'ees'' de ces objet. Nous appellerons donc {\em complexe homotopiquement simplicial augment\'e} de $\AC$ la donn\'ee 
 d'une suite $(J^{p,\bullet})_{p\geq -1}$ munie d'une famille $(d(\alpha))_{\alpha\in\fl(\o\Delta)}$ de morphismes satisfaisant les m\^emes propri\'et\'es formelles que ci-dessus, avec $\o\Delta$ \`a la place de $\Delta$. On lui associe
aussi un complexe de cochaines $\sb(J)^\bullet$ par la m\^eme formule que ci-dessus. En notant $J_{|\Delta}$ le complexe homotopiquement simplicial sous-jacent, on peut d\'efinir une augmentation 
 $$ \gamma:=\sum_{p\geq 0} d(\ve_p) :\,\, J^{-1,\bullet} \To{} \sb(J_{|\Delta})^\bullet,$$
\`a condition de prendre sur $J^{-1,\bullet}$ l'oppos\'e de la diff\'erentielle $d(\id_{-1})$.
On a alors un triangle de complexes dans $\AC$ 
$$ J^{-1,\bullet} \To{\gamma} \sb(J_{|\Delta})^\bullet \injo \sb(J)^\bullet \twoheadrightarrow J^{-1,\bullet}[1] $$
dont l'image dans la cat\'egorie homotopique est un triangle distingu\'e.
% note encore $\sb(J)^\bullet$ le complexe de cochaines associ\'e comme ci-dessus  au complexe homotopiquement simplicial {\em non-augment\'e} sous-jacent. 
%Nous laisserons le lecteur v\'erifier que les axiomes munissent ce dernier d'une 

\alin{Le lemme crucial}
Avant de continuer, fixons quelques notations d'alg\`ebre simpliciale : 
%\begin{nota} 
\begin{enumerate}
\item On note $\partial : \o{\Delta}\To{} \Delta$ le foncteur d\'ecalage d\'efini sur les objets par $\partial(p):=p+1$ et sur les fl\`eches par $\partial(\alpha)(i):= \left\{\begin{array}{ll} 0 & \hbox{si } i=0 \\ 
\alpha(i-1)+1 & \hbox{si } 0<i\leq s(\alpha)+1 \end{array}\right.$

\item On note $\sigma_p \in \hom{p-1}{p}{\o\Delta}$ l'application d\'efinie par $i\mapsto i+1$ lorsque $p>0$ et par $\sigma_0 :=\ve_0$ pour $p=0$. Ainsi l'"op\'erateur de face" $[0,p-1]\To{} [0,p]$ qui saute l'entier $i\in [0,p]$ est donn\'e par la formule $\partial^i\sigma_{p-i}$.
%les $p+1$ fl\`eches $p-1\To{} p$ dans $\hbox{FL}(\Delta)$ (appel\'es op\'erateurs de faces) sont donn\'ees par 
\end{enumerate}
%\end{nota}

\begin{lemme} \label{relevexpl}
Il existe une famille $d(\alpha)_{\alpha\in\fl\o\Delta}$ de morphismes de $A^{dis}$-modules gradu\'es
$$ d(\alpha):\; (e^*e_*)^{s(\alpha)+1}I_K^\bullet \To{} (e^*e_*)^{b(\alpha)+1}I_K^\bullet [1-|\alpha|] $$
v\'erifiant les propri\'et\'es suivantes :
\begin{enumerate}
        \item $\forall \alpha\in \fl(\o\Delta),\;\; d(\partial \alpha) = -(e^*e_*)(d(\alpha)).$
        \item $d(\id_{-1})$ est la diff\'erentielle du complexe $I_K^\bullet$.
        \item $d(\sigma_0)=\kappa$ et pour tout $p>0$, $d(\sigma_p)=e^*\hbox{Adj}_{(e_*e^*)^{p-1}e_*I_K^\bullet}$ o\`u $\hbox{Adj}_X : X \To{} e_*e^*X $ est le morphisme d'adjonction.
        \item Si $|\alpha|\geq 2$, alors $d(\alpha)\neq 0 \Rightarrow \alpha = \partial^{s(\alpha)+1}(\ve_{|\alpha|-1})$.
\item $\forall \alpha\in \fl(\o\Delta),\;\; S(\alpha):=\sum_{\alpha=\beta\gamma} d(\beta)d(\gamma)=0 $.
\end{enumerate}
\end{lemme}
%Remarquons que si une telle famille existe, elle est enti\`erement d\'etermin\'ee par la famille des $d(\ve_p)$ pour $p\geq 0$, en vertu des axiomes i) et iv), car $\alpha\in \im\partial^{s(\alpha)+1} \equ \alpha= \partial^{s(\alpha) +1}(\ve_{|\alpha|-1})$.

\begin{proof} Les propri\'et\'es i), ii) et iii) imposent tous les $d(\alpha)$ pour $|\alpha|\leq 1$. La propri\'et\'e contraignante est bien-s\^ur v). Lorsque $|\alpha|=0$, v) demande simplement que $d(\alpha)^2=0$, ce qui est bien le cas. Lorsque $|\alpha|=1$, v) demande que $d(\alpha)$ soit un morphisme de complexes (au signe pr\`es), ce qui est encore le cas. 

Nous prouvons maintenant l'existence des $d(\alpha)$ pour $|\alpha|\geq 2$ v\'erifiant les propri\'et\'es i), iv) et v) par r\'ecurrence sur $|\alpha|$. 
Supposons donc construites les $d(\beta)$ pour $|\beta|<|\alpha|$. Deux cas se pr\'esentent.

Si $\alpha\neq \partial^{s(\alpha)+1}(\ve_{|\alpha|-1})$, alors la propr\'et\'e iv) impose $d(\alpha)=0$ et pour satisfaire v) il nous faut donc v\'erifier que la somme $S'(\alpha):= \sum_{\beta\gamma=\alpha ; \beta,\gamma\neq \alpha}d(\beta)d(\gamma)$ est nulle. Par la propri\'et\'e i) de l'hypoth\`ese de r\'ecurrence, et le fait que le foncteur $\partial$ induit une bijection
$$ (\partial,\partial):\; \{(\beta,\gamma)\; \beta\gamma=\alpha\} \simto \{(\beta,\gamma)\; \beta\gamma=\partial\alpha\},$$
on peut supposer que $\alpha\notin \im(\partial)$. Par la propri\'et\'e iv) de l'hypoth\`ese de r\'ecurrence, la somme $S'(\alpha)$ n'a alors que deux termes non-nuls :
$$ S'(\alpha) = d(\sigma_{b(\alpha)}) d(\partial^{s(\alpha)+1} \ve_{|\alpha|-2}) + d(\partial^{s(\alpha)+2}\ve_{|\alpha|-2})d(\sigma_{s(\alpha)+1}) $$
et sa nullit\'e r\'esulte de la fonctorialit\'e de $\hbox{Adj}_X$ en $X$.

 Si $\alpha = \partial^{s(\alpha)+1}(\ve_{|\alpha|-1})$, alors comme ci-dessus, par la propri\'et\'e i) de l'hypoth\`ese de r\'ecurrence, il suffit de trouver $d(\ve_{|\alpha|-1})$ tel que $S(\ve_{|\alpha|-1})=0$. Or, comme dans \cite[p. 90]{BBD}, l'hypoth\`ese de r\'ecurrence implique $S'(\ve_{|\alpha|-1})d(\id_{-1}) = d(\id_{|\alpha|-1}) S'(\ve_{|\alpha|-1})$, autrement dit, $S'(\ve_{|\alpha|-1})$ est un morphisme de complexes $I_K^\bullet \To{} (e^*e_*)^{|\alpha|}(I_K^\bullet)[2-|\alpha|]$. Deux cas se pr\'esentent  \`a nouveau :
 \begin{itemize}
 \item si $|\alpha|=2$, on invoque la propri\'et\'e ii) de $\kappa$ qui assure que $S'(\ve_1)$ est nulle dans la cat\'egorie homotopique. 
 \item si $|\alpha|>2$, alors on invoque la propri\'et\'e ``$\CC^{dis}$ admissible et stable par $e^*Re_*$" qui implique $\hom{I_K^\bullet}{(e^*e_*)^{|\alpha|}I_K^\bullet[2-|\alpha|]}{D^{dis}}=0$ et donc que $S'(\ve_{|\alpha|-1})$ est aussi nulle dans la cat\'egorie homotopique.
 \end{itemize}
Dans chacun des cas, il ne reste plus qu'\`a choisir pour $d(\ve_{|\alpha|-1})$ une homotopie entre $S'(\ve_{|\alpha|-1})$ et $0$.
\end{proof}

\begin{lemme}
Avec les notations du lemme pr\'ec\'edent, d\'efinissons une famille $c(\alpha)_{\alpha\in \fl\Delta}$ de morphismes de $A^{top}$-modules gradu\'es 
$$ c(\alpha):\; e_*(e^*e_*)^{s(\alpha)}I_K^\bullet \To{} e_*(e^*e_*)^{b(\alpha)}I_K^\bullet [1-|\alpha|] $$
par les r\`egles suivantes :
\begin{enumerate}
        \item Si $\alpha\in\im(\partial)$, alors $c(\alpha):=-e_*(d(\partial^{-1}\alpha))$.
        \item Si $\alpha\notin\im(\partial)$, alors 
        \begin{itemize}
        \item $c(\alpha):=\hbox{Adj}_{(e_*e^*)^{s(\alpha)}e_*(I_K^\bullet)}$ si $|\alpha|=1$
        \item $c(\alpha)=0$ sinon.
\end{itemize}
\end{enumerate}
 Alors, pour tout $\alpha\in\fl(\Delta)$, on a $\sum_{\beta\gamma=\alpha} c(\beta)c(\gamma)=0$.
\end{lemme}
\begin{proof}
On remarque que pour toute fl\`eche $\alpha$ de $\Delta$, on a $e^* c(\alpha) = d(\alpha)$. Or, $e^*$ est fid\`ele sur les $A^{top}$-modules.
\end{proof}

Le syst\`eme des  $J^{p,\bullet}_{top}:=e_*(e^*e_*)^pI_K^\bullet$,  $p\geq 0$ muni de la famille des $c(\alpha)_{\alpha\in\fl\Delta}$ du lemme ci-dessus est un complexe homotopiquement simplicial de $A^{top}$-modules, tandis que le syst\`eme des $J^{p,\bullet}_{dis}:=(e^*e_*)^pI_K^\bullet$, $p\geq -1$ muni des $(d(\alpha))_{\alpha\in\fl\o\Delta}$ du lemme  \ref{relevexpl} est un complexe homotopiquement simplicial augment\'e de $A^{dis}$-modules. On a par construction ${J_{dis}}_{|\Delta}=e^*(J_{top})$ ;
on a donc un morphisme de complexes de $A^{dis}$-modules $I_K^\bullet \To{\gamma} e^*\sb(J_{top})^\bullet$. 
Le lemme suivant montre que le complexe $\sb(J_{top})^\bullet $ de $D^{top}$ rel\`eve l'objet $(K,\kappa)$ de $D^{naif}$ et r\'esoud donc la question de l'essentielle surjectivit\'e.
\begin{lemme} 
\begin{enumerate}
        \item $\gamma:\;I_K^\bullet \To{} e^*\sb(J_{top})^\bullet$ est un quasi-isomorphisme.
        \item le diagramme suivant est commutatif dans la cat\'egorie d\'eriv\'ee $D^{dis}$ 
        $$\xymatrix{ I_K^\bullet \ar[r]^-\gamma \ar[d]_\kappa & e^*\sb(J_{top})^\bullet \ar[d]^{e^*\hbox{Adj}_{\sb(J_{top})^\bullet}} \\
           e^*e_* I_K^\bullet \ar[r]_-{e^*e_*(\gamma)} & e^*e_*e^*\sb(J_{top})^\bullet}$$
\end{enumerate}
\end{lemme}
\begin{proof}
Nous prouverons d'abord ii). Pour cela, nous commen{\c c}ons par quelques g\'en\'eralit\'es sur les complexes homotopiquement simpliciaux augment\'es (d'une cat\'egorie ab\'elienne $\AC$ quelconque). Un morphisme $f$ entre deux tels objets $(J^{*\bullet}_i, d_i(\alpha)_\alpha)$ pour $i=1,2 $ consiste en une famille $f=(f(\alpha))_{\alpha\in \o\fl(\Delta)}$ de morphismes $f(\alpha):\; J_1^{s(\alpha),\bullet} \To{} J_2^{b(\alpha),\bullet}[-|\alpha|]$ d'objets gradu\'es de $\AC$ v\'erifiant la propri\'et\'e :
$$ \forall \alpha\in \fl(\o\Delta),\;\; \sum_{\alpha=\beta\gamma} f(\beta)d_1(\gamma) = \sum_{\alpha=\beta\gamma} d_2(\beta)f(\gamma). $$
On v\'erifie sans peine qu'un tel syst\`eme induit un morphisme de complexes $\sb(J_1)^\bullet \To{\sb(f)} \sb(J_2)^\bullet$ et que dans le morphisme de triangles :
\ini\begin{equation} \label{mortri}
\xymatrix{ J_1^{-1,\bullet} \ar[r]^{\gamma_1} \ar[d]_{f(\id_{-1})} & \sb({J_1}_{|\Delta})^\bullet \ar[r] \ar[d]_{\sb(f_{|\Delta})} & \sb(J_1)^\bullet \ar[r] \ar[d]_{\sb(f)} & J^{-1,\bullet}_1[1] \ar[d]^{f(\id_{-1})[1]} \\
 J_2^{-1,\bullet} \ar[r]^{\gamma_1} & \sb({J_2}_{|\Delta})^\bullet \ar[r] & \sb(J_2)^\bullet \ar[r] & J^{-1,\bullet}_2[1] },
 \end{equation}
les deux carr\'es de droite sont commutatifs, et donc le premier est commutatif {\em dans la cat\'egorie homotopique} (mais g\'en\'eralement pas dans la cat\'egorie des complexes).

Ceci \'etant, on d\'efinit le {\em d\'ecal\'e} d'un complexe homotopiquement simplicial augment\'e $(J^{*\bullet},d(\alpha)_\alpha)$ par les formules :
$$ (\partial J)^{p,\bullet}:= J^{p+1,\bullet},\;\;\hbox{ et }\;\; \forall \alpha\in \fl(\o\Delta),\; \partial d(\alpha) := -d(\partial \alpha). $$
On d\'efinit aussi un morphisme $(J,d) \To{f^\partial} (\partial J,\partial d)$ par 
$$ \forall \alpha\in \fl(\o\Delta),\;\; f^\partial(\alpha):= d(\sigma_{b(\alpha)+1}\circ\alpha) .$$
(Pour v\'erifier que ce syst\`eme est bien un morphisme, on utilise l'identit\'e
$$ \sum_{\alpha=\beta\gamma} d(\sigma_{b(\beta)+1}\beta)d(\gamma) + \sum_{\alpha=\beta\gamma} d(\partial \beta)d(\sigma_{b(\gamma)+1}\gamma) = \sum_{\sigma_{b(\alpha)+1}\alpha= \beta'\gamma'} d(\beta')d(\gamma') $$
qui repose sur le fait qu'une fl\`eche $\delta$ se factorise sous $\sigma_{b(\delta)}$ \ssi\ elle n'est pas dans l'image de $\partial$.)

Appliquons ceci \`a $J_{dis}^{*\bullet}$. Par l'axiome i) du lemme \ref{relevexpl} on a $\partial J_{dis} = e^*e_*(J_{dis})$. Par l'axiome iv), les fl\`eches $f^\partial(\alpha)$ sont nulles d\`es que $|\alpha|>0$, et par l'axiome iii) on a $f^\partial(\id_p)=e^*\hbox{Adj}$ pour $p\geq 0$, tandis que $f^\partial(\id_{-1}) = \kappa$. La commutativit\'e \`a homotopie pr\`es du carr\'e du point ii) de l'\'enonc\'e vient donc de celle du premier carr\'e du diagramme \ref{mortri} appliqu\'e au morphisme $f^\partial$, compte tenu de l'\'egalit\'e ${J_{dis}}_{|\Delta}= e^*J_{top}$.

\bigskip

Passons \`a la preuve du point i). En vertu du triangle distingu\'e
$$ I_K^\bullet \To{\gamma} e^*\sb(J_{top})^\bullet \To{} \sb(J_{dis})^\bullet \To{} I_K^\bullet[1], $$
il suffit de montrer l'acyclicit\'e de $\sb(J_{dis})^\bullet$. Par l'axiome i) impos\'e \`a $\kappa$, la compos\'ee
$$ \sb(J_{dis})^\bullet \To{f^\partial} \sb(\partial J_{dis})^\bullet=e^*e_*\sb(J_{dis})^\bullet \To{Adj} \sb(J_{dis})^\bullet$$
est un isomorphisme dans la cat\'egorie d\'eriv\'ee $D^{dis}$. Il nous suffira donc de prouver l'acyclicit\'e de $e^*e_*\sb(J_{dis})^\bullet$. 

\`A ce point, il faut se rappeler que, dans la cat\'egorie homotopique le syst\`eme des  $J_{dis}^{p,\bullet}$, $p\geq -1$ muni des op\'erateurs de face $d(\partial^i\sigma_{p-i})$ se prolonge en un syst\`eme cosimplicial complet ({\em i.e.} avec op\'erateurs de d\'eg\'en\'erescence, lesquels sont donn\'es par l'axiome i) de $\kappa$). Comme l'objet cosimplicial $e^*e_*(J_{dis}^{*,\bullet})$ de $D^{dis}$ est le d\'ecal\'e de $J_{dis}^{*,\bullet}$, on sait qu'il est {\em homotopiquement trivial}, {\em cf} par exemple \cite[Prop VI.1.4]{IllCC2}.

Il nous reste plus qu'\`a invoquer le r\'esultat g\'en\'eral :
\begin{fact}
Soit $(J^{*\bullet},d(\alpha)_\alpha)$ un complexe homotopiquement simplicial augment\'e d'une cat\'egorie ab\'e\-lienne $\AC$ dont le complexe cosimplicial associ\'e dans $D(\AC)$ est homotopiquement trivial. Alors $\sb(J)^\bullet$ est acyclique. 
\end{fact}
\begin{proof}
Tout complexe homotopiquement simplicial augment\'e est muni d'une filtration 
$\tau_{\bullet \leq n}(J)^{*\bullet}$ d\'efinie par : 
$$ \tau_{\bullet \leq n}(J)^{pq} = \left\{\begin{array}{ll}
J^{pq} & \hbox{si } q<n \\
\ker(d(\id_p)_{|J^{pn}}) \hbox{si } q=n  \\
0 & \hbox{si } q>n 
\end{array}\right.
$$
\'equip\'e du syst\`eme des  restrictions des $(d(\alpha))_\alpha$.
On en d\'eduit une filtration croissante $T_n(\sb(J))^\bullet$ du complexe simple associ\'e.
On a aussi une filtration $\sigma_{\bullet \leq n}(J)^{*\bullet}$ d\'efinie par : 
$$ \sigma_{\bullet \leq n}(J)^{pq} = \left\{\begin{array}{ll}
J^{pq} & \hbox{si } q\leq n \\
\im(d(\id_p)_{|J^{pn}}) \hbox{si } q=n+1  \\
0 & \hbox{si } q>n+1 
\end{array}\right.
$$
\'equip\'e du syst\`eme des  restrictions des $(d(\alpha))_\alpha$, et 
dont on d\'eduit une filtration croissante $S_n(\sb(J))^\bullet$ du complexe simple associ\'e.
Il nous suffira de prouver que les gradu\'es $\gr^S_m\gr^T_n(\sb(J))$ sont acycliques.
Remarquons qu'ils sont nuls pour $m\neq n, n+1$.

Pour $m=n$, le gradu\'e est, au d\'ecalage de $n$ pr\`es, le complexe de cochaines augment\'e
$$ \HC^n(J^{-1,\bullet}) \To{\delta_0} \HC^n(J^{0,\bullet}) \To{} \cdots \To{} \HC^n(J^{p,\bullet}) \To{} \cdots $$
%dont la diff\'erentielle $\delta_p$ est 
associ\'e au complexe cosimplicial de $\AC$ d\'eduit du complexe cosimplicial $J^{*\bullet}$  de $D(\AC)$ par application du foncteur $\HC^n$. Il est donc homotopiquement trivial et {\em a fortiori} acyclique.

Pour $m=n-1$, le gradu\'e obtenu est le complexe simple associ\'e \`a un complexe homotopiquement simplicial augment\'e $(J^{*,\bullet },d(\alpha)_\alpha)$ dont les complexes $J^{p,\bullet}$ sont acycliques. Par la filtration d\'ecroissante b\^ete en l'indice $p$, on voit qu'un tel complexe est acyclique.
\end{proof}

\end{proof}

\end{document}